\documentclass[12pt]{article}

\usepackage{amsmath, amsthm, amssymb}

\usepackage{graphicx}
\def\R{{\mathbb R}}
\def\N{{\mathbb N}}
\theoremstyle{problem}
\newtheorem{problem}{Problem}[section]
\newtheorem{proposition}{Proposition}[section]
\newtheorem{theorem}{Theorem}[section]

\newtheorem{definition}{Definition}[section]
\newtheorem{lemma}{Lemma}[section]
\newtheorem{corollary}{Corollary}[section]

\def\vec#1{{\mathchoice{\mbox{\boldmath$\displaystyle#1$}}
{\mbox{\boldmath$\textstyle#1$}}
{\mbox{\boldmath$\scriptstyle#1$}}
{\mbox{\boldmath$\scriptscriptstyle#1$}}}}

\newcommand{\xib}{\vec{\xi}}
\def \gamb{\vec{\gamma}} 
\def \varrhob{\vec{\varrho}}
 \newcommand{\etab}{\vec{\eta}}
 \def\cA{{\cal A}}

 \newcommand{\fb}{\vec{f}}
 \def\ZR{{\mathbb R}}
 \newcommand{\Ebb}{{\bf E}} 
\newcommand{\Abb}{{\bf A}}
\def\Sym{\mathop{\rm Sym}\nolimits}

\begin{document}

\title{\sc{Existence of a weak solution to a nonlinear fluid-structure interaction problem modeling the flow
of an incompressible, viscous fluid in a cylinder with deformable walls}}
\author{ Boris Muha \thanks{Department of Mathematics,
    University of Houston, Houston, Texas 77204-3476,
    borism@math.hr} \and Sun\v{c}ica \v{C}ani\'{c}\thanks{Department of Mathematics,
    University of Houston, Houston, Texas 77204-3476,
    canic@math.uh.edu} }

\date{}

\maketitle

\begin{abstract}
We study a nonlinear, unsteady, moving boundary, fluid-structure interaction (FSI) problem arising in modeling blood flow through elastic and viscoelastic arteries. 
The fluid flow, which is driven by the time-dependent pressure data, 
is governed by 2D incompressible Navier-Stokes equations, while the elastodynamics of the cylindrical wall 
is modeled by  the 1D cylindrical Koiter  shell model. Two cases are considered: the linearly viscoelastic and the linearly elastic Koiter shell.
The fluid and structure are fully coupled (2-way coupling) via the kinematic and dynamic lateral boundary conditions
describing continuity of velocity (the no-slip condition), and balance of contact forces at the fluid-structure interface.
We prove existence
of weak solutions to the two FSI problems (the viscoelastic and the elastic case) as long as the cylinder radius is greater than zero.

The proof is based on a novel semi-discrete, operator splitting numerical scheme, 
known as the kinematically coupled scheme,
introduced in \cite{GioSun} to numerically solve the underlying FSI problems. 
The backbone of the kinematically coupled scheme is the well-known Marchuk-Yanenko scheme, 
also known as the Lie splitting scheme.
We effectively prove convergence of that numerical scheme to a solution of the corresponding FSI problem.  
\end{abstract}

\section{Introduction}
\begin{figure}[ht]
\centering{
\includegraphics[scale=0.35]{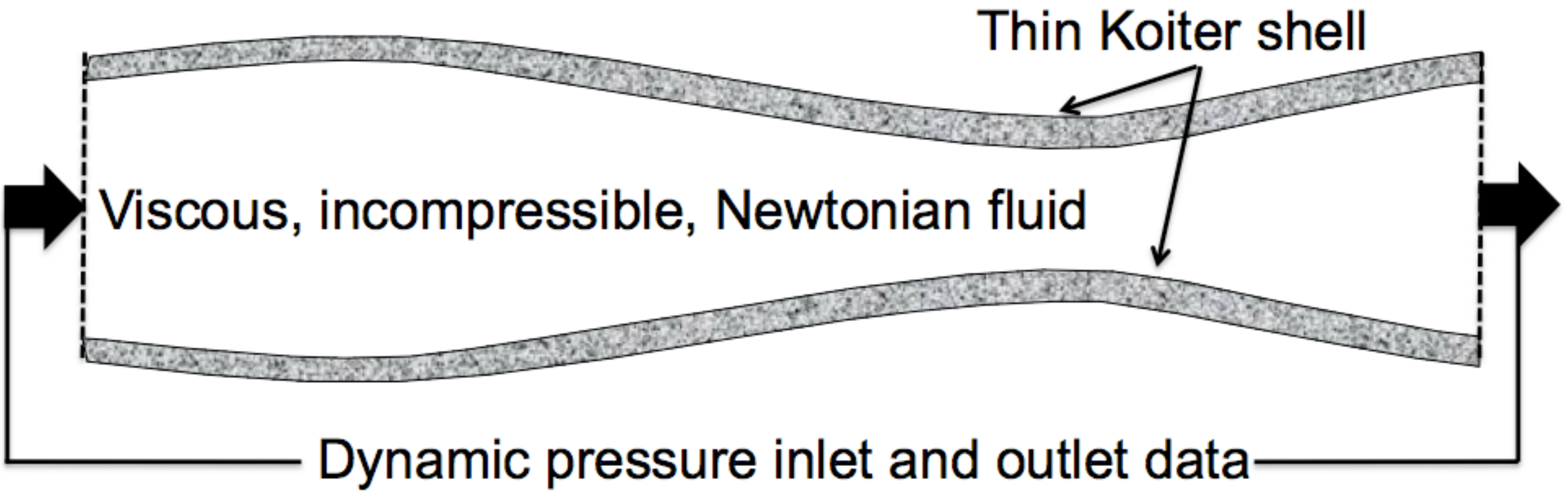}
}
\caption{Domain sketch.}
\label{fig:domain}
\end{figure}

We study the existence of a weak solution to a nonlinear moving boundary, unsteady, fluid-structure interaction (FSI) problem
between an incompressible, viscous, Newtonian fluid, flowing through a cylindrical 2D domain,
whose lateral boundary is modeled as a  cylindrical Koiter shell. See Figure~\ref{fig:domain}.
Two Koiter shell models are considered: the linearly viscoelastic and the linearly elastic model.

The fluid flow is driven by the time-dependent inlet and outlet dynamic pressure data.
The fluid and structure are fully coupled via the kinematic and dynamic lateral boundary conditions
describing continuity of velocity (the no-slip condition), and balance of contact forces at the 
fluid-structure interface.

This work was motivated by a study of blood flow through medium-to-large human arteries. 
In medium-to-large human arteries, such as the aorta or coronary arteries, blood can be modeled as an incompressible,
viscous, Newtonian fluid. 
%Modeling elastodynamics of arterial walls, on the other hand, has been a topic of much research. 
Due to the high complexity and nonlinearity of the underlying FSI problem, various simplifications 
of the arterial wall mechanical properties have  to be taken into account.
The viscoelastic Koiter shell model, considered in this manuscript, and its many versions that appear in literature,
have been shown to be a good approximation of the mechanical properties of arterial walls \cite{SunTam,SunAnnals}.
In fact, a  version of the problem considered in this manuscript has been used as a {\sl benchmark problem} for 
FSI solvers in hemodynamics, see e.g., \cite{formaggia2001coupling}.
Therefore, studying its well-posedness has implications for both analysis and numerics.

The development of existence theory for moving boundary, fluid-structure interaction problems, has become particularly active since the late 1990's.
The first existence results were obtained for the cases  in which the structure is completely immersed in the fluid, 
and the structure was considered to be either a rigid body, or described by a finite number of modal functions. 
See e.g.,
\cite{Boulakia,CSJ,CumTak,DE2,DEGLT,Feireisl,Gal,Starovoitov}, and the references therein.
More recently, the coupling between the $2D$ or $3D$ Navier-Stokes equations and $2D$ or $3D$ linear elasticity on fixed domains,
was considered for linear models in \cite{Gunzburger}, and for nonlinear models in \cite{BarGruLasTuff2,BarGruLasTuff,KukavicaTuffahaZiane}.

Concerning compliant (elastic or viscoelastic) structures,
the first FSI existence result, locally in time, was obtained in \cite{BdV1},
where a strong solution for an interaction between an incompressible, viscous fluid in $2D$ and a $1D$ viscoelastic string
was obtained, assuming periodic boundary conditions.
This result was extended by
Lequeurre in \cite{Lequeurre}, where the existence of a unique, local in time, strong solution
for any data, and the existence of a global strong solution for small data, was proved in the case when the structure was modeled as a clamped viscoelastic beam.
D.~Coutand and S.~Shkoller proved existence, locally in time, of a unique, regular solution for
an interaction between a viscous, incompressible fluid in $3D$ and a $3D$ structure, immersed in the fluid,
where the structure was modeled by the equations of linear \cite{CSS1}, or quasi-linear \cite{CSS2} elasticity. 
In the case when the structure (solid) is modeled by a linear wave equation, 
I. Kukavica and A. Tufahha proved the existence, locally in time, of a strong solution, 
assuming lower regularity  for the initial data \cite{Kuk}.
 A fuid-structure interaction between a viscous, incompressible fluid in $3D$, and $2D$ elastic shells was considered in 
\cite{ChenShkoller,ChengShkollerCoutand} where existence, locally in time,
of the unique regular solution was proved. 
All the above mentioned existence results for strong solutions are local in time.
We also mention that the works of Shkoller et al., and Kukavica at al. were obtained in the context 
of Lagrangian coordinates, which were used for both the structure and fluid problems.

In the context of weak solutions, the following results have been obtained.
Continuous dependence of  weak solutions on initial data for a
fluid structure interaction problem with a free boundary type coupling condition was studied in \cite{GioPad}.
Existence of a weak solution for a FSI problem between a $3D$ incompressible, viscous fluid 
and a $2D$ viscoelastic plate was considered by Chambolle et al. in \cite{CDEM}, while
Grandmont improved this result in \cite{CG}  to hold for
 a $2D$ elastic plate.
 In these works existence of a weak solution was proved for as long as the elastic boundary does not touch "the bottom" (rigid)
portion of the fluid domain boundary. 

In the present manuscript we prove the existence of a weak solution to a FSI problem modeling the flow
of an incompressible, viscous, Newtonian fluid flowing through a cylinder whose lateral wall is described by 
the linearly viscoelastic, or by the linearly elastic Koiter shell equations. The fluid domain is two-dimensional, while the structure equations are in 1D.
The two existence results (the viscoelastic case and the elastic case) hold for as long as the compliant  tube walls
do not touch each other. 
The {\bf main novelty} of this work is in the methodology of proof. 
The proof is based on a semi-discrete, operator splitting Lie scheme, 
which was used in \cite{GioSun} for a design of a 
stable, loosely coupled numerical scheme, called the kinematically coupled scheme (see also \cite{MarSun}). 
Therefore, in this work, we effectively prove that the kinematically coupled scheme converges to a weak solution of the 
underlying FSI problem. To the best of our knowledge, this is the first result of this kind in the area of nonlinear FSI problems.
Semi-discretization is a well known method for proving existence of weak solutions to the Navier-Stokes equations, see e.g. \cite{Tem}, Ch. III.4.
The Lie operator splitting scheme, also known as the Marchuk-Yanenko scheme, has been widely used in numerical computations, see \cite{glowinski2003finite} and the references therein. Temam
used a combination of these approaches in \cite{TemCar} to prove the existence of a solution of the nonlinear Carleman equation.
The present manuscript represents the first use of this methodology
in the area of nonlinear FSI problems. This method is robust in the sense that 
it can be applied to the viscoelastic case and to the elastic case, independently.
Namely, our existence result in the case when the structure is purely elastic in {\em not} obtained in the limit, as the regularization provided by structural viscosity tends to zero. 

Another novelty of this work is in the the boundary conditions, which are not periodic, but are motivated by the blood flow application and are given by the 
prescribed inlet and outlet dynamic pressure data. This means that the Lagrangian framework for the treatment of the entire coupled FSI problem cannot be
used in this context, and so we employed the Arbitrary Langrangian-Eulerian (ALE) method to deal with the motion of the fluid domain.

Our proof is constructive, and its main steps have already been implemented in the design of several stable computational FSI schemes for 
the simulation of blood flow in human arteries \cite{MarSun,Fernandez,GioSun,Lukacova}.
The main steps in the proof include the {\bf ALE  weak formulation}  and the {\bf time-discretization via Lie operator splitting}.
Solutions to the time-discretized problems define a sequence of approximate solutions to the continuous time-dependent problem.
The time-dependent FSI problem is discretized in time in such a way that at each time step,
 this multi-physics problem is split into two sub-problems: the fluid and the structure sub-problem.
However, to achieve stability and convergence of the corresponding numerical scheme \cite{GioSun}, the splitting
had to be performed in a special way in which the fluid sub-problem 
includes structure inertia (and structural viscosity in the viscoelastic case)
via a  "Robin-type" boundary condition. See \cite{BorSunStability} for more details.
The fact that structure inertia (and structural viscosity)  are included implicitly in the fluid sub-problem, enabled us, in the present work, to get
appropriate energy estimates for the approximate solutions, independently of the size of time discretization.
Passing to the limit, as the size of the time step converges to zero, is achieved by the use of compactness
arguments and a careful construction of the appropriate test functions associated with moving domains.
The main difference between the viscoelastic and the elastic case is in the compactness argument.

The main body of the manuscript is dedicated to the proof in the case when the structure is modeled as a linearly viscoelastic Koiter shell.
In Section~\ref{elastic_case} we summarize the main steps of the proof in the case when the structure is modeled as a linearly elastic Koiter shell. 

\if 1 = 0
This manuscript is organized as follows. In Section~\ref{sec:problem}
we describe the underlying FSI problem and present the energy of the coupled problem.
In Section~\ref{sec:ALE_and Lie} we introduce the 
corresponding ALE formulation and  perform the Lie operator splitting, thereby introducing two sub-problems, one
for the fluid and one for the structure.
In Section~\ref{sec:weak_solutions} we define solution spaces and weak solutions.
The time discretization and the definition of approximate solutions for each $\Delta t > 0$, are introduced in
Section~\ref{sec:approximate_solutions}. In this section the main uniform energy estimates, independent of $\Delta t$, are obtained. 
Convergence of approximate solutions, as $\Delta t \to 0$, is discussed in Section~\ref{sec:convergence}.
Finally, in Section~\ref{sec:limit} we show that the limiting functions satisfy the weak form of the coupled FSI problem.
\fi

%We believe that this method can be used for the existence proof for more
%general FSI problem that have application in blood flow, for example fluid-multilayered structure interaction.

\section{Problem description}\label{sec:problem}
\begin{figure}[ht]
\centering{
\includegraphics[scale=0.35]{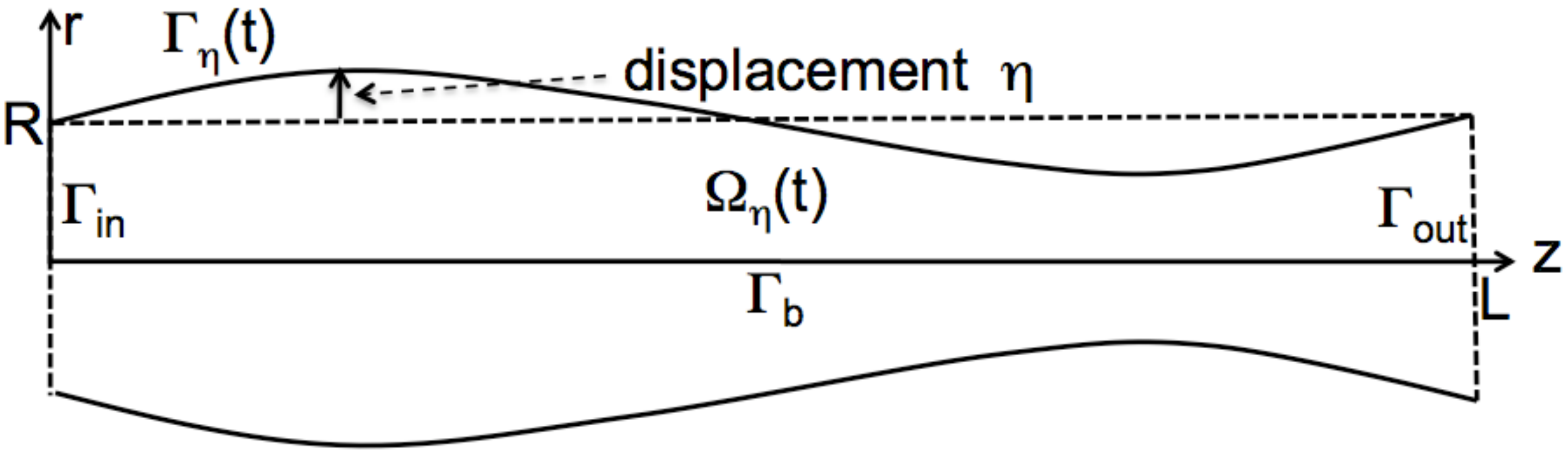}
}
\caption{Domain sketch and notation.}
\label{fig:domain2}
\end{figure}

We consider the flow of an incompressible, viscous fluid in a two-dimensional, symmetric cylinder (or channel) of reference length $L$, and reference width $2R$,
see Figure~\ref{fig:domain2}. 

Without loss of generality, we consider only the upper half of the fluid domain supplemented by a symmetry boundary condition at the bottom boundary. 
Thus, the reference domain is
$(0,L)\times (0,R)$ with the lateral (top) boundary given by
$
(0,L)\times\{R\}.
$

We will be assuming that  the lateral boundary of the cylinder is deformable
and that its location is not known {\sl a priori}, but is fully coupled to the motion of 
the viscous, incompressible fluid occupying the fluid domain. Furthermore,
it will be assumed that the lateral boundary is a thin, isotropic, homogeneous structure,
whose dynamics is modeled by the cylindrical, linearly viscoelastic, or  by the cylindrical linearly elastic Koiter shell equations.
Additionally, for simplicity, we will be assuming that only the displacement in vertical (radial) direction is non-negligible.
The (vertical) displacement from the reference configuration will be denoted by $\eta(t,z)$.
See Figure~\ref{fig:domain2}.
Models of this kind are common in blood flow applications \cite{SunTam,Pontrelli,QTV00}, where
the lateral boundary of the cylinder corresponds to arterial walls.
%Until recently, it was assumed that arterial walls deform predominantly in the radial direction.

The fluid domain, which depends on time and is not known {\sl a priori},
will be denoted by
$$
\Omega_{\eta}(t)=\{(z,r)\in\R^2:z\in (0,L),\ r\in (0,R+\eta(t,x)\},
$$
and the corresponding lateral (top) boundary by
$$
\Gamma_{\eta}(t)=\{(z,r)\in\R^2:r=R+\eta(t,z),\ z\in (0,L)\}.
$$
The ``bottom'' (symmetry) boundary of the fluid domain will be denoted by
$\Gamma_b=(0,L)\times\{0\}$, 
while the inlet and outlet sections 
of the fluid domain boundary 
by $\Gamma_{in}=\{0\}\times(0,R)$, $\Gamma_{out}=\{L\}\times (0,R)$.
See Figure~\ref{fig:domain2}.

{\bf The fluid problem:} We are interested in studying a dynamic pressure-driven flow through $\Omega_\eta(t)$
of an incompressible, viscous fluid modeled by 
the Naiver-Stokes equations:
\begin{equation}
\left .
\begin{array}{rcl}
\rho_f (\partial_t{\bf u}+{\bf u}\cdot\nabla{\bf u})&=&\nabla\cdot\sigma,
\\
\nabla\cdot{\bf u}&=&0,
\end{array}
\right \}\ {\rm in}\ \Omega_{\eta}(t),\ t\in (0,T),
\label{NS}
\end{equation}
where $\rho_f$ denotes the fluid density, ${\bf u}$ fluid velocity, $p$ fluid pressure, $\sigma=-p{\bf I}+2\mu{\bf D}({\bf u})$ is the Cauchy stress tensor of the fluid, 
$\mu$ is the kinematic viscosity coefficient,  and ${\bf D}({\bf u})=\frac 1 2(\nabla{\bf u}+\nabla^{\tau}{\bf u})$ is the symmetrized gradient of ${\bf u}$.

At the inlet and outlet boundaries we prescribe zero tangential velocity and dynamic pressure $p+\frac{\rho_f}{2}|u|^2$ (see e.g. \cite{CMP}):
\begin{equation}\label{PD}
\left. \begin{array}{rcl}
\displaystyle{p+\frac{\rho_f}{2}|u|^2}&=&P_{in/out}(t), \\
u_r &=& 0,\\
\end{array}
\right\} \quad {\rm on}\ \Gamma_{in/out},
\end{equation}
where $P_{in/out}\in L^{2}_{loc}(0,\infty)$ are given. 
Therefore the fluid flow is driven by a prescribed dynamic pressure drop, 
and the flow enters and leaves the fluid domain orthogonally to the inlet and outlet boundary. 

At the bottom boundary we prescribe the symmetry boundary condition:
\begin{equation}
u_r=\partial_r u_z=0,\quad {\rm on}\ \Gamma_b.
\label{SymBC}
\end{equation}

The {\bf structure problem}, namely, the dynamics of the lateral boundary,
is defined by the linearly viscoelastic cylindrical Koiter shell equations
capturing radial displacement $\eta$ (for the purely elastic problem see Section~\ref{elastic_case}):
\begin{equation}\label{Koiter}
\rho_s h {\partial^2_t \eta}+C_0 \eta -C_1{\partial^2_z \eta}+C_2\partial^4_z\eta+D_0{\partial_t \eta}-
D_1{\partial_t\partial^2_z\eta}+D_2\partial_t\partial^4_z\eta  =  f. 
\end{equation}
Here, $\rho_s$ is the structure density, $h$ is the structure thickness, and $f$  is the force density in the 
radial (vertical) ${\bf e}_r$ direction acting on the structure.
The constants $C_i$ and $D_i > 0$ are the material constants describing structural elasticity and viscosity, respectively,
which are given in terms of  four material parameters: the Young's modulus of elasticity $E$,
and the Poisson ratio $\sigma$,
and their viscoelastic couter-parts $E_v$ and $\sigma_v$
 (for a derivation of this model and the exact form of the coefficients, please see Appendix, and \cite{MarSun,SunTam}).
 The purely elastic case, i.e., $D_i = 0, i = 0,1,2,$ will be considered in Section~\ref{elastic_case}.

We consider the dynamics of a clamped Koiter shell with 
the boundary conditions 
$$
\eta(0)=\partial_z\eta(0)=\eta(L)=\partial_z\eta(L)= 0.
$$

The {\bf coupling} between the fluid and structure is defined by two sets of boundary conditions
satisfied at the lateral boundary $\Gamma_\eta(t)$. They are the kinematic and dynamic lateral boundary conditions
describing continuity of velocity (the no-slip condition), and continuity of normal stress, respectively.
Written in Lagrangian framework, with $z\in (0,L)$, and for  $t\in (0,T)$, they read:
\begin{itemize}
\item {\bf The\  kinematic\  condition:}
\begin{equation}\label{Coupling1a}
\displaystyle{({\partial_t \eta}(t,z),0)} = \boldsymbol{u}(t,z,R+\eta(t,z)),
\end{equation}
\item {\bf The\  dynamic\  condition:}
\begin{equation}
\label{Coupling1b}
f(t,z) = -J(t,z)(\sigma{\bf n})|_{(t,z,R+\eta(t,z))}\cdot{\bf e}_r.
\end{equation}
Here, $f=f(t,z)$ corresponds to the right hand-side of equation \eqref{Koiter}, and
$J(t,z)=\displaystyle{\sqrt{1+({\partial_z \eta}(t,z))^2}}$
denotes the Jacobian of the transformation from  Eulerian to Lagrangian coordinates.
\end{itemize}

System \eqref{NS}--\eqref{Coupling1b} 
is supplemented with the following {\bf initial conditions}:
\begin{equation}
{\bf u}(0,.)={\bf u}_0,\ \eta(0,.)=\eta_0,\ \partial_t\eta(0,.)=v_0.
\label{IC}
\end{equation}

Additionally,  we will be assuming 
that the initial data satisfies the following compatibility conditions:
\begin{equation}
\begin{array}{c}
{\bf u}_0(z,R+\eta_0(z))=v_0(z){\bf e}_r,\quad z\in (0,L),
\\
\eta_0(0)=\eta_0(L)=v_0(0)=v_0(L)=0,
\\
R+\eta_0(z) > 0,\quad z\in [0,L].
\end{array}
\label{CC}
\end{equation}
Notice that the last condition requires that the initial displacement is such that
the lateral boundary does not touch the bottom of the domain. 
This is an important condition which will be used at several places
throughout this manuscript.

In summary, we study the following fluid-structure interaction problem:
\begin{problem}\label{FSI}
Find ${\bf u}=(u_z(t,z,r),u_r(t,z,r)),p(t,z,r$), and $\eta(t,z)$ such that
\end{problem}

\begin{equation}
\left.\begin{array}{rcl}
\rho_f\big ({\partial_t{\bf u}}+({\bf u}\cdot\nabla){\bf u}\big)&=&\nabla\cdot\sigma    \\
\nabla\cdot{\bf u}&=&0 
\end{array}
\right\}\ \textrm{in}\; \Omega_{\eta}(t),\ t\in (0,T), 
\end{equation}

\begin{equation}
\left.
\begin{array}{rcl}
{\bf u}&=&\displaystyle{{\partial_t\eta}{\bf e}_r}, \\
\rho_s h \partial^2_t\eta+C_0 \eta -C_1 {\partial^2_z \eta}+C_2\partial^4_z\eta &&\\
+D_0\partial_t\eta-D_1\partial_t\partial^2_z\eta+D_2\partial_t\partial^4_z\eta&=&-J\sigma{\bf n}\cdot {\bf e}_r,  
\end{array}
\right\}\ \textrm{on}\;  (0,T)\times(0,L),
\end{equation}

\begin{equation}
\left.
\begin{array}{rcl}
u_r&=&0,\\
{\partial_r u_z}&=&0,   
\end{array}
\right\}\ \textrm{on}\;  (0,T)\times\Gamma_b,
\end{equation}

\begin{equation}
\left.
\begin{array}{rcl}
p+\frac{\rho_f}{2}|u|^2&=&P_{in/out}(t),\\
u_r&=&0,
\end{array}
\right\} \textrm{on}\; (0,T)\times \Gamma_{in/out},
\end{equation}

\begin{equation}
\left.
\begin{array}{rcl}
{\bf u}(0,.)&=&{\bf u}_0,\\
\eta(0,.)&=&\eta_0,\\
\partial_t\eta(0,.)&=&v_0.
\end{array}
\right\} \textrm{at}\;  t = 0.
\end{equation}

This is a nonlinear, moving-boundary problem, which captures the full, two-way fluid-structure interaction coupling.
The nonlinearity in the problem is represented by the quadratic term in the fluid equations, and by
the nonlinear coupling between the fluid and structure defined at the lateral boundary $\Gamma_\eta(t)$,
which is one of the unknowns in the problem.

\subsection{The energy of the problem}
Problem \eqref{FSI} satisfies the following energy inequality:
\begin{equation}\label{EE}
\frac{d}{dt}E(t)+D(t)\leq C(P_{in}(t),P_{out}(t)),
\end{equation}
where $E(t)$ denotes the sum of the kinetic energy of the fluid and of the structure, and the elastic energy of the Koiter shell:
\begin{equation}\label{energy}
\begin{array}{rcl}
E(t)&=&\displaystyle{\frac{\rho_f}{2}\|{\bf u}\|^2_{L^2(\Omega_{\eta}(t))}+\frac{\rho_sh}{2}\|\partial_t\eta\|^2_{L^2(0,L)}}\\
&+&
\displaystyle{\frac{1}{2}\big (C_0\|\eta\|^2_{L^2(0,L)}+C_1\|\partial_z\eta\|^2_{L^2(0,L)}+C_2\|\partial^2_z\eta\|^2_{L^2(0,L)}\big )},
\end{array}
\end{equation}
the term $D(t)$ captures dissipation due to structural and fluid viscosity:
\begin{equation}\label{dissipation}
D(t)=\mu\|{\bf D}({\bf u})\|^2_{L^2(\Omega_{\eta}(t)))}+
D_0\|\partial_t\eta\|^2_{L^2(0,L)}+D_1\|\partial_t\partial_z\eta\|^2_{L^2(0,L)}+D_2\|\partial_t\partial^2_z\eta\|^2_{L^2(0,L)},
\end{equation}
and $C(P_{in}(t),P_{out}(t)))$ is a constant which depends only on the inlet and outlet pressure data, which are both functions of time.

To show that \eqref{EE} holds, we first multiply equation (\ref{NS}) by ${\bf u}$, integrate over $\Omega_{\eta}(t)$, and formally integrate by parts to obtain:
$$
\int_{\Omega_{\eta}(t)}\rho_f\big (\partial_t{\bf u}\cdot{\bf u}+({\bf u}\cdot\nabla){\bf u}\cdot{\bf u}\big )+2\mu\int_{\Omega_{\eta}(t)}|{\bf D}{\bf u}|^2
-\int_{\partial\Omega_{\eta}(t)}(-p{\bf I}+2\mu{\bf D}({\bf u})){\bf n}(t)\cdot{\bf u}=0.
$$
To deal with the inertia term we first recall that  $\Omega_\eta(t)$ is moving in time and that the velocity of the 
lateral boundary is given by ${\bf u}|_{\Gamma(t)}$.
The transport theorem applied to the first term on the left hand-side of the above equation then gives:
$$
\int_{\Omega_{\eta}(t)}\partial_t{\bf u}\cdot{\bf u} = \frac 1 2\frac{d}{dt}\int_{\Omega_{\eta}(t)}|{\bf u}|^2-
\frac 1 2\int_{\Gamma_{\eta}(t)}|{\bf u}|^2{\bf u}\cdot{\bf n}(t).
$$
The second term on the left hand side can can be rewritten by using integration by parts, and the divergence-free condition, to obtain:
$$
\int_{\Omega_{\eta}(t)}({\bf u}\cdot\nabla){\bf u}\cdot{\bf u}=\frac 1 2\int_{\partial\Omega_{\eta}(t)}|{\bf u}|^2{\bf u}\cdot{\bf n}(t)=
\frac 1 2\big (\int_{\Gamma_{\eta}(t)}|{\bf u}|^2{\bf u}\cdot{\bf n}(t) 
$$
$$
-\int_{\Gamma_{in}}|{\bf u}|^2u_z+\int_{\Gamma_{out}}|{\bf u}|^2u_z.
\big )
$$
These two terms added together give
\begin{equation}\label{inertia}
\int_{\Omega_{\eta}(t)}\partial_t{\bf u}\cdot{\bf u}+\int_{\Omega_{\eta}(t)}({\bf u}\cdot\nabla){\bf u}\cdot{\bf u}=
\frac 1 2\frac{d}{dt}\int_{\Omega_{\eta}(t)}|{\bf u}|^2 -\frac 1 2\int_{\Gamma_{in}}|{\bf u}|^2u_z+\frac 1 2\int_{\Gamma_{out}}|{\bf u}|^2u_z.
\end{equation}

To deal with the boundary integral over $\partial\Omega_\eta(t)$, we first notice that on $\Gamma_{in/out}$
the boundary condition (\ref{PD}) implies $u_r=0$. Combined with the divergence-free condition we obtain
$\partial_z u_z=-\partial_r u_r=0$.
Now, using the fact that the normal to $\Gamma_{in/out}$ is ${\bf n}=(\mp 1,0)$ we get:
\begin{equation}\label{in_out}
\int_{\Gamma_{in/out}}(-p{\bf I}+2\mu{\bf D}({\bf u})){\bf n}\cdot{\bf u}=\int_{\Gamma_{in}}P_{in}u_z-\int_{\Gamma_{out}}P_{out}u_z.
\end{equation}
In a similar way, using the symmetry boundary conditions (\ref{SymBC}), we get:
$$
\int_{\Gamma_b}(-p{\bf I}+2\mu{\bf D}({\bf u})){\bf n}\cdot{\bf u}=0.
$$

What is left is to calculate the remaining boundary integral over $\Gamma_\eta(t)$. 
For this we consider  the Koiter shell equation \eqref{Koiter}, multiply it by $\partial_t \eta$,
and integrate by parts to obtain
\begin{eqnarray}\nonumber
&\displaystyle{\int_0^Lf\partial_t\eta = \frac{\rho_s h}{2}\frac{d}{dt}\|\partial_t\eta\|^2_{L^2(0,L)}}\\
\nonumber 
\\
\label{Koiter_energy}
&\displaystyle{+\frac{1}{2} \frac{d}{dt}\left( C_0\|\eta\|^2_{L^2(0,L)} 
+C_1\|\partial_z\eta\|^2_{L^2(0,L)}+C_2\|\partial^2_z\eta\|^2_{L^2(0,L)}\right)}\\
\nonumber
\\
&+ \displaystyle{D_0\|\partial_t\eta\|^2_{L^2(0,L)}
+D_1\|\partial_t\partial_z\eta\|^2_{L^2(0,L)} 
+D_2\|\partial_t\partial^2_z\eta\|^2_{L^2(0,L)}}.
\nonumber
\end{eqnarray}

By enforcing the dynamic coupling condition 
\eqref{Coupling1b} we obtain 
\begin{equation}\label{dynamic_energy}
-\int_{\Gamma_{\eta}(t)}\sigma{\bf n}(t)\cdot{\bf u}=-\int_0^LJ\sigma{\bf n}\cdot{\bf u}=\int_0^Lf\partial_t\eta.
\end{equation}

Finally, by combining  \eqref{dynamic_energy} with \eqref{Koiter_energy}, and by adding the remaining
contributions to the energy of the FSI problem calculated in equations \eqref{inertia} and \eqref{in_out}, one obtains the following energy equality: 
\begin{eqnarray}
\nonumber
\displaystyle{\frac{\rho_f}{2}\frac{d}{dt}\int_{\Omega_{\eta}(t)}|{\bf u}|^2
+\frac{\rho_s h}{2}\frac{d}{dt}\|\partial_t\eta\|^2_{L^2(0,L)}
+2\mu\int_{\Omega_{\eta}(t)}|{\bf D}{\bf u}|^2
+\frac{1}{2}
\frac{d}{dt}\big (C_0\|\eta\|^2_{L^2(0,L)}}
 \\ 
 \nonumber
 \\
 \label{EnergyEq}
+ \displaystyle{C_1\|\partial_z\eta\|^2_{L^2(0,L)}+C_2\|\partial^2_z\eta\|^2_{L^2(0,L)}\big )+
D_0\|\partial_t\eta\|^2_{L^2(0,L)}}\hskip 0.8in 
\\ 
\nonumber
\\
\displaystyle{+D_1\|\partial_t\partial_z\eta\|^2_{L^2(0,L)} 
+D_2\|\partial_t\partial^2_z\eta\|^2_{L^2(0,L)}=\pm P_{in/out}(t)\int_{\Sigma_{in/out}}u_z.}
\nonumber
\end{eqnarray}
By using the trace inequality and Korn inequality one can estimate:
$$
|P_{in/out}(t)\int_{\Sigma_{in/out}}u_z|\leq C |P_{in/out}|\|{\bf u}\|_{H^1(\Omega_{\eta}(t))}\leq
\frac{C}{2\epsilon}|P_{in/out}|^2+\frac{\epsilon C}{2}\|{\bf D}({\bf u})\|^2_{L^2(\Omega_{\eta}(t)}.
$$
By choosing $\epsilon$ such that $\frac{\epsilon C}{2}\leq \mu$ we get the energy inequality \eqref{EE}.

\section{The ALE formulation and Lie splitting}\label{sec:ALE_and Lie}
\subsection{First order ALE formulation}
To prove the existence of a weak solution to Problem~\ref{FSI} it is convenient to map Problem~\ref{FSI} onto a fixed
domain $\Omega$. In our approach we choose $\Omega$ to be the reference domain $\Omega=(0,L)\times(0,1)$.
We follow the approach typical of numerical methods for fluid-structure interaction problems 
and map our fluid  domain $\Omega(t)$ onto $\Omega$ by using an Arbitrary Lagrangian-Eulerian (ALE) mapping
\cite{MarSun,GioSun,donea1983arbitrary,quaini2007semi,QTV00}. 
We remark here that in our problem it is not convenient to use the Lagrangian formulation of the fluid sub-problem, 
as is done in e.g., \cite{CSS2,ChenShkoller,Kuk}, since, in our problem,
the fluid domain consists of a fixed, control volume of a cylinder, which does not follow Largangian flow.

We begin by defining a family of ALE mappings $A_{\eta}$ parameterized by $\eta$:
\begin{equation}
A_{\eta}(t):\Omega\rightarrow\Omega_{\eta}(t),\quad
A_{\eta}(t)(\tilde{z},\tilde{r}):=\left (\begin{array}{c}\tilde{z}\\(R+\eta(t,\tilde z))\tilde{r}\end{array}\right ),\quad (\tilde{z},\tilde{r})\in\Omega,
\label{RefTrans}
\end{equation} 
where $(\tilde z,\tilde r)$ denote the coordinates in the reference domain $\Omega=(0,L)\times (0,1)$.
Mapping $A_\eta(t)$ is a bijection, 
and its Jacobian is given by
\begin{equation}\label{ALE_Jacobian}
|{\rm det} \nabla A_\eta(t)| = |R + \eta(t,\tilde z)|.
\end{equation}
Composite functions with the ALE mapping will be denoted by
$$
{\bf u}^{\eta}(t,.)={\bf u}(t,.)\circ A_{\eta}(t) \quad {\rm and} \quad p^{\eta}(t,.)=p(t,.)\circ A_{\eta}(t).
$$
The derivatives of composite functions satisfy:
$$
\partial_t {\bf u}=\partial_t{\bf u}^{\eta}-({\bf w}^{\eta}\cdot\nabla^{\eta}){\bf u}^{\eta},\quad \nabla{\bf u}=\nabla^{\eta}{\bf u}^{\eta},
$$
where the ALE domain velocity, ${\bf w}^{\eta}$, and the transformed gradient, $\nabla^{\eta}$, are given by:
\begin{equation}
\displaystyle{{\bf w}^{\eta}=\partial_t\eta{\tilde r}{\bf e}_r,\quad \nabla^{\eta}=}
\left (\begin{array}{c}
\displaystyle{\partial_{\tilde{z}}-\tilde{r}\frac{\partial_z\eta}{R+\eta}\partial_{\tilde{r}}}
\\
\displaystyle{
\frac{1}{R+\eta}\partial_{\tilde{r}}}
\end{array}\right ).
\label{nablaeta}
\end{equation}
Note that 
\begin{equation}\label{nablaeta1}
\nabla^{\eta}{\bf v}=\nabla{\bf v}(\nabla A_{\eta})^{-1}.
\end{equation}
The following notation will also be useful:
$$
 \sigma^{\eta}=-p^{\eta}{\bf I}+2\mu{\bf D}^{\eta}({\bf u}^{\eta}),\quad  
{\bf D}^{\eta}({\bf u}^{\eta})=\frac 1 2(\nabla^{\eta}{\bf u}^{\eta}+(\nabla^{\eta})^{\tau}{\bf u}^{\eta}).
$$
We are now ready to rewrite Problem~\ref{FSI} in the ALE formulation. 
However, before we do that, we will make one more important step in our strategy to prove the existence of a weak solution to Problem~\ref{FSI}.
Namely, as mentioned in the Introduction, we would like to ``solve''  the coupled FSI problem by approximating the problem using
time-discretization via operator splitting, and then prove that the solution to the semi-discrete problem converges to a weak solution 
of the continuous problem, as the time-discretization step tends to zero.
To perform time discretization via operator splitting, which will be described in the next section, we need to write our FSI problem
as a first-order system in time. This will be done by replacing the second-order time-derivative of $\eta$, with the first-order time-derivative 
of the structure velocity. 
To do this, we further notice that in the coupled FSI problem, the kinematic coupling condition \eqref{Coupling1a} implies that
the structure velocity is equal to the normal trace of the fluid velocity on $\Gamma_\eta(t)$. Thus, we will introduce
a new variable, $v$, to denote this trace, and will replace $\partial_t \eta$ by $v$ {\bf everywhere} in the structure equation.
This has deep consequences both for the existence proof presented in this manuscript, as well as for the 
proof of  stability of the underlying numerical scheme, presented in \cite{BorSunStability},
as it enforces the kinematic coupling condition implicitly in all the steps of the scheme. 

Thus, Problem \ref{FSI} can be reformulated in the ALE framework,
on the reference domain $\Omega$, and written as a first-order system in time, in the following way:
\begin{problem}\label{FSIref}
Find ${\bf u}(t,\tilde{z},\tilde{r}),p(t,\tilde{z},\tilde{r}),\eta(t,\tilde{z})$, and $v(t,\tilde{z})$ such that
\end{problem}
 
\begin{equation}\label{FSIeqRef}
\left.
\begin{array}{rcl}
\rho_f\big ({\partial_t{\bf u}}+(({\bf u}-{\bf w}^{\eta})\cdot\nabla^{\eta}){\bf u}\big)&=&\nabla^{\eta}\cdot\sigma^{\eta}, \\
\nabla^{\eta}\cdot{\bf u}&=&0,
\end{array}
\right\} \textrm{in}\; (0,T)\times\Omega,  
\end{equation}

\begin{equation}
\left.
\begin{array}{rcl}
u_r&=&0,\\
{\partial_r u_z}&=&0
\end{array}
\right\} \textrm{on}\; (0,T)\times\Gamma_b,
\end{equation}

\begin{equation}
\left.
\begin{array}{rcl}
p+\frac{\rho_f}{2}|u|^2&=&P_{in/out}(t),\\
u_r&=&0,
\end{array}
\right\} \textrm{on}\; (0,T)\times\Gamma_{in/out},
\end{equation}

\begin{equation}\label{structure_ref}
\left.
\begin{array}{rcl}
{\bf u}&=&v{\bf e}_r,  \\
\partial_t \eta&=&v,    \\
\rho_s h \partial_t v+C_0 \eta -C_1 {\partial^2_z \eta}+C_2\partial^4_z\eta \\
+D_0v-D_1\partial^2_zv+D_2\partial^4_z v&=&-J\sigma{\bf n}\cdot {\bf e}_r
\end{array}
\right\} \textrm{on}\; (0,T)\times(0,L) ,
\end{equation}

\begin{equation}
{\bf u}(0,.)={\bf u}_0,\eta(0,.)=\eta_0,  v(0,.)=v_0, \quad {\rm at}\quad t = 0.
\end{equation}
Here, we have dropped the superscript $\eta$ in ${\bf u}^\eta$ for 
easier reading.

We are now ready to define the time discretization by operator splitting.
The underlying multi-physics problem will be split into the fluid and structure sub-problems,
following the different ``physics'' in the problem, but the splitting will be performed
in a particularly careful manner, so that the resulting problem defines a scheme 
which converges to a weak solution of the continuous problem
(and provides a numerical
scheme which is unconditionally stable).

\subsection{The operator splitting scheme}\label{sec:LieSplitting}
We use the Lie splitting, also known as the Marchuk-Yanenko splitting scheme.
The splitting can be summarized as follows.
Let $N\in\N$, $\Delta t=T/N$ and $t_n=n\Delta t$.
Consider the following
initial-value problem:
$$
\frac{d\phi}{dt}+A\phi=0\quad {\rm in}\ (0,T),\quad \phi(0)=\phi_0,
$$
where $A$ is an operator defined on a Hilbert space, and $A$ can be written as $A=A_1+A_2$. 
Set $\phi^0=\phi_0$, and, for $n=0,\dots, N-1$ and $i=1,2$, compute $\phi^{n+\frac i 2}$ by solving
\begin{equation*}
\left . 
\begin{array}{rcl}
\displaystyle{\frac{d}{dt}\phi_i+A_i\phi_i}&=&0\\
\phi_i(t_n)&=&\phi^{n+\frac{i-1}{2}}
\end{array}
\right\}\quad {\rm in}\ (t_n,t_{n+1}),
\end{equation*}
and then set
$
\phi^{n+\frac i 2}=\phi_i(t_{n+1}), \ {\rm for}\ i = 1,2.
$
It can be shown that this method is first-order accurate in time, see e.g.,  \cite{glowinski2003finite}.

We apply this approach to split Problem~\ref{FSIref} in the fluid and structure sub-problems. During this procedure the structure equation \eqref{structure_ref} 
will be split into the viscous part, i.e., the part involving the normal trace of the fluid velocity on $\Gamma_\eta(t)$, $v$, and the purely elastic part.
The viscous part of the structure problem will be used as a boundary condition for the fluid sub-problem, while the elastic part of the structure problem will be solved
separately. More precisely, we define the splitting of Problem~\ref{FSIref} in the following way:

\vskip 0.1in
{\bf Problem A1: The structure elastodynamics problem.} In this step we solve the elastodynamics problem for the location of the deformable boundary
by involving only the elastic energy of the structure. The motion of the structure is driven by the initial velocity, which
is equal to the trace of the fluid velocity on the lateral boundary, taken from the previous step. 
The fluid velocity $\bf u$ remains unchanged in this step.
More precisely, the problem reads: 
{\sl Given $({\bf u}^n,\eta^n,v^n)$ from the previous time step, find 
$({\bf u},v,\eta)$ such that:
\begin{equation}\label{Step1DF}
 \left\{\begin{array}{l@{\ }} 
{\partial_t\mathbf{u}} = 0,   \quad \textrm{in} \; (t_n, t_{n+1})\times\Omega,  \\ \\
\rho_s h \partial_t v+C_0 \eta -C_1 \partial^2_z\eta+C_2\partial^4_z\eta=0 
\quad \textrm{on} \; (t_n, t_{n+1})\times (0,L), \\ \\
\partial_t{\eta}  = v \quad \textrm{on} \; (t_n, t_{n+1})\times(0,L),\\ \\
\eta(0)=\partial_z\eta(0)=\eta(L)=\partial_z\eta(L)=0,\\ \\
{\bf u}(t_n)={\bf u}^{n},\; \eta(t_n)=\eta^n,\; v(t_n)=v^n.\end{array} \right.  
\end{equation}
Then set ${\bf u}^{n+\frac 1 2}={\bf u}(t_{n+1})$, $\eta^{n+\frac 1 2}=\eta(t_{n+1})$, $v^{n+\frac 1 2}=v(t_{n+1})$.}

\vskip 0.1in
{\bf Problem A2: The fluid problem.} In this step we solve the Navier-Stokes equations coupled with structure inertia and viscoelastic energy of the structure,
through a ``Robin-type'' boundary condition  on $\Gamma$.
The kinematic coupling condition is implicitly satisfied.
The structure displacement remains unchanged.
With a slight abuse of notation, the problem can be written as follows:
{\sl Find 
$({\bf u},v,\eta)$ such that:

\begin{eqnarray*}
 \partial_t\eta  = 0 \hskip 0.85in  &\textrm{on} \; (t_n, t_{n+1})\times(0,L),\\
\left.
\begin{array}{rcl}
 \rho_f \big (\partial_t \mathbf{u}+(({\bf u}^n-{\bf w}^{\eta^{n+\frac 1 2}})\cdot\nabla^{\eta^n}){\bf u}\big ) &=&\nabla^{\eta^n} \cdot{\sigma}^{\eta^n}   \\
 \nabla^{\eta^n} \cdot \mathbf{u}&=&0
 \end{array}
\right\} 
&\textrm{in} \; (t_n, t_{n+1})\times\Omega, 
\\
\left.
\begin{array}{rcl}
 \rho_s h\partial_t v+D_0v-D_1\partial^2_z v+D_2\partial^4_z v&=&- J\sigma{\bf n}\cdot {\bf e}_r\\
 {\bf u}&=&v{\bf e}_r
 \end{array}
\right\} 
&\textrm{on} \; (t_n, t_{n+1})\times (0,L), \\
 \left.
 \begin{array}{rcl}
 u_r&=&0\\
 {\partial_r u_z}&=&0
 \end{array}
 \right\} 
 & \textrm{on}\; (t_n, t_{n+1})\times\Gamma_b,\\ 
\left.
 \begin{array}{rcl}
 p+\frac{\rho_f}{2}|u|^2&=&P_{in/out}(t)\\
 u_r&=&0
 \end{array}
\right\}
& \textrm{on}\; (t_n, t_{n+1})\times\Gamma_{in/out},
\end{eqnarray*}
\begin{equation}\label{A2}
 {\rm with}\ 
{\bf u}(t_n,.)={\bf u}^{n+\frac 1 2},\; \eta(t_n,.)=\eta^{n+\frac 1 2}, \;  v(t_n,.)=v^{n+\frac 1 2}. 
\end{equation}
Then set $\mathbf{u}^{n+1}=\mathbf{u}(t_{n+1}), \; {\eta}^{n+1}={\eta}(t_{n+1}),\; v^{n+1}=v(t_{n+1}).$
}
\vskip 0.1in
Notice that, since  in this step $\eta$ does not change,  this problem is linear. 
Furthermore, it can be viewed as a stationary Navier-Stokes-like
problem on a fixed domain, coupled with the viscoelastic part of the structure equation
through a Robin-type boundary condition. 
In numerical simulations, one can use the ALE transformation $A_{\eta^n}$ to 
``transform'' the problem back to domain $\Omega_{\eta^n}$
and solve it there, thereby avoiding the un-necessary calculation
of the transformed gradient $\nabla^{\eta^n}$. The ALE velocity is the only extra term
that needs to be included with that approach. See, e.g., \cite{MarSun} for
more details. For the purposes of our proof, we will, however, remain in the fixed, reference domain $\Omega$.

It is important to notice that in Problem A2, the problem is ``linearized'' around the previous location of the boundary, i.e., 
we work with the domain determined by $\eta^n$,
and not by $\eta^{n+1/2}$. This is in direct relation with the implementation of the numerical scheme studied in \cite{MarSun,BorSunStability}.
However, we also notice that  ALE velocity, $w^{n+\frac 1 2}$, is taken from the just calculated Problem A1! 
This choice is {\sl crucial} for obtaining a semi-discrete version of an energy inequality, discussed in Section~\ref{sec:approximate_solutions}.

In the remainder of this paper we use the splitting scheme described above to define approximate solutions of Problem \ref{FSIref} (or equivalently Problem \ref{FSI}) and
show that the approximate solutions converge to a weak solution,  as $\Delta t\rightarrow 0$.

\section{Weak solutions}\label{sec:weak_solutions}

\subsection{Notation and function spaces}
To define weak solutions of the moving-bounday problem \ref{FSIref} we first introduce some notation which will
simplify the subsequent analysis.
We begin by introducing the following bilinear forms
associated with the elastic and viscoelastic energy of the Koiter shell:
\begin{equation}\label{Elastic}
a_S(\eta,\psi)=\int_0^L\big (C_0\eta\psi+C_1\partial_z\eta\partial_z\psi+C_2\partial^2_z\eta\partial^2_z\psi\big ),
\end{equation}
\begin{equation}\label{Viscoelastic}
a'_S(\eta,\psi)=\int_0^L\big (D_0\eta\psi+D_1\partial_{z}\eta\partial_z\psi+D_2\partial^2_{z}\eta\partial^2_{z}\psi\big ).
\end{equation}
Furthermore, we will be using $b$ to denote the following trilinear form  corresponding to the 
(symmetrized) nonlinear term in the Navier-Stokes equations:
\begin{equation}\label{transport}
b(t,{\bf u},{\bf v},{\bf w})=\frac 1 2\int_{\Omega_{\eta}(t)}({\bf u}\cdot\nabla){\bf v}\cdot{\bf w}-
\frac 1 2\int_{\Omega_{\eta}(t)}({\bf u}\cdot\nabla){\bf w}\cdot{\bf v}.
\end{equation}
Finally, we define a linear functional which associates the inlet and outlet dynamic pressure boundary data to
a test function $\bf v$ in the following way:
$$
\langle F(t),{\bf v}\rangle_{\Gamma_{in/out}}=P_{in}(t)\int_{\Gamma_{in}}v_z-P_{out}(t)\int_{\Gamma_{out}}v_z.
$$

To define a weak solution to Problem \ref{FSIref} we introduce the necessary function spaces.
For the fluid velocity we will need the classical function space
\begin{equation}\label{V_F}
\begin{array}{rl}
{\cal V}_F(t)=&\{{\mathbf u}=(u_z,u_r)\in H^1(\Omega_{\eta}(t))^2:\nabla\cdot{\bf u}=0,\\ 
 &u_z=0\ {\rm on}\ \Gamma(t),
\ u_r=0\ {\rm on}\ \Omega_{\eta}(t)\setminus\Gamma(t)\}.
\end{array}
\end{equation}
The function space associated with weak solutions of the Koiter shell is given by
\begin{equation}\label{V_S}
{\cal V}_S=H^2_0(0,L).
\end{equation}
Motivated by the energy inequality we also define the corresponding evolution spaces
for the fluid and structure sub-problems, respectively:
\begin{equation}\label{vel_test}
{\cal W}_F(0,T)=L^{\infty}(0,T;L^2(\Omega_{\eta}(t))\cap L^2(0,T;{\cal V}_F(t))
\end{equation}
\begin{equation}\label{struc_test}
{\cal W}_S(0,T)=W^{1,\infty}(0,T;L^2(0,L))\cap H^1(0,T;{\cal V}_S).
\end{equation}
The solution space for the coupled fluid-structure interaction problem must involve
the kinematic coupling condition. Thus, we define
\begin{equation}\label{W}
{\cal W}(0,T)=\{({\bf u},\eta)\in {\cal W}_F(0,T)\times{\cal W}_S(0,T):{\bf u}(t,z,R+\eta(t,z))=\partial_t\eta(t,z){\bf e}_r\}.
\end{equation}
The corresponding test space will be denoted by
\begin{equation}\label{Q}
{\cal Q}(0,T)=\{({\bf q},\psi)\in C^1_c([0,T);{\cal V}_F\times{\cal V}_S):{\bf q}(t,z,R+\eta(t,z))=\psi(t,z){\bf e}_r\}.
\end{equation}

\subsection{Weak solution on the moving domain}
We are now in a position to define weak solutions of our moving-boundary problem,
defined on the moving domain $\Omega_\eta(t)$.
\begin{definition}\label{DefWS}
We say that $({\bf u},\eta)\in{\cal W}(0,T)$ is a weak solution of Problem~\ref{FSI} if
for every $({\bf q},\psi)\in{\cal Q}(0,T)$  the  following equality holds:
\begin{equation}
\begin{array}{c}
\displaystyle{\rho_f\big (-\int_0^T\int_{\Omega_{\eta}(t)}{\bf u}\cdot\partial_t{\bf q}+\int_0^T b(t,{\bf u},{\bf u},{\bf q})\big )+
2\mu\int_0^T\int_{\Omega_{\eta}(t)}{\bf D}({\bf u}):{\bf D}({\bf q})}
\\ \\
\displaystyle{-\frac{\rho_f}{2}\int_0^T\int_0^L(\partial_t\eta)^2\psi
-\rho_sh\int_0^T\int_0^L \partial_t\eta\partial_t\psi+\int_0^T \left(a_S(\eta,\psi)+a'_S(\partial_t\eta,\psi)\right)}
\\ \\
\displaystyle{=\int_0^T\langle F(t),{\bf q}\rangle_{\Gamma_{in/out}}+\rho_f\int_{\Omega_{\eta_0}}{\bf u}_0\cdot{\bf q}(0)+\rho_s h\int_0^Lv_0\psi(0).}
\end{array}
\label{VF}
\end{equation}
\end{definition}

In deriving the weak formulation we used integration by parts in a classical way, and the following equalities which hold for smooth functions:
\begin{eqnarray*}
\int_{\Omega_{\eta}(t)}({\bf u}\cdot\nabla){\bf u}\cdot{\bf q}=&
\displaystyle{\frac 1 2\int_{\Omega_{\eta}(t)}({\bf u}\cdot\nabla){\bf u}\cdot{\bf q}-
\frac 1 2\int_{\Omega_{\eta}(t)}({\bf u}\cdot\nabla){\bf q}\cdot{\bf u}}\\
&\displaystyle{+\frac 1 2\int_0^L(\partial_t\eta)^2\psi\pm\frac 1 2\int_{\Gamma_{out/in}}|u_r|^2v_r},
\end{eqnarray*}
$$
\int_0^T\int_{\Omega_{\eta}(t)}\partial_t{\bf u}\cdot{\bf q}=-\int_0^T\int_{\Omega_{\eta}(t)}{\bf u}\cdot\partial_t{\bf q}-\int_{\Omega_{\eta_0}}{\bf u}_0\cdot{\bf q}(0)
-\int_0^T\int_0^L(\partial_t\eta)^2\psi.
$$

\subsection{Weak solution on a fixed, reference domain}
Since most of our analysis will be performed on the problem defined on the fixed, reference domain $\Omega$, 
we rewrite the above definition in terms of  $\Omega$ using the ALE mapping
$A_\eta(t)$ defined in \eqref{RefTrans}.
For this purpose, we introduce the notation for the transformed trilinear functional $b^\eta$, and the function spaces for 
the composit, transformed functions defined on the fixed domain $\Omega$.

The transformed trilinear form $b^\eta$ is defined as:
\begin{equation}\label{b_eta}
b^{\eta}({\bf u},{\bf u},{\bf q})=\frac 1 2\int_{\Omega}(R+\eta)\big ((({\bf u}-{\bf w}^{\eta})\cdot\nabla^{\eta}){\bf u}\cdot{\bf q}
-(({\bf u}-{\bf w}^{\eta})\cdot\nabla^{\eta}){\bf q}\cdot{\bf u}\big ),
\end{equation}
where $R+\eta$ is the Jacobian of the ALE mapping, calculated in \eqref{ALE_Jacobian}. 
Notice that we have included the ALE domain velocity ${\bf w}^{\eta}$ into $b^{\eta}$.

It is important to point out that
the transformed fluid velocity ${\bf u}^\eta$ is not divergence-free anymore.
Rather,  it satisfies the transformed divergence-free condition $\nabla^{\eta}\cdot{\bf u}^\eta=0$.
Therefore we need to redefine the function spaces for the fluid velocity by introducing
$$
{\cal V}_F^{\eta}=\{{\bf u}=(u_z,u_r)\in H^1(\Omega)^2:\nabla^{\eta}\cdot{\bf u}=0,\ u_z=0\ {\rm on}\ \Gamma,\ u_r=0\ {\rm on}\ \Omega\setminus\Gamma\}.
$$
The function spaces ${\cal W}_F^{\eta}(0,T)$ and ${\cal W}^{\eta}(0,T)$ are defined the same as before, but with ${\cal V}_F^{\eta}$ instead ${\cal V}_F(t)$.
More precisely:
\begin{equation}\label{W_F_eta}
{\cal W}_F^\eta(0,T)=L^{\infty}(0,T;L^2(\Omega)\cap L^2(0,T;{\cal V}_F^\eta(t)),
\end{equation}
\begin{equation}\label{W_eta}
{\cal W^\eta}(0,T)=\{({\bf u},\eta)\in {\cal W}_F^\eta(0,T)\times{\cal W}_S(0,T):{\bf u}(t,z,1)=\partial_t\eta(t,z){\bf e}_r\}.
\end{equation}
The corresponding test space is defined by
\begin{equation}\label{Q_eta}
{\cal Q^\eta}(0,T)=\{({\bf q},\psi)\in C^1_c([0,T);{\cal V}_F^\eta\times{\cal V}_S):{\bf q}(t,z,1)=\psi(t,z){\bf e}_r\}.
\end{equation}

\begin{definition}\label{DefWSRef}
We say that $({\bf u},\eta)\in{\cal W}^{\eta}(0,T)$ is a weak solution of Problem~\ref{FSIref} defined on the reference domain $\Omega$,
if for every $({\bf q},\psi)\in {\cal Q^\eta}(0,T)$ the following equality holds:
\begin{equation}
\begin{array}{c}
\displaystyle{\rho_f\big (-\int_0^T\int_{\Omega}(R+\eta){\bf u}\cdot\partial_t{\bf q}+\int_0^T b^{\eta}({\bf u},{\bf u},{\bf q})\big )}\\
\\
\displaystyle{+2\mu\int_0^T\int_{\Omega}(R+\eta){\bf D}^{\eta}({\bf u}):{\bf D}^{\eta}({\bf q})}
\displaystyle{-\frac{\rho_f}{2}\int_0^T\int_{\Omega}(\partial_t\eta){\bf u}\cdot{\bf q}}\\
\\
\displaystyle{-\rho_s h\int_0^T\int_0^L\partial_t\eta\partial_t\psi+\int_0^T\big (a_S(\eta,\psi)+a'_S(\partial_t\eta,\psi)\big )}
\\ \\
\displaystyle{=R\int_0^T\big (P_{in}(t)\int_0^1(q_z)_{|z=0}-P_{out}(t)\int_0^1(q_z)_{|z=L}\big )}\\
\\
\displaystyle{+\rho_f\int_{\Omega_{\eta_0}}{\bf u}_0\cdot{\bf q}(0)+\rho_s h\int_0^Lv_0\psi(0).}
\end{array}
\label{VFRef}
\end{equation}
\end{definition}

To see that this is consistent with the weak solution defined in Definition~\ref{DefWS}, we present the main steps in the 
transformation of the first integral on the left hand-side in \eqref{VF}, responsible for the fluid kinetic energy. 
Namely, we formally calculate:
$$
-\int_{\Omega_{\eta}}{\bf u}\cdot\partial_t{\bf q}=
-\int_{\Omega}(R+\eta){\bf u}^{\eta}\cdot(\partial_t{\bf q}-({\bf w}^{\eta}\cdot{\nabla}^{\eta}){\bf q})
=-\int_{\Omega}(R+\eta){\bf u}^{\eta}\cdot\partial_t{\bf q}
$$
$$
+ \frac 1 2\int_{\Omega}(R+\eta)({\bf w}^{\eta}\cdot{\nabla}^{\eta}){\bf q}\cdot{\bf u}^{\eta}
+ \frac 1 2\int_{\Omega}(R+\eta)({\bf w}^{\eta}\cdot{\nabla}^{\eta}){\bf q}\cdot{\bf u}^{\eta}.$$
In the last integral on the right hand-side we use the definition of ${\bf w}^\eta$ and of $\nabla^\eta$,
given in \eqref{nablaeta}, to obtain
$$
\int_{\Omega}(R+\eta)({\bf w}^{\eta}\cdot{\nabla}^{\eta}){\bf q}\cdot{\bf u}^{\eta}=\int_{\Omega}\partial_t\eta \ {\tilde r}\ \partial_{\tilde r}{\bf q}\cdot{\bf u}^{\eta}.
$$
Using integration by parts with respect to $r$, keeping in mind that $\eta$ does not depend on $r$,
we obtain
$$
-\int_{\Omega_{\eta}}{\bf u}\cdot\partial_t{\bf q}=
-\int_{\Omega}(R+\eta){\bf u}^{\eta}\cdot(\partial_t{\bf q}-({\bf w}^{\eta}\cdot{\nabla}^{\eta}){\bf q})
=-\int_{\Omega}(R+\eta){\bf u}^{\eta}\cdot\partial_t{\bf q}
$$
$$
+ \frac 1 2\int_{\Omega}(R+\eta)({\bf w}^{\eta}\cdot{\nabla}^{\eta}){\bf q}\cdot{\bf u}^{\eta}-\frac 1 2\int_{\Omega}(R+\eta)({\bf w}^{\eta}\cdot{\nabla}^{\eta}){\bf u}^{\eta}\cdot{\bf q}-\frac 1 2\int_{\Omega}\partial_t\eta{\bf u}^{\eta}\cdot{\bf q}
+\frac 1 2\int_0^L(\partial_t\eta)^2\psi,
$$
By using this identity in \eqref{VF}, and by recalling the definitions for $b$ and $b^\eta$, we obtain exactly the weak form \eqref{VFRef}.

In the remainder of this manuscript we will be working on the fluid-structure interaction problem defined on the fixed domain $\Omega$,
satisfying the weak formulation presented in Definition~\ref{DefWSRef}.
For brevity of notation, since no confusion is possible, we omit the superscript ``tilde'' which is used to denote the coordinates of points in $\Omega$.

\section{Approximate solutions}\label{sec:approximate_solutions}
In this section we use the Lie operator splitting scheme and semi-discretization to define a sequence of approximations
of a weak solution to Problem \ref{FSIref}. 
Each of the sub-problems defined by the Lie splitting in Section~\ref{sec:LieSplitting} (Problem A1 and Problem A2), 
will be discretized in time using the Backward Euler scheme. 
This approach defines a time step, which will be denoted by $\Delta t$, and a number of time sub-intervals $N\in\N$, so that
$$
(0,T) = \cup_{n=0}^{N-1} (t^n,t^{n+1}), \quad t^n = n\Delta t, \ n=0,...,N-1.
$$
For every subdivision containing $N\in\N$ sub-intervals, we recursively define the vector of unknown approximate solutions 
\begin{equation}
{\bf X}_N^{n+\frac i 2}=\left (\begin{array}{c} {\bf u}_N^{n+\frac i 2} \\ v_N^{n+\frac i 2} \\ \eta_N^{n+\frac i 2} \end{array}\right ), n=0,1,\dots,N-1,\,\ i=1,2,
\label{X}
\end{equation}
where $i = 1,2$ denotes the solution of sub-problem A1 or A2, respectively.
The initial condition will be denoted by
$$
{\bf X}^0=\left (\begin{array}{c} {\bf u}_0 \\ v_0 \\ \eta_0 \end{array} \right ).
$$

The semi-discretization and the splitting of the problem will be performed in 
such a way that the discrete version of the energy inequality \eqref{EE} is preserved 
at every time step.
This is a crucial ingredient for the existence proof.

We define the semi-discrete versions of the kinetic and elastic energy, originally defined in \eqref{energy}, and of dissipation, originally defined in \eqref{dissipation}, by
the following:
 \begin{equation}
\begin{array}{c}
\displaystyle{E_N^{n+\frac i 2}=\frac 1 2\Big (\rho_f\int_{\Omega}(R+\eta^{n-1+i})|{\bf u}^{n+\frac i 2}_N|^2+
\rho_s h\|v^{n+\frac i 2}_N\|^2_{L^2(0,L)}}\\
\displaystyle{+C_0\|\eta^{n+\frac i 2}_N\|^2_{L^2(0,L)}+C_1\|\partial_z\eta^{n+\frac i 2}_N\|^2_{L^2(0,L)}+
C_2\|\partial^2_{z}\eta^{n+\frac i 2}_N\|^2_{L^2(0,L)}\Big ),}
\end{array}
\label{kenergija}
\end{equation}
\begin{equation}
\begin{array}{c}
\displaystyle{D_N^{n+1}=\Delta t\Big (\mu\int_{\Omega}(R+\eta^n)|D^{\eta^n}({\bf u}_N^{n+1})|^2+
D_0\|v^{n+1}_N\|^2_{L^2(0,L)}
+D_1\|\partial_zv^{n+1}_N\|^2_{L^2(0,L)}}
\\
\displaystyle{+D_2\|\partial^2_{z}v^{n+1}_N\|^2_{L^2(0,L)}\Big ),\; n=0,\dots,N-1,\; i=0,1.}
\end{array}
\label{kdisipacija}
\end{equation}
Throughout the rest of this section, we fix the time step $\Delta t$, i.e., we keep $N\in\N$ fixed, and study 
the semi-discretized sub-problems defined by the Lie splitting.
To simplify notation, we will omit the subscript $N$ and write
$({\bf u}^{n+\frac i 2},v^{n+\frac i 2},\eta^{n+\frac i 2})$ instead of 
$({\bf u}^{n+\frac i 2}_N,v^{n+\frac i 2}_N,\eta^{n+\frac i 2}_N)$.

%We will show bellow that each sub-problem satisfies a discrete version of the energy inequality \eqref{EE},
%involving the discrete energy and dissipation, defined in \eqref{kenergija} and \eqref{kdisipacija}.
%This will be used later in the proof of convergence of approximate solutions to a weak solution of Problem~\ref{FSIref}.

\subsection{Semi-discretization of Problem A1}
We write a semi-discrete version of Problem A1 (Structure Elastodynamics), defined by the Lie splitting in \eqref{Step1DF}.
In this step ${\bf u}$ does not change, and so
$${\bf u}^{n+\frac 1 2}={\bf u}^n.$$
We define
$(v^{n+\frac 1 2},\eta^{n+\frac 1 2})\in H^2_0(0,L)\times H_0^2(0,L)$ as a solution of the following problem, written in weak form:
\begin{equation}
\begin{array}{c}
\displaystyle{\int_0^L\frac{\eta^{n+\frac 1 2}-\eta^{n}}{\Delta t}\phi=\int_0^Lv^{n+\frac 1 2}\phi},\quad \phi\in L^2(0,L),
\\ \\
\displaystyle{\rho_sh\int_0^L\frac{v^{n+\frac 1 2}-v^{n}}{\Delta t}\psi+a_S(\eta^{n+\frac 1 2},\psi)=0,}\quad \psi\in H^2_0(0,L).
\end{array}
\label{DProb3}
\end{equation}
The first equation is a weak form of the semi-discretized kinematic coupling condition, while the second equation
corresponds to a weak form of the semi-discretized elastodynamics equation.
\begin{proposition}
For each fixed $\Delta t > 0$, problem \eqref{DProb3} has a unique solution $(v^{n+\frac 1 2},\eta^{n+\frac 1 2})\in H^2_0(0,L)\times H_0^2(0,L)$.
\end{proposition}

\proof 
The proof is a direct consequence of the Lax-Milgram Lemma applied to the weak form
\begin{equation}\nonumber
\int_0^L\eta^{n+\frac 1 2}\psi+(\Delta t)^2a_S(\eta^{n+\frac 1 2},\psi)=\int_0^L\big (\Delta tv^n+\eta^n\big )\psi,\ \psi\in H^2_0(0,L),
\end{equation}
which is obtained after elimination of $v^{n+\frac 1 2}$ in the second equation,
by using the kinematic coupling condition given by the first equation.
\qed

\begin{proposition}\label{prop:Energy1}
For each fixed $\Delta t > 0$, solution of problem \eqref{DProb3} satisfies the following discrete energy equality:
\begin{equation}
\begin{array}{c}
\displaystyle{E_N^{n+\frac 1 2}+\frac 1 2\big (\rho_sh\|v^{n+\frac 1 2}-v^{n}\|^2+
C_0\|\eta^{n+\frac 1 2}-\eta^{n}\|^2}\\
\displaystyle{+C_1\|\partial_z(\eta^{n+\frac 1 2}-\eta^{n})\|^2+C_2\|\partial^2_{z}(\eta^{n+\frac 1 2}-\eta^{n})\|^2 \big)=E_N^n,}
\end{array}
\label{DEO31}
\end{equation}
where the kinetic energy $E_N^n$ is defined in \eqref{kenergija}.
\end{proposition}

\proof
From the first equation in (\ref{DProb3}) we immediately get 
$$v^{n+\frac 1 2}=\frac{\eta^{n+\frac 1 2}-\eta^{n}}{\Delta t}\in H_0^2(0,L).$$ 
Therefore we can take $v^{n+\frac 1 2}$ as a test function in the second equation in \eqref{DProb3}.
We replace the test function $\psi$ by $v^{n+\frac 1 2}$ in the first term on the left hand-side, 
and replace $\psi$ by $({\eta^{n+\frac 1 2}-\eta^{n}})/{\Delta t}$ in the bilinear form $a_S$. 
We then use the algebraic identity
$(a-b)\cdot a=\frac 1 2(|a|^2+|a-b|^2-|b|^2)$
to deal with the terms $(v^{n+1/2}-v^n)v^{n+1/2}$ and $(\eta^{n+1/2}-\eta^n)\eta^{n+1/2}$.
After multiplying the entire equation by $\Delta t$, 
the second equation in  \eqref{DProb3} can be written as:
$$
\rho_sh(\|v^{n+\frac 1 2}\|^2+\|v^{n+\frac 1 2}-v^{n}\|^2)+a_S(\eta^{n+\frac 1 2},\eta^{n+\frac 1 2})+
a_S(\eta^{n+\frac 1 2}-\eta^{n},\eta^{n+\frac 1 2}-\eta^{n})
$$
$$
=\rho_sh\|v^{n}\|^2+a_S(\eta^n,\eta^n).
$$
We then recall that ${\bf u}^{n+\frac 1 2}={\bf u}^n$ in this sub-problem, and so we can add 
$\rho_f \int_\Omega (1+\eta^n){\bf u}^{n+1/2}$ on the left hand-side, and
$\rho_f \int_\Omega (1+\eta^n){\bf u}^{n}$ on the right hand-side of the equation,
to obtain exactly the energy equality \eqref{DEO31}.
\qed

\subsection{Semi-discretization of Problem A2}
We write a semi-discrete version of Problem A2 (The Fluid Problem), defined by the Lie splitting in \eqref{A2}.
In this step $\eta$ does not change, and so 
$$
\eta^{n+1}=\eta^{n+\frac 1 2}.
$$
Define
$({\bf u}^{n+1},v^{n+1})\in {\cal V}_F^{\eta^n}\times H^2_0 (0,L)$ by requiring that
%\begin{problem}\label{Step2Prob}
for each $({\bf q},\psi)\in{\cal V}_F^{\eta^n}\times H^2_0 (0,L)$ such that
${\bf q}_{|\Gamma}=\psi{\bf e}_r$, the following weak formulation of problem \eqref{A2} holds:
%\end{problem}
\begin{equation}
\begin{array}{c}
\displaystyle{\rho_f\int_{\Omega}(R+\eta^{n})   \left( \frac{{\bf u}^{n+1}-{\bf u}^{n+\frac 1 2}}{\Delta t}\cdot{\bf q}+
\frac 1 2\left[({\bf u}^n-v^{n+\frac 1 2}r{\bf e}_r)\cdot\nabla^{\eta^n}\right]{\bf u}^{n+1}\cdot{\bf q}\right.}\\ \\
\displaystyle{ \left. -\frac 1 2 \left[({\bf u}^n-v^{n+\frac 1 2}r{\bf e}_r)\cdot\nabla^{\eta^n} \right] {\bf q}\cdot{\bf u}^{n+1}\right)
+\frac{\rho_f}{2} \int_{\Omega}{v^{n+\frac 1 2}}{\bf u}^{n+1}\cdot{\bf q}}\\ \\
+2\mu\int_{\Omega}(R+\eta^n){\bf D}^{\eta^n}({\bf u}):{\bf D}^{\eta^n}({\bf q})\\ \\
\displaystyle{+\rho_sh\int_0^L\frac{v^{n+1}-v^{n+\frac 1 2}}{\Delta t}\psi
+a'_S(v^{n+1},\psi)=R\big (
P^n_{in}\int_0^1(q_z)_{|z=0}-P^n_{out}\int_0^1(q_z)_{|z=L}\big ),}
\\
\\
{\rm with}\ \nabla^{\eta^n}\cdot{\bf u}^{n+1}=0,\quad {\bf u}^{n+1}_{|\Gamma}=v^{n+1}{\bf e}_r,
\end{array}
\label{D1Prob1}
\end{equation}
{\rm where}  $\displaystyle{P_{in/out}^n=\frac 1 {\Delta t}\int_{n\Delta t}^{(n+1)\Delta t}P_{in/out}(t)dt}$.

\begin{proposition}\label{existenceA1}
Let $\Delta t > 0$, and assume that $\eta^n$ are such that $R+\eta^n  \ge R_{\rm min} > 0, n=0,...,N$. Then, the fluid sub-problem
defined by (\ref{D1Prob1}) has a  
unique weak solution $({\bf u}^{n+1},v^{n+1})\in {\cal V}_F^{\eta^n}\times H^2_0 (0,L)$.
\end{proposition}

\proof
The proof is again a consequence of the Lax-Milgram Lemma.
More precisely, donote by ${\cal U}$ the Hilbert space
\begin{equation}
{\cal U}=\{({\bf u},v)\in {\cal V}_F^{\eta^n}\times H^2_0 (0,L):{\bf u}_{|\Gamma}=v{\bf e}_z\},
\label{FSProb2}
\end{equation}
and define the bilinear form associated with problem \eqref{D1Prob1}:
\begin{equation*}
\begin{array}{rcl}
\displaystyle{a(({\bf u},v),({\bf q},\psi)) } &:=& 
\displaystyle{ \rho_f\int_{\Omega}(R+\eta^n)\left( {\bf u}\cdot{\bf q}
+\frac{\Delta t}{2}\left[ ({\bf u}^n-v^{n+\frac 1 2}r{\bf e}_r)\cdot\nabla^{\eta^n}\right]  {\bf u}\cdot{\bf q}\right.}
\\
&-& \displaystyle{ \left. \frac{\Delta t}{2} \left[({\bf u}^n-v^{n+\frac 1 2}r{\bf e}_r)\cdot\nabla^{\eta^n} \right] {\bf q}\cdot{\bf u}\right)}
\\
&+& \displaystyle{\Delta t\frac{\rho_f}{2} \int_{\Omega}{v^{n+\frac 1 2}}{\bf u}\cdot{\bf q}+\Delta t2\mu\int_{\Omega}(R+\eta^n){\bf D}^{\eta^n}({\bf u}):{\bf D}^{\eta^n}({\bf q})}
\\
&+&\displaystyle{ \rho_s h\int_0^Lv\psi+\Delta t a'_S(v,\psi),\quad ({\bf u},v),({\bf q},\psi)\in{\cal U}.}
\end{array}
\end{equation*}
We need to prove that this bilinear form $a$ is coercive and continuous on ${\cal U}$.
To see that $a$ is coercive, we write
$$
a(({\bf u},v),({\bf u},v))=\rho_f\int_{\Omega}(R+\eta^n+\frac{\Delta t}{2}v^{n+\frac 1 2})|{\bf u}|^2+\rho_s h\int_0^Lv^2
$$
$$
+\Delta t(2\mu\int_{\Omega}(R+\eta^n)|{\bf D}^{\eta^n}({\bf u})|^2+a'_S(v,v)).
$$
Coercivity follows immediately after recalling that $\eta^n$ are such that $R+\eta^n \ge R_{\rm min} > 0$,
which implies that $R+\eta^n+\frac{\Delta t}{2}v^{n+\frac 1 2}=R+\frac 1 2(\eta^n+\eta^{n+\frac 1 2}) \ge R_{\rm min} >0$. 

Before we prove continuity notice that from (\ref{nablaeta}) we have:
$$
\|\nabla^{\eta^n} {\bf u}\|_{L^2(\Omega)}\leq C\|\eta^n\|_{H^2(0,L)}\|{\bf u}\|_{H^1(\Omega)}.
$$
Therefore, by applying the generalized H\"{o}lder inequality and the continuous embedding of $H^1$ into  $L^4$, we obtain
$$
a(({\bf u},v),({\bf q},\psi))\leq C\left( \rho_f\|{\bf u}\|_{L^2(\Omega)}\|{\bf q}\|_{L^2(\Omega)}+\rho_s h\|v\|_{L^2(0,L)}\|\psi\|_{L^2(0,L)}\right.
$$
$$+
\Delta t\|\eta^n\|_{H^2(0,L)} (\|{\bf u}^n\|_{H^1(\Omega)}+\|v^{n+\frac 1 2}\|_{H^1(0,L)})\|{\bf u}\|_{H^1(\Omega)}\|{\bf q}\|_{H^1(\Omega)}
$$
$$
+\left. \Delta t\mu\|\eta^n\|_{H^2(0,L)}^2\|{\bf u}\|_{H^1(\Omega)}\|{\bf q}\|_{H^1(\Omega)}+\Delta t \|v\|_{H^2(0,L)}\|\psi\|_{H^2(0,L)}\right).
$$
This shows that $a$ is continuous.
The Lax-Milgram lemma now implies the existence of a unique solution $({\bf u}^{n+1},v^{n+1})$ of problem (\ref{D1Prob1}).
\qed

\begin{proposition}
For each fixed $\Delta t > 0$, solution of problem (\ref{D1Prob1}) satisfies the following discrete energy inequality:
\begin{equation}
\begin{array}{c}
E_N^{n+1}+\displaystyle{\frac{\rho_f}{2}\int_{\Omega}(R+\eta^n)|{\bf u}^{n+1}-{\bf u}^n|^2+\frac{\rho_s h}{2}\|v^{n+1}-v^{n+\frac 1 2}\|^2_{L^2(0,L)}}\\
+D^{n+1}_N\leq E_N^{n+\frac 1 2}+C\Delta t((P_{in}^n)^2+(P_{out}^n)^2),
\end{array}
\label{DEE}
\end{equation}
where the kinetic energy $E_N^n$ and dissipation $D^{n}_N$ are defined in \eqref{kenergija} and \eqref{kdisipacija},
and the constant $C$ depends only on the parameters in the problem, and not on $\Delta t$ (or $N$).
\end{proposition} 

\proof
We begin by focusing on the weak formulation \eqref{D1Prob1} in which we replace the test functions ${\bf q}$ by ${\bf u}^{n+1}$ and $\psi$ by $v^{n+1}$.
We multiply the resulting equation by $\Delta t$, and notice that the first term on the right hand-side is given by
$$
\frac{\rho_f}{2} \int_{\Omega}(R+\eta^n)|{\bf u}^{n+1}|^2.
$$
This is the term that contributes to the discrete kinetic energy at the time step $n+1$, but it does not have the correct form, since the discrete kinetic
energy at $n+1$ is given in terms of the structure location at $n+1$, and not at $n$, namely, the discrete kinetic energy at $n+1$ involves
$$
\frac{\rho_f}{2}\int_{\Omega}(R+\eta^{n+1})|{\bf u}^{n+1}|^2.
$$
To get around this difficulty it is crucial that the advection term is present in the fluid sub-problem. 
The advection term is responsible for the presence of the integral
$$
\frac{\rho_f}{2}\int_{\Omega} \Delta tv^{n+\frac 1 2}|{\bf u}^{n+1}|^2
$$
which can be re-written by noticing that  $\Delta tv^{n+\frac 1 2}:= (\eta^{n+1/2}-\eta^n)$ which is equal to $(\eta^{n+1}-\eta^n)$ since,
in this sub-problem $\eta^{n+1} =  \eta^{n+1/2}$. This implies
$$
\frac{\rho_f}{2}\Big (\int_{\Omega}(R+\eta^n)|{\bf u}^{n+1}|^2+\Delta tv^{n+\frac 1 2}|{\bf u}^{n+1}|^2\Big )=\frac{\rho_f}{2}\int_{\Omega}(R+\eta^{n+1})|{\bf u}^{n+1}|^2.
$$
Thus, these two terms combined provide the discrete kinetic energy at the time step $n+1$.
It is interesting to notice how the nonlinearity of the coupling at the deformed  boundary requires the presence of nonlinear advection
in order for the discrete kinetic energy of the fluid sub-problem to be decreasing in time,
and to thus satisfy the desired energy estimate.

To complete the proof one simply uses the algebraic identity $(a-b)\cdot a=\frac 1 2(|a|^2+|a-b|^2-|b|^2)$
in the same way as in the proof of Proposition~\ref{prop:Energy1}.
\qed

\vskip 0.1in
We pause for a second, and summarize what we have accomplished so far.
For a given $\Delta t >0$ we divided the time interval $(0,T)$ into $N=T/\Delta t$ sub-intervals $(t^n,t^{n+1}), n = 0,...,N-1$.
On each sub-interval  $(t^n,t^{n+1})$ we ``solved'' the coupled FSI problem by applying the Lie splitting scheme. 
First we solved for the structure position (Problem A1) and then for the fluid flow (Problem A2).
We have just shown that each sub-problem has a unique solution, provided that $R+\eta^n  \ge R_{\rm min} > 0, n=0,...,N$, and that its solution satisfies
an energy estimate. When combined, the two energy estimates provide a discrete version of the energy estimate \eqref{EE}. 
Thus, for each $\Delta t$ we have a time-marching, splitting scheme which defines an approximate solution on $(0,T)$ of 
our main FSI problem defined in Problem~\ref{FSIref}, and is such that for each $\Delta t$
the approximate FSI solution satisfies a discrete 
version of the energy estimate for the continuous problem.

What we would like to ultimately show is that, as $\Delta t \to 0$, the sequence of solutions parameterized by $N$ (or $\Delta t$),
converges to a weak solution of Problem~\ref{FSIref}.
Furthermore, we also need to show that $R+\eta^n  \ge R_{\rm min} > 0$ is satisfied for each $n= 0,...,N-1$.
In order to obtain this result, it is crucial to show that the discrete energy of the approximate FSI solutions defined for each $\Delta t$,
is {\sl uniformly bounded}, independently of $\Delta t$ (or $N$).
This result is obtained by the following Lemma.

\begin{lemma}\label{stabilnost}{\bf(The uniform energy estimates)}
Let $\Delta t > 0$ and $N=T/\Delta t > 0$. 
Furthermore, let $E_N^{n+\frac 1 2}, E_N^{n+1}$, and $D_N^j$ be the kinetic energy and dissipation 
given by \eqref{kenergija} and \eqref{kdisipacija}, respectively.

There exists a constant $C>0$ independent of $\Delta t$ (and $N$), which depends only on the parameters in the problem, 
on the kinetic energy of the initial data $E_0$, and on the energy norm of the inlet and outlet data $\|P_{in/out}\|_{L^2(0,T)}^2$,
such that the following estimates hold:
\begin{enumerate}
\item
$E_N^{n+\frac 1 2}\leq C, E_N^{n+1}\leq C$, for all $ n = 0,...,N-1, $
\item $ \sum_{j=1}^ND_N^j\leq C,$
\item
$\displaystyle{\sum_{n=0}^{N-1}\left(\int_{\Omega}(R+\eta^n)|{\bf u}^{n+1}-{\bf u}^n|^2+\|v^{n+1}-v^{n+\frac 1 2}\|^2_{L^2(0,L)}\right.}$\\
\phantom {} \hskip 2.05in $\left. +\|v^{n+\frac 1 2}-v^{n}\|^2_{L^2(0,L)}\right)\leq C,$
\item
$\displaystyle{\sum_{n=0}^{N-1}\left((C_0\|\eta^{n+1}-\eta^{n}\|^2_{L^2(0,L)}+C_1\|\partial_z(\eta^{n+1}-\eta^{n})\|^2_{L^2(0,L)}\right.}$\\
\phantom {} \hskip 1.85in $\left.  + C_2\|\partial^2_z(\eta^{n+1}-\eta^{n})\|^2_{L^2(0,L)}\right)\le C.$
\end{enumerate}
In fact, $C = E_0 + \tilde{C} \left(\|P_{in}\|_{L^2(0,T)}^2 + \|P_{out}\|_{L^2(0,T)}^2\right)$, where $\tilde{C}$ is the constant from 
\eqref{DEE},
which depends only on the parameters in the problem.
\end{lemma}
\proof
We begin by adding the energy estimates (\ref{DEO31}) and (\ref{DEE}) to obtain
$$
E_N^{n+1}+D_N^{n+1}+\frac 1 2\Big (\rho_f\int_{\Omega}(R+\eta^n)|{\bf u}^{n+1}-{\bf u}^n|^2+
\rho_sh\|v^{n+1}-v^{n+\frac 1 2}\|^2_{L^2(0,L)}+
$$
$$
+\rho_sh\|v^{n+\frac 1 2}-v^{n}\|^2_{L^2(0,L)}+C_0\|\eta^{n+\frac 1 2}-\eta^{n}\|^2_{L^2(0,L)}+
C_1\|\partial_z(\eta^{n+\frac 1 2}-\eta^{n})\|^2_{L^2(0,L)}+
$$
$$
+C_2\|\partial^2_{z}(\eta^{n+\frac 1 2}-\eta^{n})\|^2_{L^2(0,L)} \Big )\leq 
E_N^{n}+ \tilde{C} \Delta t((P_{in}^n)^2+(P_{out}^n)^2),\quad n=0,\dots,N-1.
$$
Then we calculate the sum, on both sides, and cancel the same terms in the kinetic energy that appear on both sides of the inequality
to obtain
$$
E_N^{N}+\sum_{n=0}^{N-1} D_N^{n+1}
+\frac 1 2\sum_{n=0}^{N-1} \Big (\rho_f\int_{\Omega}(R+\eta^n)|{\bf u}^{n+1}-{\bf u}^n|^2+
\rho_sh\|v^{n+1}-v^{n+\frac 1 2}\|^2_{L^2(0,L)}+
$$
$$
+\rho_s h \|v^{n+\frac 1 2}-v^{n}\|^2_{L^2(0,L)}+C_0\|\eta^{n+\frac 1 2}-\eta^{n}\|^2_{L^2(0,L)}+
C_1\|\partial_z(\eta^{n+\frac 1 2}-\eta^{n})\|^2_{L^2(0,L)}+
$$
$$
+C_2\|\partial^2_{z}(\eta^{n+\frac 1 2}-\eta^{n})\|^2_{L^2(0,L)} \Big )\leq 
E_0+ \tilde{C} \Delta t  \sum_{n=0}^{N-1} ((P_{in}^n)^2+(P_{out}^n)^2).
$$
To estimate the term involving the inlet and outlet pressure
we recall that on every sub-interval $(t^n,t^{n+1})$ the pressure data is approximated 
by a constant which is equal to the average value of the pressure over that time interval.
Therefore, we have, after using H\"{o}lder's inequality: 
$$
\Delta t\sum_{n=0}^{N-1}(P_{in}^n)^2 = \Delta t\sum_{n=0}^{N-1}\left( \frac{1}{\Delta t}\int_{n\Delta t}^{(n+1)\Delta t}P_{in}(t)dt\right)^2
\le \|P_{in}\|_{L^2(0,T)}^2.
$$
By using the pressure estimate to bound the 
right hand-side in the above energy estimate, we have obtained all the statements in the Lemma,
with the constant $C$ given by
$C = E_0 + \tilde{C}  \|P_{in/out}\|_{L^2(0,T)}^2 $.

Notice that Statement 1 can be obtained in the same way by summing from $0$ to $n-1$, for each $n$,
 instead of from $0$ to $N-1$.
\qed
\vskip 0.1in
We will use this Lemma in the next section
to show convergence of approximate solutions.

\section{Convergence of approximate solutions}\label{sec:convergence}
We define approximate solutions of Problem \ref{FSIref} on $(0,T)$ to be the functions which are piece-wise constant 
on each sub-interval $((n-1)\Delta t,n\Delta t],\ n=1\dots N$ of $(0,T)$, such that for 
$t\in ((n-1)\Delta t,n\Delta t],\ n=1\dots N,$
\begin{equation}
{\bf u}_N(t,.)={\bf u}_N^n,\ \eta_N(t,.)=\eta_N^n,\ v_N(t,.)=v_N^n,\ v^*_N(t,.)=v^{n-\frac 1 2}_N.
\label{aproxNS}
\end{equation}
\begin{figure}[ht]
\centering{
\includegraphics[scale=0.6]{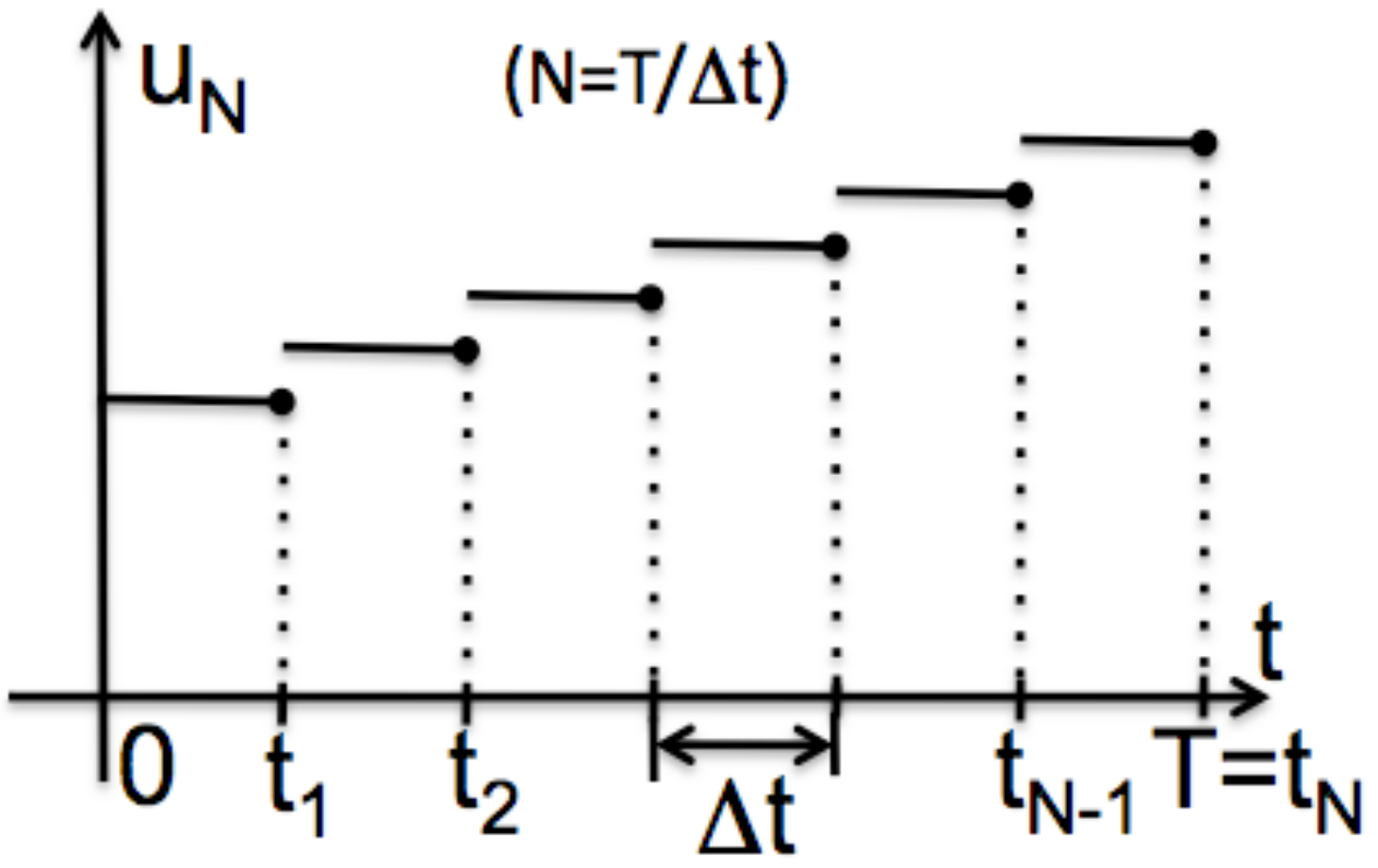}
}
\caption{A sketch of $u_N$.}
\label{fig:u_N}
\end{figure}
See Figure~\ref{fig:u_N}.
Notice that functions $v^*_N=v_N^{n-1/2}$ are determined by Step A1 (the elastodynamics sub-problem), while functions $v_N=v^{n}_N$ are determined by 
Step A2 (the fluid sub-problem). As a consequence, functions $v_N$ are equal to the normal trace of the fluid velocity on $\Gamma$,
i.e., ${\bf{u_N}} = v_N \bf e_r$. This is not necessarily the case for the functions $v^*_N$. However, we will show later that the difference between the
two sequences converges to zero in $L^2$.

Using Lemma~\ref{stabilnost} we now show that these sequences are uniformly bounded in the appropriate solution spaces.

We begin by showing that $(\eta_N)_{N\in\N}$ is uniformly bounded in $L^\infty(0,T;H_0^2(0,L))$, and that there exists a $T > 0$
for which 
$R+\eta_N^n>0$ holds independently of $N$ and $n$. This implies, among other things, that our approximate solutions are, indeed, 
well-defined on a non-zero time interval $(0,T)$.

\begin{proposition}\label{eta_bound}
Sequence $(\eta_N)_{N\in\N}$ is uniformly bounded in 
$$
L^\infty(0,T;H_0^2(0,L)).
$$
Moreover, for $T$ small enough, we have
\begin{equation}\label{eta_bounds}
0 < R_{\rm min} \le R + \eta_N(t,z) \le R_{\rm max}, \ \forall N\in \N, z\in(0,L), t\in(0,T).
\end{equation}
\end{proposition}

\proof 
%The proof relies on Lemma~\ref{stabilnost} whose proof does not require the assumption  $R+\eta_N^n>0$.
From Lemma~\ref{stabilnost} we have that $E_N^n \le C$, where $C$ is independent of $N$. This implies
$$
\|\eta_N(t)\|^2_{L^2(0,L)}, 
\|\partial_z \eta_N(t)\|^2_{L^2(0,L)}, 
\|\partial^2_{zz} \eta_N(t)\|^2_{L^2(0,L)} \le C, \ \forall t\in [0,T].
$$
Therefore, 
$$
\|\eta_N\|_{L^\infty(0,T;H_0^2(0,L))} \le C.
$$
To show that the radius $R+\eta_N$ is uniformly bounded away from zero for $T$ small enough,
we first  notice that the above inequality implies
$$
\|\eta^n_N-\eta_0\|_{H_0^2(0,L)}\leq 2C,\; n=1,\dots, N,\; N\in\N.
$$
Furthermore, we calculate
$$
\|\eta^n_N-\eta_0\|_{L^2(0,L)}\leq \sum_{i=0}^{n-1}\|\eta^{i+1}_N-\eta^{i}_N\|_{L^2(0,L)}
=\Delta t\sum_{i=0}^{n-1}\|v^{i+\frac 1 2}_N\|_{L^2(0,L)},
$$
where we recall that $\eta^{0}_N=\eta_0$.
From Lemma~\ref{stabilnost} we have that $E_N^{n+\frac 1 2}\leq C$, where $C$ is independent of $N$. This 
combined with the above inequality implies
$$
\|\eta^n_N-\eta_0\|_{L^2(0,L)}\leq C n \Delta t\leq CT,\; n=1,\dots, N,\; N\in\N.
$$
Now, we have uniform bounds for $\|\eta^n_N-\eta_0\|_{L^2(0,L)}$ and $\|\eta^n_N-\eta_0\|_{H_0^2(0,L)}$.
Therefore, we can use the interpolation inequality for Sobolev spaces (see for example \cite{ADA}, Thm. 4.17, p. 79) to get
$$
\|\eta^n_N-\eta_0\|_{H^1(0,L)}\leq 2C\sqrt{T},\; n=1,\dots, N,\; N\in\N.
$$ 
From Lemma~\ref{stabilnost} we see that $C$ depends on $T$ through the norms of the inlet and outlet data in such a way that $C$ is an increasing function of $T$.
Therefore by choosing $T$ small, we can make $\|\eta^n_N-\eta_0\|_{H^1(0,L)}$ arbitrary small for $n=1,.\dots,N$, $N\in\N$. Because of 
the Sobolev embedding of $H^1(0,L)$ into $C[0,L]$ we can also make $\|\eta^n_N-\eta_0\|_{C[0,L]}$ arbitrary small. 
Since the initial data $\eta_0$ is such that $R+\eta_0(z) > 0$ (due to the conditions listed in \eqref{CC}),
we see that for a $T>0$ small enough, there exist $R_{\min}, R_{\rm max} > 0$, such that 
$$0 < R_{\rm min} \le R + \eta_N(t,z) \le R_{\rm max}, \ \forall N\in \N, z\in(0,L), t\in(0,T).$$
\qed
\vskip 0.1in

We will show in the end that our existence result holds not only locally in time, i.e., for small $T > 0$, 
but rather, it can be extended all the way until either $T = \infty$, or until the lateral walls of the channel touch each other.

From this Proposition we see that the $L^2$-norm $\|f\|^2_{L^2(\Omega)}=\int f^2$, and the weighted $L^2$-norm
$\|f\|^2_{L^2(\Omega)}=\int (R+\eta_N) f^2$ are equivalent. More precisely, for every $f\in L^2(\Omega)$, there exist
constants $C_1, C_2 > 0$, which depend only on $R_{\rm min}, R_{\rm max}$, and not on $f$ or $N$, such that
\begin{equation}\label{equivalence}
C_1\int_{\Omega}(R+\eta_N)f^2\leq\|f\|^2_{L^2(\Omega)}\leq C_2\int_{\Omega}(R+\eta_N)f^2.
\end{equation}
We will be using this property in the next section to prove strong convergence of approximate functions.

Next we show that the sequences of approximate solutions for the velocity and its trace on the lateral boundary,
are uniformly bounded.
\begin{proposition}\label{velocity_bounds}
The following statements hold:
\begin{enumerate}
\item $(v_N)_{n\in\N}$  is uniformly bounded in $L^\infty(0,T;L^2(0,L))\cap L^2(0,T;H_0^2(0,L))$.
\item $(v_N^*)_{n\in\N}$ is uniformly bounded in $L^\infty(0,T;L^2(0,L))$.
\item $({\bf u}_N)_{n\in\N}$ is uniformly bounded in $L^\infty(0,T;L^2(\Omega))\cap L^2(0,T;H^1(\Omega))$.
\end{enumerate}

\end{proposition}
\proof
The uniform boundedness of $(v_N)_{N\in\N}, (v_N^*)_{N\in\N}$, and the uniform boundedness of $({\bf u}_N)_{N\in\N}$ in $L^\infty(0,T;L^2(\Omega))$
follow directly from Statements 1 and 2 of Lemma~\ref{stabilnost}, and from the definition of 
$(v_N)_{n\in\N}, (v_N^*)_{N\in\N}$ and $({\bf u}_N)_{N\in\N}$ as step-functions in $t$ so that 
$$
\int_0^T \|v_N\|^2_{L^2(0,L)} dt = \sum_{n=0}^{N-1} \|v_N^n\|^2_{L^2(0,L)} \Delta t.
$$
To show uniform boundedness of  $({\bf u}_N)_{N\in\N}$ in $L^2(0,T;H^1(\Omega))$ 
we need to explore the boundedness of $(\nabla {\bf u}_N)_{N \in \N}$.
From Lemma~\ref{stabilnost} we only know that the symmetrized gradient is bounded 
in the following way:
\begin{equation}\label{estKorn}
\displaystyle{\sum_{n=1}^N\int_{\Omega}(R+\eta^{n-1}_N)|{\bf D}^{\eta^{n-1}}_N({\bf u}^n_N)|^2 \Delta t\leq C.}
\end{equation}
We cannot immediately apply Korn's inequality since estimate \eqref{estKorn}
is given in terms of the transformed symmetrized gradient. Thus, there are some technical difficulties that need to be
overcome due to the fact that our problem is defined on a sequence of moving domains, and we would like
to obtain a uniform bound for the gradient $(\nabla {\bf u}_N)_{N \in \N}$.
To get around this difficulty we take the following approach.
We first transform the problem back into the original domain 
$\Omega_{\eta_N^{n-1}}$ on which $u_N$ is defined, and apply the Korn inequality in the usual way.
However, since the Korn constant depends on the domain, we will need a result which provides a 
universal Korn constant, independent of the family of domains under consideration. 
Indeed, a result of this kind was obtained in \cite{CDEM,Velcic}, assuming certain domain regularity,  which,
as we show below, 
holds for our case due to the regularity and uniform boundedness of $(\eta_N^{n-1})_{N\in \N}$.
Details are presented next.

For each fixed $N \in \N$, and for all $n = 1,\dots,N$, transform the function ${\bf u}^n_N$ back to the original domain 
which, at time step $n$, is determined by the location of the boundary $\eta_N$ at time step $n-1$, i.e., by $\eta_N^{n-1}$:
$${\bf u}(n,N)={\bf u}^n_N\circ A_{\eta^{n-1}_N}, \  n=1,\dots,N,\; N\in\N.$$
By using formula (\ref{nablaeta}) we get
$$
\int_{\Omega}(R+\eta^{n-1}_N)|{\bf D}^{\eta^{n-1}}_N({\bf u}^n_N)|^2=\int_{\Omega_{\eta^{n-1}_N}}|{\bf D}({\bf u}(n,N))|^2=
\|{\bf D}({\bf u}(n,N))\|^2_{L^2(\Omega_{\eta^{n-1}_N})}.
$$
We can now apply Korn's inequality on $\Omega_{\eta^{n-1}_N}$ to get
$$
\|\nabla{\bf u}(n,N)\|^2_{\Omega_{\eta^{n-1}_N}}\leq C(\eta^{n-1}_N)\|{\bf D}({\bf u}(j,N))\|^2_{L^2(\Omega_{\eta^{n-1}_N})},\; n=1,\dots,N,\; N\in\N,
$$
where $C(\eta^{n-1}_N)$ is the Korn's constant associated with domain $\Omega_{\eta^{n-1}_N}$.
Next, we transform everything back to $\Omega$ by using the inverse mapping, and employ (\ref{equivalence}) to obtain:
$$
\|\nabla^{\eta^{n-1}}_N{\bf u}^n_N\|_{L^2(\Omega)}\leq C(\eta^{n-1}_N)\|{\bf D}^{\eta^{n-1}}_N({\bf u}^n_N)\|_{L^2(\Omega)},\; n=1,\dots,N,\; N\in\N.
$$
Now, on the left hand-side we still have the transformed gradient $\nabla^{\eta^{n-1}}_N$ and not $\nabla$, and so we employ (\ref{nablaeta1})
to calculate the relationship between the two:
$$
\displaystyle{\nabla{\bf u}^n_N=\left(\nabla_N^{\eta^{n-1}}{\bf u}^n_N\right) \left(\nabla A_{\eta^{n-1}_N}\right),\; n=1,\dots,N,\; N\in\N.}
$$ 
Since
$\eta_N$ are bounded in $L^{\infty}(0,T;H^2(0,L))$, the gradient of the ALE mapping is bounded:
$$
\|\nabla A_{\eta^{n-1}_N}\|_{L^{\infty}(\Omega)}\leq C,\; n=1,\dots,N,\; N\in\N.
$$
Using this estimate, and by summing from $n = 1,\dots,N$, we obtain the following estimate for $\nabla {\bf u}_N^n$:
$$
\displaystyle{\sum_{n=1}^N\|\nabla{\bf u}^n_N\|^2_{L^2(\Omega)} \Delta t
\leq\sum_{n=1}^N 
C(\eta^{n-1}_N)\int_{\Omega}(R+\eta^{n-1}_N)|{\bf D}^{\eta^{n-1}}_N({\bf u}^n_N)|^2 \Delta t.}
$$
If we could show that 
$$
C(\eta^{n-1}_N)\leq K,\; n=1,\dots,N,\; N\in\N,
$$
we would have proved the uniform boundedness of ${\bf u}^n_N$ in $L^2(0,T; H^1(\Omega))$.
The existence of a uniform Korn constant $K$ follows from a result by Vel\v{c}i\'{c} \cite{Velcic} Lemma 1,  Remark 6,
summarized here in the following proposition:
\begin{proposition}{\rm \cite{Velcic}}\label{velcic}
Let $\Omega_s \subset \R^2$ be a family of open, bounded sets with  Lipschitz boundaries.
Furthermore, let us assume that the sets $\Omega_s$ are such that
$\Omega_s = F_s(\Omega)$, where $F_s$ is a family of bi-Lipschitz mappings whose bi-Lipschitz constants of
$F_s$ and $F_s^{-1}$ are uniform in $s$, and such that the family $\{F_s\}$  is strongly compact in $W^{1,\infty}(\Omega,\R^2)$.

Let ${\bf u}_s \in L^2(\Omega_s,\R^2)$ be such that the symmetrized gradient 
$D({\bf u}_s) = \frac{1}{2}\left(\nabla {\bf u}_s + (\nabla {\bf u}_s)^T \right)$ is in $L^2(\Omega_s,\R^2)$.
Then, there exists a constant $K>0$, independent of $s$, such that 
\begin{equation}
\|{\bf u}_s \|_{W^{1,2}(\Omega_s;\R^2)} \le K 
\left( \left| \int_{\Omega_s} {\bf u} dx_1 dx_2\right| + 
\left| \int_{\Omega_s} (x_1{\bf u}_2 - x_2 {\bf u}_1) dx_1 dx_2\right| + \|D({\bf u})\|_{L^2(\Omega_s)}\right).
\end{equation}
\end{proposition}
We apply this result to our problem by recalling that
Statement 1 in Lemma~\ref{stabilnost} implies
$$\|\eta^{n-1}_N\|_{H^2(0,L)}\leq C,\; n=1,\dots,N,\; N\in\N.$$
Because of the compactness of the embedding $H^2(0,L)\subset\subset W^{1,\infty}(0,L)$
and the definition of $A_{\eta^{j-1}_N}$ given in (\ref{RefTrans}),
we have 
$$
\|A_{\eta^{n-1}_N}\|_{W^{1,\infty}}\leq C,\; \|A^{-1}_{\eta^{n-1}_N}\|_{W^{1,\infty}}\leq C, n=1,\dots,N,\; N\in\N.
$$
Furthermore, the set $\{A_{\eta^{n-1}_N}:n=1,\dots,N,\; N\in\N\}$ is relatively compact in $W^{1,\infty}(\Omega)$.
Thus, by Proposition~\ref{velcic} there exists  a universal Korn constant $K > 0$ such that
$$
\displaystyle{\sum_{n=1}^N\|\nabla{\bf u}^n_N\|^2_{L^2(\Omega)} \Delta t
\leq  K \sum_{n=1}^N 
\int_{\Omega}(1+\eta^{n-1}_N)|{\bf D}^{\eta^{n-1}}_N({\bf u}^n_N)|^2 \Delta t},
$$
which implies that 
the sequence $(\nabla{\bf u}_N)_{N\in\N}$ is is uniformly bounded in $L^2((0,T)\times\Omega)$, and so
the sequence
$({\bf u}_N)_{N\in\N}$ is is uniformly bounded in $L^2(0,T;H^1(\Omega))$. 
\qed

We remark that instead of Vel\v{c}i\'{c}'s result, we could have also used the result  by Chambolle et al. in \cite{CDEM}, Lemma 6, pg. 377,
which is somewhat less general, and 
in which the main ingredient of the proof is the 
fact that the fluid velocity is divergence free, and that the displacement of
the domain boundary is only in the radial (vertical) direction.

From the uniform boundedness of approximate sequences we can now conclude that for each approximate solution sequence
there exists a subsequence which, with a slight abuse of notation, we denote the same way as the original sequence,
and  which converges weakly, or weakly*, 
depending on the function space.
More precisely, we have the following result.
\begin{lemma}{\bf (Weak and weak* convergence results)} \label{weak_convergence}
There exist subsequences $(\eta_N)_{N\in\N}, (v_N)_{N\in\N}, (v^*_N)_{N\in\N}, $ and $({\bf u}_N)_{N\in\N}$, 
and the functions $\eta \in L^{\infty}(0,T;H^2_0(0,L))$, 
$v\in L^{\infty}(0,T;L^2(0,L))\cap L^2(0,T;H^2_0(0,L))$, 
$v^*\in L^{\infty}(0,T;L^2(0,L))$, and ${\bf u}\in L^{\infty}(0,T;L^2(\Omega))\cap L^2(0,T;H^1(\Omega))$,
such that
\begin{equation}\label{weakconv}
\begin{array}{rcl}
\eta_N &\rightharpoonup & \eta \; {\rm weakly*}\; {\rm in}\; L^{\infty}(0,T;H^2_0(0,L)),
\\
v_N &\rightharpoonup & v\; {\rm weakly}\; {\rm in}\; L^2(0,T;H^2_0(0,L)),
\\
v_N &\rightharpoonup & v\; {\rm weakly*}\; {\rm in}\; L^{\infty}(0,T;L^2(0,L)),
\\
v^*_N &\rightharpoonup & v^*\; {\rm weakly*}\; {\rm in}\; L^{\infty}(0,T;L^2(0,L)),
\\
{\bf u}_N &\rightharpoonup & {\bf u}\; {\rm weakly*}\; {\rm in}\; L^{\infty}(0,T;L^2(\Omega)),
\\
{\bf u}_N &\rightharpoonup & {\bf u}\; {\rm weakly}\; {\rm in}\; L^{2}(0,T;H^1(\Omega)).
\end{array}
\end{equation}
Furthermore,
\begin{equation}\label{v_star}
v = v^*.
\end{equation}
\end{lemma}
\proof
The only thing left to show is that $v = v^*$. To show this, 
we multiply the second statement in Lemma~\ref{stabilnost} by $\Delta t$,
and 
notice again that 
$\|v_N\|_{L^2((0,T)\times(0,L))}^2=\Delta t\sum_{n=1}^N\|v^{n}_N\|^2_{L^2(0,L)}$.
This implies
$\|v_N-v^*_N\|_{L^2((0,T)\times (0,L))}\leq C\sqrt{\Delta t}$, and we have that in the limit, as $\Delta t \to 0$, $v=v^*$.
\qed

\subsection{Strong convergence of approximate sequences}

To show that the limits obtained in the previous Lemma satisfy the weak form of Problem~\ref{FSIref},
we will need to show that our sequences converge strongly in the appropriate function spaces.
To do that, we introduce the following notation which will be useful in 
the remainder of this manuscript:
denote by $\tau_h$ the translation in time by $h$ of a function $f$
\begin{equation}\label{shift}
\tau_h f(t,.)=f(t-h,.),\  h\in\R. 
\end{equation}
The strong convergence results will be achieved by using 
Corollary \ref{Corollary} listed below, of the following compactness theorem~\cite{Brezis}:
\begin{theorem}{ {\bf (Riesz-Fr\' echet-Kolmogorov Theorem)}}
Let $\Omega \subset \R^n$ be an open subset of $\R^n$, and $\omega \subset \Omega$. 
Let ${\cal{F}}$ be a bounded subset in $L^p(\Omega)$ with $1\le p < \infty$. 
Assume that 
\begin{equation*}
\forall \varepsilon > 0\quad \exists \delta > 0, \delta < {\rm dist}(\omega,\partial \Omega),\quad {\rm such \ that}
\end{equation*}
\begin{equation}\label{equicontinuity_first}
\| \tau_hf - f \|_{L^p(\omega)} < \varepsilon, \forall h \in \R^n\ {\rm such \ that}\ |h| < \delta\ {\rm and}\  \forall f \in {\cal{F}}.
\end{equation}
Then {\cal{F}} is relatively compact in $L^p(\omega)$.
\end{theorem}
Notice that $\omega\subset \Omega$ is introduced here so that the shifts $f(x\pm h)$ would be well-defined.

\begin{corollary}\label{Corollary}{\rm (Corollary 4.37, p.72 in \cite{Brezis})}
Let $\Omega \subset \R^n$ be an open subset, and let ${\cal{F}}$ be a bounded subset in $L^p(\Omega)$, for $1\le p < \infty$.
Assume that
\begin{equation}\label{equicontinuity}
\left\{
\begin{array}{c}
\forall \varepsilon > 0, \forall \omega \subset\subset \Omega, \exists \delta > 0, \delta < {\rm dist}(\omega,\partial \Omega),\quad {\rm such \ that}\\
\| \tau_hf - f \|_{L^p(\omega)} < \varepsilon, \forall h \in \R^n\ {\rm such \ that}\ |h| < \delta\ {\rm and}\  \forall f \in {\cal{F}},
\end{array}
\right.
\end{equation}
and
\begin{equation}\label{the_rest}
\forall \varepsilon > 0\quad \exists \omega \subset\subset \Omega\ {\rm such \ that}\  \|f\|_{L^p(\Omega\setminus\omega)}< \varepsilon, \ \forall f\in {\cal{F}}.
\end{equation}
Then, ${\cal{F}}$ is relatively compact in $L^p(\Omega)$.
\end{corollary}
The main ingredient in getting the ``integral equicontinuity'' estimate \eqref{equicontinuity} is
Lemma \ref{stabilnost}. Namely, if we 
multiply the third equality of Lemma \ref{stabilnost} by $\Delta t$ we get:
\begin{equation}\label{tauetaN}
\|\tau_{\Delta t}{\bf u}_{N}-{\bf u}_{N}\|^2_{L^2((0,T)\times\Omega)}
+\|\tau_{\Delta t} v_{N}-v_{N}\|^2_{L^2((0,T)\times (0,L))}\leq C {\Delta t}.
\end{equation}
This is ``almost'' \eqref{equicontinuity} except that in this estimate $\varepsilon$ depends on $\Delta t$ (i.e., $N$),
which is not sufficient to show equicontinuity \eqref{equicontinuity}. 
We need to show that estimate \eqref{equicontinuity} holds for all the functions  $(v_N)_{N\in\N}$, $({\bf u}_{N})_{N\in\N}$,
independently of $N\in\N$. This is why we need to work a little harder to get the following compactness result.

\begin{theorem}\label{u_v_convergence}
Sequences $(v_N)_{N\in\N}$, $({\bf u}_{N})_{N\in\N}$ are relatively compact in \\
$L^2(0,T;L^2(0,L))$ and $L^2(0,T;L^2(\Omega))$ respectively.
\label{LKompaktnost}
\end{theorem}
\proof The proof is based on Corollary \ref{Corollary}.
We start by showing that \eqref{equicontinuity} holds. 
 To do that, first notice that by Lemma \ref{stabilnost} both sequences are bounded (uniformly in $N\in \N$)
in the corresponding spaces. Furthermore, spatial derivatives $(\partial_z v_N)_{n\in\N}$, $(\nabla{\bf u}_{N})_{N\in\N}$ are also bounded in 
$L^2(0,T;L^2(0,L))$ and $L^2(0,T;L^2(\Omega))$ respectively,
which guarantees equicontinuity with respect to the spatial variables. 
Thus, to show that sequences $(v_N)_{N\in\N}$, $({\bf u}_{N})_{N\in\N}$ are ``equicontinuous'' in 
$L^2(0,T;L^2(0,L))$ and $L^2(0,T;L^2(\Omega))$, respectively,
we only need to consider translations in time, $\tau_h$.

We proceed by proving "integral equicontinuity'' 
\eqref{equicontinuity} for sequence $(v_N)_{N\in\N}$. Relative compactness of $({\bf u}_{N})_{N\in\N}$ can be proved analogously.

Let $\varepsilon >0$ and let $C$ be the constant from Lemma \ref{stabilnost}. 
We recall, one more time, that $C$  is independent on $N$, and thus of $\Delta t$. 

Let $\omega \subset \Omega$ be an arbitrary compact subset of $\Omega$. 
Define 
$$\delta :=\min\{{\rm dist}(\omega,\partial \Omega) /2,\varepsilon/(2C)\}.$$
We will show that
\begin{equation}\label{equicontinuity1}
\|\tau_h v_N-v_N\|^2_{L^2(\omega;L^2(0,L))}< \varepsilon, \quad \forall |h| < \delta, \ {\rm independently\ of}\ N\in\N.
\end{equation}
Thus, for each $N \in \N$, namely, for each $\Delta t = T/N$, we want to show that \eqref{equicontinuity1} holds, independently of $N$,
for each $h$ such that $|h| < \delta$.

Let  $h$ be an arbitrary real number whose absolute value is less than $\delta$.
We want to show that \eqref{equicontinuity1} holds for all $\Delta t = T/N$.
This will be shown in two steps. First, we will show that \eqref{equicontinuity1} holds for the case when $\Delta t \ge h$ (Case 1), and then for the case when $\Delta t < h$ (Case 2).

A short remark is in order: For a given $\delta > 0$, we will have $\Delta t < \delta$ for infinitely many $N$,  and both cases will apply.
For a finite number of functions $(v_N)$, we will, however, have that $\Delta t \ge \delta$. For those functions \eqref{equicontinuity1}
needs to be proved for all $\Delta t$ such that $|h| < \delta \le \Delta t$, which falls into Case 1 bellow. Thus, Cases 1 and 2 cover all
the possibilities.

\begin{figure}[ht]
\centering{
\includegraphics[scale=0.6]{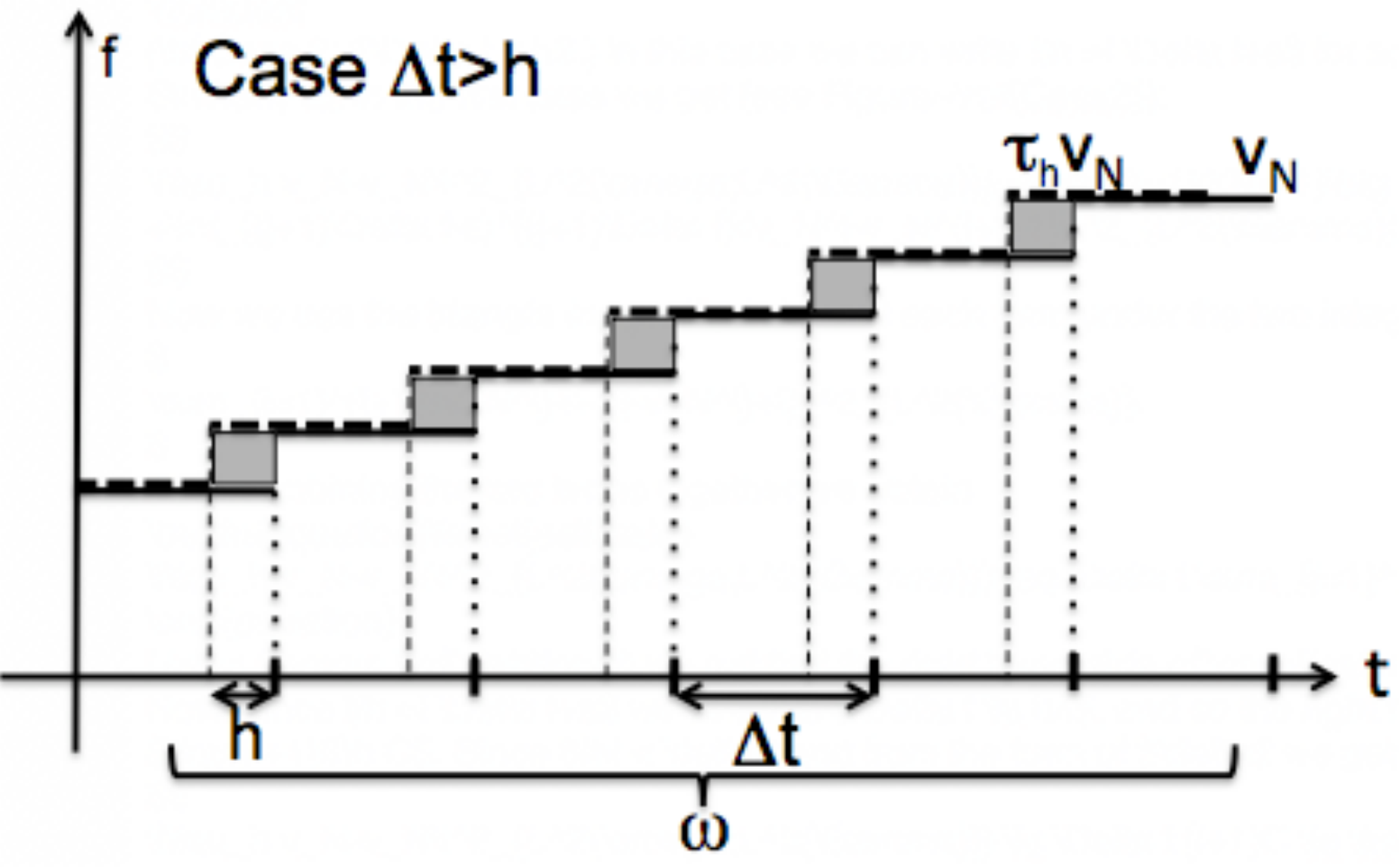}
}
\caption{Case 1: $\Delta t \ge h$. The graph of $v_N$ is shown in solid line, while the
graph of the shifted function $\tau_h v_N$ is shown in the dashed line.
The shaded area denotes the non-zero contributions to the norm $\|\tau_h v_N-v_N\|^2_{L^2}$.}
\label{fig:case1}
\end{figure}

\vskip 0.2in 
\noindent
{\bf Case 1:  $\Delta t \ge h$.} We calculate the shift by $h$ to obtain (see Figure~\ref{fig:case1}):
$$
\|\tau_h v_N-v_N\|^2_{L^2(\omega;L^2(0,L))}\leq\sum_{j=1}^{N-1}\int_{j\Delta t-h}^{j\Delta t}\|v_N^j-v_N^{j+1}\|_{L^2(0,L)}^2=
$$
$$
=h\sum_{j=1}^{N-1}\|v_N^j-v_N^{j+1}\|^2_{L^2(0,L)}\leq hC<\varepsilon/2 < \varepsilon.
$$
The last inequality follows from $|h| < \delta \le \varepsilon/(2C)$.

\begin{figure}[ht]
\centering{
\includegraphics[scale=0.6]{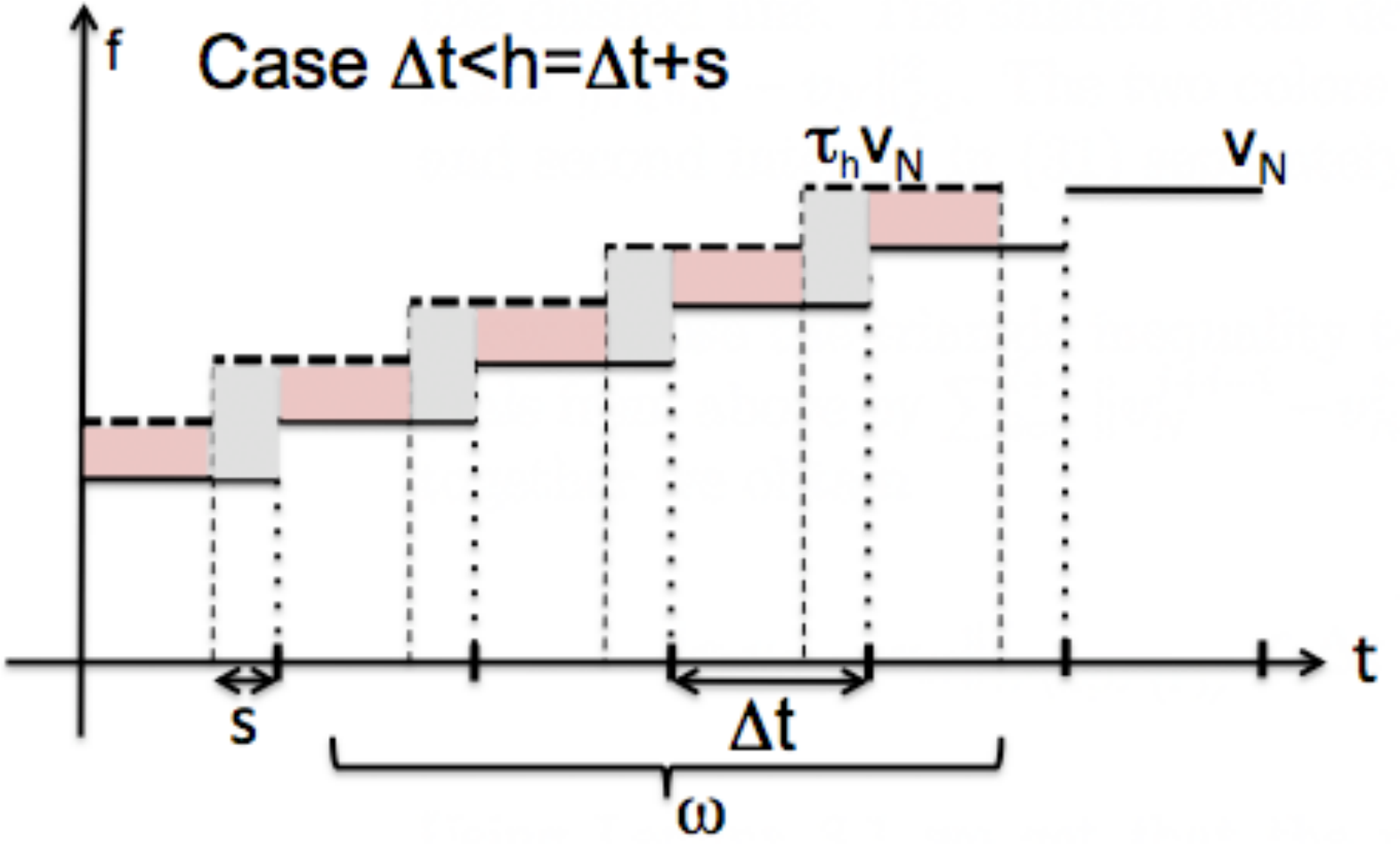}
}
\caption{Case 2: $\Delta t < h = \Delta t + s, 0<s<\Delta t$. The graph of $v_N$ is shown in solid line, while the
graph of the shifted function $\tau_h v_N$ is shown in the dashed line.
The shaded areas denote the non-zero contributions to the norm $\|\tau_h v_N-v_N\|^2_{L^2}$.
The two colors represent the contributions to the first and second integral in \eqref{est} separately.}
\label{fig:case2}
\end{figure}

\vskip 0.2in
\noindent
{\bf Case 2: $\Delta t < h$.} In this case we can write $h =l \Delta t+s$ for some $\l\in\N$, $0<s\leq \Delta t$. 
Similarly, as in the first case, we get (see Figure~\ref{fig:case2}):
\begin{equation} \label{est}
\begin{array}{rcl}
\|\tau_h v_N-v_N\|^2_{L^2(\omega;L^2(0,L))}&=&
\displaystyle{\sum_{j=1}^{N-l-1}\big (\int_{j \Delta t}^{(j+1)\Delta t-s}\|v_N^j-v_N^{j+l}\|^2_{L^2(0,L)}}\\
&+&\displaystyle{\int_{(j+1)\Delta t-s}^{(j+1)\Delta t}\|v_N^j-v_N^{j+l+1}\|^2_{L^2(0,L)}\big ).}
\end{array}
\end{equation}
Now we use the triangle inequality to bound each term under the two integrals from above by
$
\sum_{i=1}^{l+1}\|v_N^{j+i-1}-v_N^{j+i}\|^2_{L^2(0,L)}.
$
After combining the two terms together we obtain
\begin{equation}\label{estimate}
\|\tau_h v_N-v_N\|^2_{L^2(\omega;L^2(0,L))}\leq \Delta t \sum_{j=1}^{N-l-1}\sum_{i=1}^{l+1}\|v_N^{j+i-1}-v_N^{j+i}\|^2_{L^2(0,L)}.
\end{equation}
Using Lemma \ref{stabilnost} we get that the right hand-side of \eqref{estimate} is bounded by $\Delta t (l+1)C$.
Now, since $h =l \Delta t+s$ we see that $\Delta t \le h/l$, and so the right hand-side of \eqref{estimate} is bounded by 
$\frac{l+1}{l}h C$. Since $|h| < \delta$ and from the form of $\delta$ we get
$$
\|\tau_h v_N-v_N\|^2_{L^2(\omega;L^2(0,L))} \le \Delta t (l+1)C \le \frac{l+1}{l}h C \le \frac{l+1}{l} \frac{\varepsilon}{2} < \varepsilon.
$$
Thus, we have shown that \eqref{equicontinuity} holds.

To show that \eqref{the_rest} holds, let $\varepsilon > 0$. Define 
$\omega=[\varepsilon/(4C),T-\varepsilon/(4C)]$. We see that $\omega$ is obviously compact in $(0,T)$ and from the first inequality in Lemma \ref{stabilnost} 
(boundedness of $v_N^{n+\frac{i}{2}}, i = 1,2$ in $L^2(0,L)$) we have
$$\int_{(0,T)\setminus\omega}\|v_N\|^2_{L^2(0,L)}\leq \frac{\varepsilon}{2C}C<\varepsilon,\quad N\in\N.$$

\begin{figure}[ht]
\centering{
\includegraphics[scale=0.5]{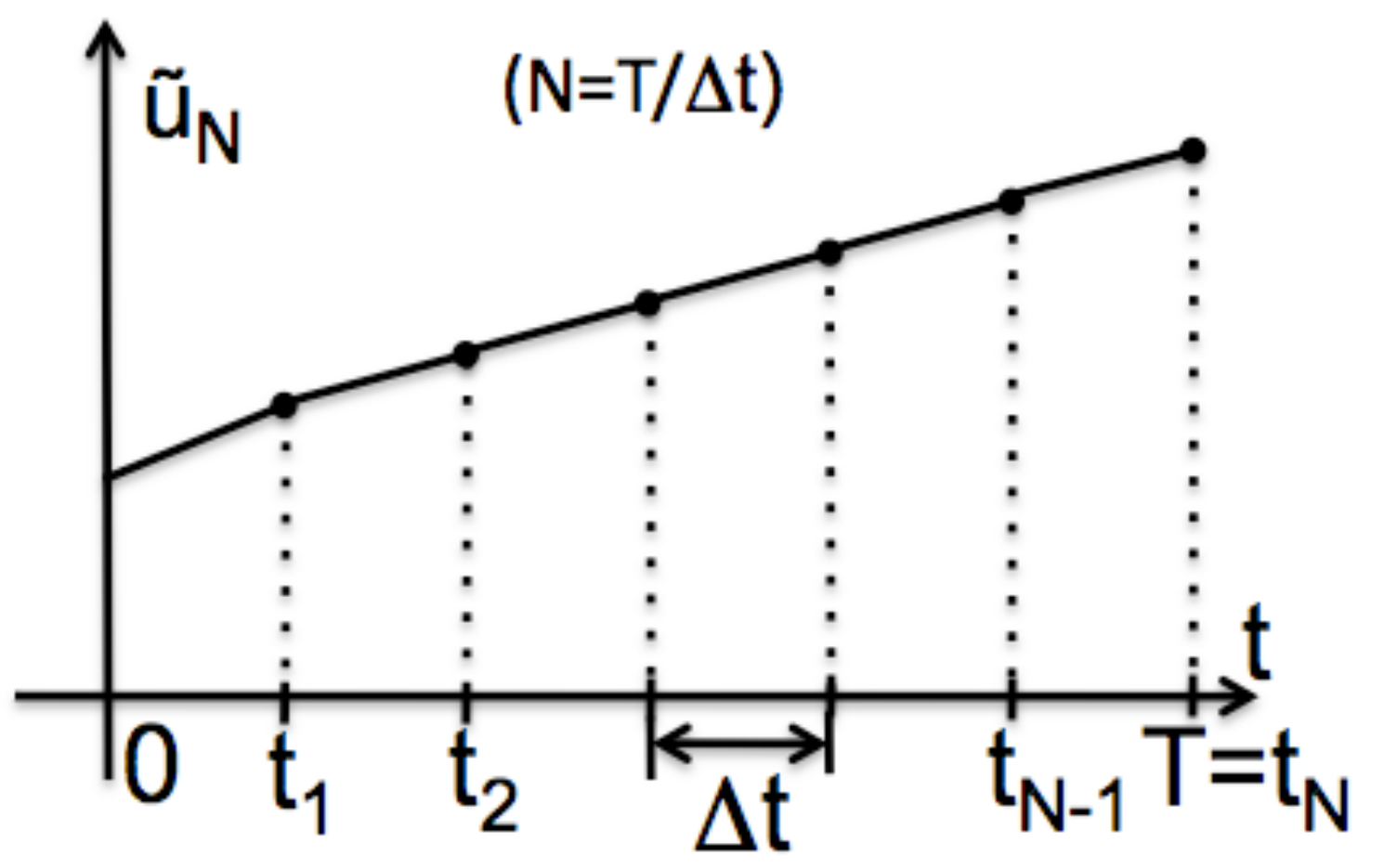}
}
\caption{A sketch of $\tilde{\bf u}_N$.}
\label{fig:u_N_tilde}
\end{figure}

By Corollary \ref{Corollary}, the compactness result for $(v_N)_{N\in\N}$ follows.
Similar arguments imply compactness of $(\mathbf{u}_N)_{N\in\N}$.
\qed
\vskip 0.1in
To show compactness of $(\eta_N)_{N\in\N}$ we introduce a slightly different set of approximate functions of $\mathbf u$, $v$, and $\eta$.
Namely, for each fixed $\Delta t$ (or $N \in \N$), define $\tilde{\bf u}_N$, $\tilde{\eta}_N$ and $\tilde{v}_N$ to be continuous, {\sl linear} on
each sub-interval $[(n-1)\Delta t,n\Delta t]$, and such that
\begin{equation}\label{tilde}
\tilde{\bf u}_N(n\Delta t,.)={\bf u}_N(n\Delta t,.),\ \tilde{v}_N(n\Delta t,.)={v}_N(n\Delta t,.),\ \tilde{\eta}_N(n\Delta t,.)={\eta}_N(n\Delta t,.), 
\end{equation}
where $n=0,\dots,N$. See Figure~\ref{fig:u_N_tilde}.
We now observe that 
$$\displaystyle{\partial_t\tilde{\eta}_N(t)=\frac{\eta^{n+1}-\eta^n}{\Delta t}= \frac{\eta^{n+1/2}-\eta^n}{\Delta t}=v^{n+\frac 1 2}},\ t\in (n\Delta t,(n+1)\Delta t),$$ 
and so, since $v^*_N$ was defined in  \eqref{aproxNS} as
a piece-wise constant function defined via $v^*_N(t,\cdot)=v^{n+\frac 1 2}$, for $t\in(n\Delta t,(n+1)\Delta t]$,
we see that
\begin{equation}\label{derivativeeta}
\partial_t\tilde{\eta}_N=v^*_N\ a.e.\ {\rm on}\ (0,T).
\end{equation}
%Functions $\tilde{\rm u}_N$ and $\tilde{v}_N$ will be used later in the proof when analysis of time-derivatives is performed. TODO

By using Lemma \ref{stabilnost} (the boundedness of $E_N^{n+\frac{i}{2}}$), we get
$$
(\tilde{\eta}_N)_{N\in\N}\; {\rm is\; bounded\; in}\; L^{\infty}(0,T;H^2_0(0,L))\cap W^{1,\infty}(0,T;L^2(0,L)).
$$
We now use the following  result on continuous embeddings:
\begin{equation}\label{compact_embedding}
L^{\infty}(0,T;H^2_0(0,L))\cap W^{1,\infty}(0,T;L^2(0,L)) \hookrightarrow C^{0,1-\alpha}([0,T]; H^{2\alpha}(0,L)),
\end{equation}
for $0<\alpha<1$.
This result follows from the standard Hilbert interpolation inequalities, see \cite{LionsMagenes}.
It was also used in \cite{CG} to deal with a set of mollifying functions approximating
a solution to a moving-boundary problem between a viscous fluid and an elastic plate.

From \eqref{compact_embedding} we see that $(\tilde{\eta}_N)_{N\in\N}$ is also bounded (uniformly in $N$) in $C^{0,1-\alpha}([0,T]; H^{2\alpha}(0,L))$.
Now, from the continuous embedding of $H^{2\alpha}(0,L)$ into $H^{2\alpha-\epsilon}$, and
by applying the Arzel\`{a}-Ascoli Theorem, we conclude that sequence $(\tilde{\eta}_N)_{N\in\N}$ has a convergent subsequence,
which we will again denote by $(\tilde{\eta}_N)_{N\in\N}$, such that
$$\tilde{\eta}_N\rightarrow \tilde{\eta}\ {\rm in} \  C([0,T];H^s(0,L)),\ 0< s<2.$$
Since sequences
$(\tilde{\eta}_N)_{N\in\N}$ and $(\eta_N)_{N\in\N}$ have the same limit, we have $\eta=\tilde{\eta}\in C([0,T];H^s(0,L))$,
where $\eta$ is the weak* limit of $(\eta_N)_{N\in\N}$, discussed in \eqref{weakconv}.
Thus, we have
$$\tilde{\eta}_N\rightarrow {\eta}\ {\rm in} \  C([0,T];H^s(0,L)),\ 0< s<2.$$
We can now prove the following Lemma:
\begin{lemma}\label{etaconvergence}
$\eta_N\rightarrow \eta$ in $L^{\infty}(0,T;H^s(0,L))$, $0 < s<2$.
\label{konveta}
\end{lemma}

\proof
The proof follows from the continuity in time of $\eta$, and from the fact that $\tilde{\eta}_N\rightarrow {\eta}\ {\rm in} \  C([0,T];H^s(0,L)),\ 0 < s<2.$
Namely, let $\varepsilon>0$. 
Then, from the continuity of $\eta$ in time we have that there exists a $\delta t>0$ such that 
$$
\|\eta(t_1)-\eta(t_2)\|_{H^s(0,L)}<\frac{\varepsilon}{2}, \  {\rm for} \ t_1,t_2\in [0,T], \ {\rm and} \ |t_1-t_2|\leq\delta t.
$$
Furthermore, from the convergence $\tilde{\eta}_N\rightarrow {\eta}\ {\rm in} \  C([0,T];H^s(0,L)),\ 0 < s<2$, we know that 
there exists an $N^* \in \N$ such that 
$$
\|\tilde{\eta}_N-\eta\|_{C([0,T];H^s(0,L))}< \frac{\varepsilon}{2}, \ \forall N\geq N^*.
$$ 
Now, let $N$ be any natural number such that $N>\max\{N^*,T/\delta t\}$. Denote by $\Delta t = T/N$, and let $t \in [0,T]$.
Furthermore, let $n\in\N$ be such that $(n-1)\Delta t<t\leq n\Delta t$.
Recall that $\tilde{\eta}_N(n\Delta t)=\eta_N(n \Delta t)=\eta_N(t)$ from the definition of $\tilde{\eta}_N$ and $\eta_N$.
By using this, and by combining the two estimates above, we get
$$
\|\eta_N(t)-\eta(t)\|_{H^s(0,L)}=
\| \eta_N(t) - \eta(n\Delta t) + \eta(n\Delta t) -\eta(t)\|_{H^s(0,L)}
$$
$$
=\| \eta_N(n\Delta t) - \eta(n\Delta t) + \eta(n\Delta t) -\eta(t)\|_{H^s(0,L)}
$$
$$
\le \| \eta_N(n\Delta t) - \eta(n\Delta t)\| + \| \eta(n\Delta t) -\eta(t)\|_{H^s(0,L)}
$$
$$
= 
\| \tilde\eta_N(n\Delta t) - \eta(n\Delta t)\|_{H^s(0,L)}  + \| \eta(n\Delta t) -\eta(t)\|_{H^s(0,L)}
<\varepsilon.
$$
Here, the first term is bounded by $\varepsilon/2$ due to the convergence $\tilde{\eta}_N \to \eta$, while
the second term is bounded by $\varepsilon/2$ due to the continuity of $\eta$. 
Since the obtained estimate is uniform in $N$ and $t$, the statement of the Lemma is proved.
\qed
\vskip 0.1in

We summarize the strong convergence results obtained in Theorem~\ref{u_v_convergence} and Lemma~\ref{etaconvergence}.
We have shown that there exist subsequences $(\mathbf u_N)_{N\in\N}$, $(v_ N)_{N\in\N}$ and $(\eta_N)_{N\in\N}$
such that
\begin{equation}
\begin{array}{rcl}
\mathbf u_N &\to& {\bf u} \  {\rm in}\  L^2(0,T;L^2(\Omega)),\\
v_ N &\to& v \  {\rm in}\  L^2(0,T;L^2(0,L)),\\
\tau_{\Delta t} \mathbf  u_N &\to& u \  {\rm in}\  L^2(0,T;L^2(\Omega)),\\
\tau_{\Delta t} v_ N &\to& v \  {\rm in}\  L^2(0,T;L^2(0,L)),\\
\eta_N &\rightarrow& \eta  \  {\rm in}\  L^{\infty}(0,T;H^s(0,L)), \ 0\leq s<2.
\end{array}
\end{equation}
Because of the uniqueness of derivatives, we also have $v=\partial_t\eta$ in the sense of distributions.
The statements about convergence of $(\tau_{\Delta t} \mathbf  u_N)_{N\in\N}$ and $(\tau_{\Delta t} v_ N)_{N\in\N}$ follow directly from \eqref{tauetaN}.

Furthermore, one can also show that 
subsequences $(\tilde{v}_N)_N$ and $(\tilde{\bf u}_N)_N$ also converge to $v$ and ${\bf u}$, respectively. 
More precisely,
\begin{equation}
\begin{array}{rcl}
\tilde{\mathbf u}_N &\to& {\bf u} \  {\rm in}\  L^2(0,T;L^2(\Omega)),\\
\tilde v_ N &\to& v \  {\rm in}\  L^2(0,T;L^2(0,L)).
\end{array}
\end{equation}
This statement follows directly from the following inequalities (see \cite{Tem}, p. 328)
$$
\|v_N-\tilde{v}_N\|^2_{L^2(0,T;L^2(0,L))}\leq{\frac{\Delta t}{3}}\sum_{n=1}^N\|v^{n+1}-v^{n}\|^2_{L^2(0,L)},
$$
$$
\|{\bf u}_N-\tilde{\bf u}_N\|^2_{L^2(0,T;L^2(\Omega))}\leq{\frac{\Delta t}{3}}\sum_{n=1}^N\|{\bf u}^{n+1}-{\bf u}^{n}\|^2_{L^2(\Omega)},
$$
and Lemma \ref{stabilnost} which provides uniform boundedness of the sums on the right hand-sides of the inequalities.

We conclude this section by showing one last convergence result that will be used in the next section to prove
that the limiting functions satisfy the weak formulation of the FSI problem. Namely,
we want to show that
\begin{equation}\label{zeta}
\begin{array}{rcl}
{\eta_N} & \to &\eta\; {\rm in}\; L^{\infty}(0,T;C^1[0,L]),\\
{\tau_{\Delta t}\eta_N} & \to &\eta\; {\rm in}\; L^{\infty}(0,T;C^1[0,L]).
\end{array}
\end{equation}
The first statement is a direct consequence of Lemma \ref{konveta}
in which we proved that $\eta_N\to\eta$ in $L^{\infty}(0,T;H^s(0,L))$, $0<s<2$.
This means that for $s>\frac 3 2$ we immediately have 
\begin{equation}
{\eta_N} \to \eta\; {\rm in}\; L^{\infty}(0,T;C^1[0,L]).
\end{equation}
To show convergence of the shifted displacements $\tau_{\Delta t} \eta_N$
to the same limiting function $\eta$, 
we recall that 
\begin{equation*}
\begin{array}{rcl}
\tilde{\eta}_N & \to &\eta\; {\rm in}\; C([0,T];H^s[0,L]),\; 0<s<2,
\end{array}
\end{equation*}
and that $(\tilde{\eta}_N)_{N\in\N}$ is uniformly bounded in $C^{0,1-\alpha}([0,T]; H^{2\alpha}(0,L))$, $0 < \alpha < 1$.
Uniform boundeness of $(\tilde{\eta}_N)_{N\in\N}$ in $C^{0,1-\alpha}([0,T]; H^{2\alpha}(0,L))$ implies that 
there exists a constant $C > 0$, independent of $N$, such that
$$
 \|  \tilde{\eta}_N((n-1)\Delta t) - \tilde\eta_N(n \Delta t)\|_{H^{2\alpha}(0,L)}  \le C |\Delta t|^{1-\alpha}.
$$
This means that for each $\varepsilon > 0$, there exists an $N_1 > 0$ such that 
$$
 \|  \tilde{\eta}_N((n-1)\Delta t) - \tilde\eta_N(n \Delta t)\|_{H^{2\alpha}(0,L)}  \le \frac{\varepsilon}{2}, {\rm for\ all}\ N \ge N_1.
$$ 
Here, $N_1$ is chosen by recalling that $\Delta t = T/N$, and so the right hand-side implies that we want an $N_1$ such that
$$
C\left(\frac{T}{N}\right)^{1-\alpha} <  \frac{\varepsilon}{2} \ {\rm for\ all}\ N\ge N_1.
$$
Now, convergence $\tilde{\eta}_N \to \eta\; {\rm in}\; C([0,T];H^s[0,L]),\; 0<s<2,$ implies that for each $\varepsilon > 0$,
there exists an $N_2 > 0$ such that
$$
 \| \tilde \eta_N (n \Delta t)  - \eta(t) \|_{H^s(0,L)} < \frac{\varepsilon}{2},\ {\rm for \ all}\ N \ge N_2.
$$
We will use this to show that for each $\varepsilon > 0$ there exists an $N^* \ge  {\rm max}\{N_1,N_2\}$, such that
$$
\displaystyle{\| \tau_{\Delta t} \tilde{\eta}_N(t) - \eta(t) \|_{H^s(0,L)} < \varepsilon, \ {\rm for \ all}\ N\ge N^*.}
$$ 
Let $t\in(0,T)$.  Then there exists an $n$  such that $t\in( (n-1)\Delta t, n \Delta t]$. We calculate
\begin{eqnarray*}
&\displaystyle{\| \tau_{\Delta t} \tilde{\eta}_N(t) - \eta(t) \|_{H^s(0,L)} =
\| \tau_{\Delta t} \tilde{\eta}_N(t) - \tilde\eta_N(n \Delta t) +\tilde \eta_N (n \Delta t)  - \eta(t) \|_{H^s(0,L)}}\\
&\displaystyle{=\|  \tilde{\eta}_N((n-1)\Delta t) - \tilde\eta_N(n \Delta t) +\tilde \eta_N (n \Delta t)  - \eta(t) \|_{H^s(0,L)} } \\
&\displaystyle{ \le \|  \tilde{\eta}_N((n-1)\Delta t) - \tilde\eta_N(n \Delta t)\|_{H^s(0,L)}   + \| \tilde \eta_N (n \Delta t)  - \eta(t) \|_{H^s(0,L)} .}
\end{eqnarray*}
The first term is less than $\varepsilon$ for all $N > N^*$ by  the uniform boundeness of $(\tilde{\eta}_N)_{N\in\N}$ in $C^{0,1-\alpha}([0,T]; H^{2\alpha}(0,L))$,
while the second term is less than $\varepsilon$ for all $N > N^*$ by the convergence of $\tilde{\eta}_N$ to $\eta$ in 
$C([0,T];H^s[0,L]),\; 0<s<2$.

Now, we notice that $\tau_{\Delta t}\tilde{\eta}_N=\widetilde{({\tau_{\Delta t}\eta}_N)}$. We  use the same argument as
in Lemma 4.1. to show that sequences
$\widetilde{({\tau_{\Delta t}\eta}_N)}$ and $\tau_{\Delta t}\eta_N$ both converge to the same limit $\eta$ in $L^{\infty}(0,T;H^s(0,L))$,
for $s<2$.

\section{The limiting problem and weak solution}\label{sec:limit}

Next we want to show that the limiting functions satisfy the weak form \eqref{VFRef} of Problem \ref{FSI}.
In this vein, one of the things that needs to be considered is what happens in the limit as $N \to \infty$, i.e., as $\Delta t \to 0$,
of problem (\ref{D1Prob1}).
Before we pass to the limit we must observe that, unfortunately, the velocity test functions in (\ref{D1Prob1}) depend of $N$!
More precisely, they depend on $\eta^n_N$ because of the requirement that
the transformed divergence-free condition $\nabla^{\eta^n_N}\cdot{\bf q} = 0$ must
be satisfied. 
This is a consequence of the fact that we mapped our problem onto a fixed domain $\Omega$.
Therefore we 
will need to take a special care in constructing the suitable velocity test function and passing to the limit in \eqref{D1Prob1}.

\subsection{Construction of the appropriate test functions}
We begin by recalling that the test functions  $({\bf q},\psi)$ for the 
limiting problem are defined by the space ${\cal Q}$, given in \eqref{Q}, which depends on $\eta$.
Similarly, the test spaces for the approximate problems depend on $N$ through the dependence on $\eta_N$.
The fact that the velocity test function depend on $N$ presents a technical difficulty when passing to the limit, as $N\to\infty$. 
To get around this difficulty, we will restrict ourselves to a dense subset ${\cal X}$ of all test functions in ${\cal Q}$ , 
which is independent of $\eta_N$ even for the approximate problems.
The set  ${\cal X}$ will consist 
of the test functions $({\bf q},\psi) \in{\cal X} = {\cal X}_F \times {\cal X}_S$, such that
the velocity components ${\bf q}$ of the test functions are smooth, independent of $N$, and $\nabla\cdot {\mathbf q} = 0$.
 
 \if 1 = 0
To get around the difficulty associated with the fact that
the velocity test functions in the approximate fluid sub-problem (\ref{D1Prob1}) depend on $N$, 
we will construct a dense subset of the velocity test function-space 
for which we will be able to show that each test function converges {\sl uniformly} on $[0,T]\times\Omega$
to a limiting test function, as $N\to\infty$. 
This will help us pass to the limit in the weak formulation.

We begin by constructing a dense  subset ${\cal X}$ of the test space ${\cal Q}$ such that  the velocity components $\mathbf q$,
of the test functions $({\bf q},\psi) \in{\cal X} = {\cal X}_F \times {\cal X}_S$, are smooth, independent of $N$, and $\nabla\cdot {\mathbf q} = 0$.
%TODO: promijeni. Tu treba samo velocity space, i ukljuciti u ${\cal X}$ sve sto je u velocity test space-u...
%
Composed with the ALE mapping $A_{\eta}$,
functions $\tilde{\mathbf q} = \mathbf q \circ A_{\eta}$, $\mathbf q \in {\cal X}_F$, define a dense subset of the test functions defined on the fixed domain $\Omega$.
Furthermore, functions $\mathbf q$ will be defined so that for each fixed $\mathbf q$, the functions
${\mathbf q}_N = \mathbf q \circ A_{{\tau_{\Delta t}\eta_N}}$ are well defined for all but possibly finitely many $N$. 
Namely, functions $\mathbf q$ will be defined so that for each  $\mathbf q$ there exists an $N_q > 0$, such that ${\mathbf q}_N$ are well defined for all $N \ge N_q$.

For each fixed $N$, i.e., $\Delta t$, functions $\mathbf q_N =  \mathbf q \circ A_{{\tau_{\Delta t}\eta_N}}, \mathbf q \in {\cal X}_F$, are also dense in the test space associated with the 
approximate weak solutions defined in (\ref{D1Prob1}). Moreover, one can easily see that $\nabla^{{\tau_{\Delta t}\eta_N}} \cdot \mathbf q_N = 0$. 

We will see that for these test functions $\mathbf q_N$ we can easily pass to the limit as $N\to\infty$ 
in the weak form \eqref{D1Prob1} and obtain the desired weak solution to Problem \ref{FSI}.
\fi
To construct the set ${\cal X}_F$ we follow ideas similar to those used in  \cite{CDEM}.
We look for the functions ${\bf q}$ which can be written as an algebraic
sum of the functions ${\bf q}_0$, which have compact support in $\Omega_{\eta}\cup\Gamma_{in}\cup\Gamma_{out}\cup\Gamma_b$, plus a function ${\bf q}_1$,
which captures the behavior of the solution at the boundary $\eta$. 
More precisely, let $\Omega_{min}$ and $\Omega_{max}$ denote the fluid domains associated with the radii $R_{min}$ and $R_{max}$,
respectively. 
\begin{enumerate}
\item {\bf Definition of test functions $({\bf q}_0,0)$ on $(0,T)\times\Omega_{max}$}: 
Consider all smooth functions ${\bf q}$ with compact support in $\Omega_{\eta}\cup\Gamma_{in}\cup\Gamma_{out}\cup\Gamma_b$,
and such that $\nabla \cdot {\bf q} = 0$.
Then we can extend ${\bf q}$ by $0$ to a divergence-free vector field on $(0,T)\times\Omega_{{max}}$.  
This defines ${\bf q}_0$. 

Notice that since ${\eta_N}$ converge uniformly to $\eta$,
there exists an $N_q > 0$ such that supp$({\bf q}_0)\subset\Omega_{{\tau_{\Delta t}\eta_N}}$, $\forall N \ge N_q$. 
Therefore, ${\bf q}_0$ is well defined  on infinitely many approximate domains $\Omega_{{\tau_{\Delta t}\eta_N}}$.
%Therefore, we can transform ${\bf q}$ to the reference domain $\Omega$ using $A_{{\tau_{\Delta t}\eta_N}}$.
\item {\bf Definition of test functions $({\bf q}_1,\psi)$ on $(0,T)\times\Omega_{max}$}: 
Consider  $\psi\in C_c^1([0,T);H^2_0(\Gamma_{\eta}))$. 
Define
\begin{equation*}
{\bf q}_1 :=\left\{
\begin{array}{l}
 \left.  
   \begin{array}{l}
       {\rm A\  constant \ extension \  in\  the\  vertical}\\
      {\rm direction\  of}  \ \psi{\bf e}_r\ {\rm on}\  \Gamma_\eta:  {\bf q}_1  := (0,\psi(z))^T;\\
       {\rm  Notice}\ {\rm div} {\bf q}_1 = 0.\\
  \end{array}
 \right\}
    {\rm on}\  \Omega_{{max}}\setminus\Omega_{{min}},\\
    \\
\left.
   \begin{array}{l}
      {\rm A \ divergence-free\  extension\  to } \  \Omega_{min}\\
     {\rm (see,\ e.g.\  \cite{GB2}, \ p.\  127).} 
  \end{array}
\right\}
    \ {\rm on}\ \Omega_{min}.
\end{array}
\right.
\end{equation*}
 From the construction
it is clear that ${\bf q}_1$ is also defined on $\Omega_{{\tau_{\Delta t}\eta_N}}$ for each $N$, and so it can be 
mapped onto the reference domain $\Omega$ by the transformation $A_{{\tau_{\Delta t}\eta_N}}$.
\end{enumerate}
For any test function $({\bf q},\psi)\in{\cal Q}$ it is easy to see that
the velocity component ${\bf q}$ can then be written as ${\bf q}={\bf q} - {\bf q}_1 + {\bf q}_1$,
where ${\bf q} - {\bf q}_1$ can be approximated by divergence-free functions ${\bf q}_0$,
which have compact support in $\Omega_{\eta}\cup\Gamma_{in}\cup\Gamma_{out}\cup\Gamma_b$.
Therefore, one can easily see that functions $({\bf q},\psi) = ({\bf q}_0 + {\bf q}_1,\psi)$ in ${\cal X}$ satisfy the following properties:
\begin{itemize}
\item
 ${\cal X}$  is dense in the space  ${\cal Q}$ of all test functions defined on the physical, moving domain $\Omega_\eta$,
defined by \eqref{Q}; furthermore, $\nabla\cdot \mathbf q = 0, \forall \mathbf q \in {\cal X}_F$.
\item Similarly, for each ${\bf q} \in {\cal X}_F$, define
$$
\tilde{\mathbf q} = \mathbf q \circ A_{\eta}.
$$
The set $\{(\tilde{\mathbf q},\psi)  |  \tilde{\mathbf q} = \mathbf q \circ A_{\eta},  \mathbf q \in {\cal X}_F, \psi \in {\cal X}_S\}$ is  dense in 
the space ${\cal Q}_\eta$ of all test functions defined on the fixed, reference domain $\Omega$,
defined by \eqref{Q_eta}.
\end{itemize}

For each ${\bf q} \in {\cal X}_F$, define
$$
 {\bf q}_N:={\bf q}\circ A_{{\tau_{\Delta t}\eta_N}}.
$$
 Functions $ {\bf q}_N$ satisfy
$\nabla^{{\tau_{\Delta t}\eta_N}}\cdot{\bf q_N}=0$.
These will serve as the test functions for approximate problems, defined on the domains determined by ${\tau_{\Delta t}\eta_N}$.
The approximate test functions have the following uniform convergence properties.
\begin{lemma}\label{testf}
For every $(\bf q,\psi)\in {\cal X}$ we have 
$${\bf q}_N\rightarrow\tilde{\bf q}\ {\rm  and}\ 
\nabla{\bf q}_N\rightarrow\nabla\tilde{\bf q}\ {\rm uniformly\  on}\  [0,T]\times\Omega.
$$
\end{lemma}
\proof
By the Mean-Value Theorem we get:
\begin{eqnarray*}
|{\bf q}_N(t,z,r)-\tilde{\bf q}(t,z,r)|&=&|{\bf q}(t,z,(R+\tau_{\Delta t}\eta_N)r)-{\bf q}(t,z,(R+\eta)r)|\\
    &=&|\partial_r{\bf q}(t,z,\zeta)r|\ |\eta(t,z)-\eta_N(t-\Delta t,z)|.
\end{eqnarray*}
The uniform convergence of ${\bf q}_N$ follows 
from the uniform convergence of $\eta_N$, since $\mathbf q$ are smooth. 

To show the uniform convergence of the gradients, one can use the chain rule to calculate
$$
\partial_z{\bf q}_N(t,z,r)=\partial_z{\bf q}(t,z,(R+\tau_{\Delta t}\eta_N)r)+\left[\partial_z\tau_{\Delta t}\eta_N(t,z)r\right]
\  \left[\partial_r{\bf q}(t,z,(R+\tau_{\Delta t}\eta_N)r)\right].
$$
The uniform convergence of $\partial_z{\bf q}_N(t,z,r)$ follows from the
uniform convergence of $\partial_z\tau_{\Delta t}\eta_N$. Combined with  the fist part of the proof we get
$\partial_z{\bf q}_N\rightarrow\partial_z\tilde{\bf q}$
uniformly on $[0,T]\times\Omega$. 
The uniform convergence of $\partial_r{\bf q}_N$ can be shown in a similar way.
\qed

Before we can pass to the limit in the weak formulation of the approximate problems, 
there is one more useful observation that we need.
Namely, notice that although ${\bf q}$ are smooth functions both in the spatial variables and in time, 
the functions ${\bf q}_N$ are discontinuous at $n\Delta t$ because ${\tau_{\Delta t}\eta_N}$ is a step function
in time. As we shall see below, it will be useful to approximate each discontinuous function ${\bf q}_N$ in time
by a piece-wise constant function,  $\bar{\bf q}_N$, so that
$$
\bar{\bf q}_N(t,.)={\bf q}(n\Delta t-,.),\quad t\in [(n-1)\Delta t,n\Delta t),\ n=1,\dots,N,
$$
where ${\bf q}_N(n\Delta t-)$ is the limit from the left of ${\bf q}_N$ at  $n\Delta t$, $n=1,\dots,N$.
By combining Lemma \ref{testf} with the argument in the proof of Lemma \ref{konveta}, we get
$$
\bar{\bf q}_N\rightarrow \tilde{\bf q}\ {\rm  uniformly \ on}\   [0,T]\times\Omega.
$$

\subsection{Passing to the limit}

To get to the weak formulation of the coupled problem, take the test functions
$({\bf q}_N(t),\psi(t))$ (where $\mathbf q_N = \mathbf q \circ A_{{\tau_{\Delta t}\eta_N}}$, $\mathbf q \in {\cal X}_F$)
in equation (\ref{D1Prob1}) and integrate with respect to $t$ from $n\Delta t$ to $(n+1)\Delta t$.
Furthermore, take $\psi(t)\in {\cal X}_S$ as the test functions in (\ref{DProb3}), and again integrate
over the same time interval. 
Add the two equations together, and take the sum from $n=0,\dots,N-1$ to get the time integrals over $(0,T)$ as follows:
\begin{equation}
\begin{array}{c}
\displaystyle{\rho_f \int_0^T\int_{\Omega}(R+\tau_{\Delta t}\ \eta_N) \Big (\partial_t\tilde{\bf u}_N\cdot {\bf q}_N+
\frac 1 2(\tau_{\Delta t}{\bf u}_N-{\bf w}_N)\cdot\nabla^{{\tau_{\Delta t}\eta_N}}{\bf u}_N\cdot{\bf q}_N}\\ \\
\displaystyle{-\frac 1 2(\tau_{\Delta t}{\bf u}_N-{\bf w}_N)\cdot\nabla^{{\tau_{\Delta t}\eta_N}}{\bf q}_N\cdot{\bf u}_N\Big )+
\frac{\rho_f}{2}\int_0^T \int_{\Omega}{v^*_N}{\bf u}_N\cdot{\bf q}_N}\\ \\
\displaystyle{+\int_0^T\int_{\Omega}(R+\tau_{\Delta t}\eta_N)2\mu{\bf D}^{{\tau_{\Delta t}\eta_N}}({\bf u_N}):{\bf D}^{{\tau_{\Delta t}\eta_N}}({\bf q}_N)+
\rho_s h\int_0^T\int_0^L\partial_{t} \tilde{v}_N\psi}\\ \\
\displaystyle{+\int_0^T \left(a_S(\eta_N,\psi)+a_S'(v_N,\psi)\right)}\\
\\
\displaystyle{ =R\big (\int_0^TP_{in}^Ndt\int_0^1q_z(t,0,r)dr-\int_0^TP_{out}^Ndt\int_0^1q_z(t,L,r)dr\big ),}
\end{array}
\label{AproxEqNS}
\end{equation} 
with
\begin{equation}
\begin{array}{c}
\nabla^{\tau_{\Delta t}\eta}\cdot{\bf u}_N=0,\quad v_N=((u_r)_N)_{|\Gamma},
\\[0.3cm]
{\bf u}_N(0,.)={\bf u}_0,\ \eta(0,.)_N=\eta_0,\ v_N(0,.)=v_0.
\end{array}
\label{AproxEqNS}
\end{equation} 
Here $\tilde{\mathbf u}_N$ and $\tilde{v}_N$ are the piecewise linear functions defined in \eqref{tilde},
$\tau_{\Delta t}$ is the shift in time by $\Delta t$ to the left, defined in \eqref{shift},
$\nabla^{{\tau_{\Delta t}\eta_N}}$ is the transformed gradient via the ALE mapping $A_{{\tau_{\Delta t}\eta_N}}$,
defined in \eqref{nablaeta}, and $v_N^*$, $\mathbf u_N$, $v_N$ and $\eta_N$ are defined in \eqref{aproxNS}.

Using the convergence results obtained for the approximate functions in Section~\ref{sec:convergence}, and 
the convergence results just obtained for the test functions $\mathbf q_N$,
we can pass to the limit directly in all the terms except in the term that contains $\partial_t\tilde{\mathbf u}_N$.
To deal with this term we notice that, 
since ${\bf q}_N$ are smooth on sub-intervals $(j\Delta t,(j+1)\Delta t)$, we can use integration by parts on these sub-intervals to obtain:
$$
\int_0^T\int_{\Omega}(R+\tau_{\Delta t}\eta_N)\partial_t\tilde{\bf u}_N\cdot {\bf q}_N=
\sum_{j=0}^{N-1}\int_{j\Delta t}^{(j+1)\Delta t}\int_{\Omega}(R+\eta^j_N)\partial_t\tilde{\bf u}_N\cdot {\bf q}_N
$$
$$
=\sum_{j=0}^{N-1}\Big (-\int_{j\Delta t}^{(j+1)\Delta t}\int_{\Omega}(R+\tau_{\Delta t}\eta_N)\tilde{\bf u}_N\cdot\partial_t{\bf q}_N
$$
\begin{equation}\label{two_terms}
+\int_{\Omega}(R+\eta^{j+1}-\eta^{j+1}+\eta^j){\bf u}_N^{j+1}\cdot{\bf q}_N((j+1)\Delta t-)
-\int_{\Omega}(R+\eta^j){\bf u}_N^{j}\cdot{\bf q}_N(j\Delta t+)\Big ).
\end{equation}
Here, we have denoted by ${\bf q}_N((j+1)\Delta t-)$ and ${\bf q}_N(j\Delta t+)$
the limits from the left and right, respectively, of ${\bf q}_N$ at the appropriate points.

The integral involving $\partial_t{\bf q}_N$ can be simplified by recalling that
$\mathbf q_N = \mathbf q \circ A_{\eta_N}$, where $\eta_N$ are constant on each sub-interval $(j\Delta t,(j+1)\Delta t)$.
Thus, by the chain rule, we see that $\partial_t{\bf q}_N= \partial_t \mathbf q$ on $(j\Delta t,(j+1)\Delta t)$.
After summing over all $j = 0,...,N-1$ we obtain
$$
- \sum_{j=0}^{N-1}\int_{j\Delta t}^{(j+1)\Delta t}\int_{\Omega}(R+\tau_{\Delta t}\eta_N)\tilde{\bf u}_N\cdot\partial_t{\bf q}_N
=-\int_0^T\int_{\Omega}(R+\tau_{\Delta t}\eta_N)\tilde{\bf u}_N\cdot\partial_t{\bf q}.
$$
To deal with the last two terms in \eqref{two_terms} we calculate
$$
\sum_{j=0}^{N-1}\Big(\int_{\Omega}(R+\eta^{j+1}_N-\eta^{j+1}_N+\eta^j_N){\bf u}_N^{j+1}\cdot{\bf q}_N((j+1)\Delta t-)
-\int_{\Omega}(R+\eta^j_N){\bf u}_N^{j}\cdot{\bf q}_N(j\Delta t+)\Big)
$$
$$
=\sum_{j=0}^{N-1} \int_{\Omega}\left((R+\eta^{j+1}_N){\bf u}_N^{j+1}\cdot{\bf q}_N((j+1)\Delta t-)
-(\eta^{j+1}_N-\eta^j_N) {\bf u}_N^{j+1}\cdot{\bf q}_N((j+1)\Delta t-)\right)
$$
$$
- \int_{\Omega}(R+\eta_0) {\bf u}_0\cdot \mathbf q(0)
-\sum_{j=1}^{N-1} \int_{\Omega}(R+\eta^j_N){\bf u}_N^{j}\cdot{\bf q}_N(j\Delta t+)\Big)
$$
Now, we can write $(\eta^{j+1}-\eta^j)$ as $\displaystyle{v^{j+\frac 1 2}\Delta t}$, and rewrite the summation indexes
in the first term to obtain that above expression is equal to
$$
=\sum_{j=1}^{N} \int_{\Omega}(R+\eta^{j}_N){\bf u}_N^{j}\cdot{\bf q}_N(j\Delta t-)
-\int_0^T\int_{\Omega}v^*_N{\bf u}_N\cdot \bar{\bf q}_N
$$
$$
- \int_{\Omega}(R+\eta_0) {\bf u}_0\cdot \mathbf q(0)
-\sum_{j=1}^{N-1}\int_{\Omega}(R+\eta_N^j){\bf u}^j_N \cdot {\bf q}_N(j\Delta t+)
$$
Since the test functions have compact support in $[0,T)$, the value of the first term at $j = N$ is zero, and so 
we can combine the two sums to obtain
$$
=\sum_{j=1}^{N} \int_{\Omega}(R+\eta^{j}_N){\bf u}_N^{j}\cdot \left( {\bf q}_N(j\Delta t-)-{\bf q}_N(j\Delta t+)\right)
$$
$$
- \int_{\Omega}(R+\eta_0) {\bf u}_0\cdot \mathbf q(0)-\int_0^T\int_{\Omega}v^*_N{\bf u}_N\cdot \bar{\bf q}_N.
$$
Now we know how to pass to the limit in all the terms expect the first one. 
We continue to rewrite the first expression by using the Mean Value Theorem to obtain:
$$
{\bf q}_N(j\Delta t-,z,r)-{\bf q}_N(j\Delta t+,z,r)={\bf q}(j\Delta t,z,(R+\eta_N^j)r)-{\bf q}(j\Delta t,z,(R+\eta_N^{j+1})r)=
$$
$$
=\partial_r{\bf q}(j\Delta t,z,\zeta)r(\eta_N^j-\eta_N^{j+1})=-\Delta t\partial_r{\bf q}(j\Delta t,z,\zeta)v_N^{j+\frac 1 2}r.
$$ 
Therefore we have:
$$
\sum_{j=1}^{N-1}\int_{\Omega}(R+\eta_N^j){\bf u}^j_N\big ({\bf q}(j\Delta t-)-{\bf q}(j\Delta t+))=-
\int_0^{T-\Delta t}\int_{\Omega}(R+\eta_N){\bf u}_Nr\tau_{-\Delta t}v^*_N\partial_r\bar{\bf q}.
$$
We can now pass to the limit in this last term to obtain:
$$
\int_0^{T-\Delta t}\int_{\Omega}(R+\eta_N){\bf u}_Nr\tau_{-\Delta t}v^*_N\partial_r\bar{\bf q}\rightarrow 
\int_0^T\int_{\Omega}(R+\eta){\bf u}r\partial_t\eta\partial_r{\bf q}.
$$
Therefore, by noticing that $\partial_t\tilde{\bf q}=\partial_t {\bf q}+r\partial_t\eta\partial_r{\bf q}$
we have finally obtained
$$
\int_0^T\int_{\Omega}(R+\tau_{\Delta t}\eta_N)\partial_t\tilde{\bf u}_N\cdot {\bf q}_N\rightarrow -\int_0^T\int_{\Omega}(R+\eta){\bf u}\cdot\partial_t\tilde{\bf q}
-\int_0^T\int_{\Omega}\partial_t\eta{\bf u}\cdot\tilde{\bf q}
$$
$$
-\int_{\Omega}(R+\eta_0){\bf u}_0\cdot\tilde{\bf q}(0),
$$
where we recall that $\tilde{\mathbf q }= \mathbf q \circ A_\eta$.

Thus, we have shown that the limiting functions $\mathbf u$ and $\eta$ satisfy the weak form of Problem \ref{FSI} in the sense of Definition \ref{DefWSRef},
for all test functions $(\tilde{\mathbf q},\psi)$, which are dense in the test space ${\cal Q}_\eta$.
The following theorem holds: 

\begin{theorem}{\rm{\bf{(Main Theorem-The Viscoelastic Case)}}}
Let $\varrho_f$, $\varrho_s$, $\mu$, $h$, $C_i$, $D_i> 0$, $i=1,2,3$. Suppose that the initial data   $v_0\in L^2(0,L)$, ${\bf u}_0\in L^2(\Omega_{\eta_0})$,
and $\eta_0\in H^2_0(0,L)$ are
such that $(R+\eta_0(z))>0$, $z\in [0,L]$. Furthermore, let $P_{in}$, $P_{out}\in L^{2}_{loc}(0,\infty)$. 

Then there exist a $T>0$ and a weak solution of $(\mathbf u,\eta)$ of
Problem \ref{FSIref} (or equivalently Problem \ref{FSI}) on $(0,T)$ in the sense of Definition \ref{DefWSRef} (or equivalently Definition \ref{DefWS}),
which satisfy the following energy estimate:
\begin{equation}
E(t)+\int_0^tD(\tau)d\tau \leq E_0+C(\|P_{in}\|_{L^{2}(0,t)}^2+\|P_{out}\|_{L^{2}(0,t)}^2),\quad t\in [0,T],
\label{EE1}
\end{equation}
where $C$ depends only on the coefficients in the problem, $E_0$ is the kinetic energy of initial data,
and $E(t)$ and $D(t)$ are given by
\begin{eqnarray*}
E(t)&=&\frac{\rho_f}{2}\|{\bf u}\|^2_{L^2(\Omega_{\eta}(t))}+\frac{\rho_sh}{2}\|\partial_t\eta\|^2_{L^2(0,L)}\\
&+&\frac{1}{2}\big (C_0\|\eta\|^2_{L^2(0,L)}+C_1\|\partial_z\eta\|^2_{L^2(0,L)}+C_2\|\partial^2_z\eta\|^2_{L^2(0,L)}\big ),\\
D(t)&=&\mu\|{\bf D}({\bf u})\|^2_{L^2(\Omega_{\eta}(t)))}+
D_0\|\partial_t\eta\|^2_{L^2(0,L)}+D_1\|\partial_t\partial_z\eta\|^2_{L^2(0,L)}+D_2\|\partial_t\partial^2_z\eta\|^2_{L^2(0,L)}.
\end{eqnarray*}
\vskip 0.1in
 Furthermore, one of the following is true: either
\begin{enumerate}
\item $T=\infty$, or
\item $\displaystyle{\lim_{t\rightarrow T}\min_{z\in [0,L]}(R+\eta(z))=0}.$
\end{enumerate}
\label{ExistenceWS}
\end{theorem}
\proof
It only remains to prove the last assertion, which states that our result is either global in time, or, in case the walls of the cylinder 
touch each other, our existence result holds until the time of touching. 
This can be proved by using similar argument as in \cite{CDEM} p.~397-398. 
For the sake of completeness we present the arguments here. 

Let $(0,T_1)$, $T_1 > 0$, be the interval on which we have constructed 
our solution $({\bf u},\eta)$,
and let $m_1=\min_{(0,T_1)\times (0,L)} (R+\eta)$. From Lemma~\ref{eta_bound} we know that $m_1>0$. Furthermore, since
$\eta\in W^{1,\infty}(0,T;L^2(0,L))\cap L^{\infty}(0,T;H^2_0(0,L))$ and ${\bf u}\in L^{\infty}(0,T;L^2(\Omega))$,
we can take $T_1$ such that $\eta(T_1)\in H^2_0(0,L)$, $\partial_t\eta(T_1)\in L^2(0,L)$ and ${\bf u}(T_1)\in L^2(\Omega)$. 
We can now use the first part of Theorem \ref{ExistenceWS} to prolong the solution $({\bf u},\eta)$ to the interval $(0,T_2)$, $T_2>T_1$.
By iteration, we can continue the construction of our solution to the interval $(0,T_{k})$, $k\in\N$, where $(T_k)_{k\in\N}$ is an increasing sequence.
We set $m_k=\min_{(0,T_k)\times (0,L)} (R+\eta)>0$. 

Since $m_k > 0$ we can continue the construction further. Without loss of generality we could choose
a $T_{k+1} > T_k$  so that $m_{k+1}\geq\frac{m_{k}}{2}$. 
From (\ref{EE1}) and (\ref{compact_embedding}), by taking $\alpha=1/2$, we have that 
the displacement $\eta$ is H\"{o}lder continuous in time, namely,
$$
\|\eta\|_{C^{0,1/2}(0,T_{k+1};C[0,L])}\leq C(T_{k+1}).
$$
Therefore, the following estimate holds:
$$
R+\eta(T_{k+1},z)\geq R+\eta(T_{k},z)-C(T_{k+1})(T_{k+1}-T_{k})^{1/2}\geq m_{k}-C(T_{k+1})(T_{k+1}-T_{k})^{1/2}.
$$
For a $T_{k+1}$ chosen so that $m_{k+1}\geq\frac{m_{k}}{2}$
this estimate implies
\begin{equation}\label{prolong}
\displaystyle{T_{k+1}-T_{k}\geq \frac{m_{k}^2}{4C(T_{k+1})^2}},\; k\in\N.
\end{equation}

Now, let us take $T^*=\sup_{k\in\N}T_k$ and set $m^*=\min_{(0,T^*)\times (0,L)} (R+\eta)$. Obviously, $m_k\geq m^*$, $k\in\N$.
There are two possibilities. Either $m^* = 0$, or $m^* > 0$. 
If $m^*=0$, this means that 
$\displaystyle{\lim_{t\rightarrow T}\min_{z\in [0,L]}(R+\eta(z))=0}$, and the second statement in the theorem is proved.
If $m^* > 0$, we need to show that $T^* = \infty$.
To do that, notice that (\ref{EE1}) gives the form of the constant $C(T)$ which is a non-decreasing function of $T$.
Therefore, we have $C(T_k)\leq C(T^*)$, $\forall k\in\N$. Using this observation and that fact that $m_k\geq m^*$, $k\in\N$,
estimate (\ref{prolong}) implies
$$
T_{k+1}-T_{k}\geq \frac{(m^*)^2}{2C(T^*)^2},\; \forall k\in\N.
$$
Since this holds for all $k \in \N$, we have that $T^*=\infty$.
\qed

\section{The purely elastic case}\label{elastic_case}

We consider problem \eqref{NS}--\eqref{Coupling1b}  with the structural viscoelasticity constants $D_0=D_1=D_2=0$.
Thus, we study a FSI problem between the flow of a viscous, incompressible fluid, and the motion of a linearly elastic Koiter shell.
We show here that an analogue of Theorem \ref{ExistenceWS} holds, with a proof which is a small modification of the proof presented for the 
viscoelastic case. 

We begin by defining weak solutions.
% analogue to Definitions \ref{DefWS} and \ref{DefWSRef}.
First notice that the following formal energy inequality holds, which is an analogue of the energy inequality \eqref{EE}:
$$
\frac{d}{dt}E(t)+D_F(t)\leq C(P_{in}(t),P_{out}(t)).
$$
Here $D_F(t)=\mu\|{\bf D}({\bf u})\|^2_{L^2(\Omega_{\eta}(t)))}$ includes  fluid dissipation,
while $E(t)$ accounts for the kinetic energy of the fluid and of the structure, and the elastic energy of the Koiter shell.
Thus, $E(t)$ is the same as in the viscoelastic case, and is given by \eqref{energy}. 
%The entire dissipation mechanism in this problem comes from fluid viscosity.

To define a weak solution we introduce the solution spaces. 
The only difference with respect to the viscoelastic case is in the solution space for
the structure, which has lower regularity in time. 
Namely, we set
\begin{equation}\label{elastic_test}
{\cal W}_{SE}(0,T)=W^{1,\infty}(0,T;L^2(0,L))\cap L^{\infty}(0,T;{\cal V}_S),
\end{equation}
to replace the space ${\cal{W}}_S$, defined in \eqref{struc_test},
where ${\cal V}_S = H_0^2(0,L)$ as before.
The velocity space is the same as before, and is given by ${\cal{W}}_F$, defined in \eqref{vel_test}. 

The solution space for the coupled FSI problem is now given by
\begin{equation}\label{WE}
{\cal W}_E(0,T)=\{({\bf u},\eta)\in {\cal W}_F(0,T)\times{\cal W}_{SE}(0,T):{\bf u}(t,z,R+\eta(t,z))=\partial_t\eta(t,z){\bf e}_r\},
\end{equation}
which replaces the space \eqref{W} for the viscoelastic case. 

For the problem defined on a fixed, reference domain, the solution space
\begin{equation}
{\cal W^\eta}_E(0,T)=\{({\bf u},\eta)\in {\cal W}_F^\eta(0,T)\times{\cal W}_{SE}(0,T):{\bf u}(t,z,1)=\partial_t\eta(t,z){\bf e}_r\}
\end{equation}
replaces ${\cal W^\eta}$ defined in \eqref{W_eta}.

With these modifications, we have the following analogue definition of weak solutions:
\begin{definition}\label{DefWSElastic}
We say that $({\bf u},\eta)\in{\cal W}_E(0,T)$ is a weak solution of Problem~\ref{FSI} with $D_0=D_1=D_2=0$, if
for every $({\bf q},\psi)\in C^1_c([0,T);{\cal V}_F\times{\cal V}_{SE})$ such that 
${\bf q}(t,x,R+\eta(t,x))=\psi(t,x){\bf e}_r$, the  following equality holds:
\begin{equation}
\begin{array}{c}
\displaystyle{\rho_f\big (-\int_0^T\int_{\Omega_{\eta}(t)}{\bf u}\cdot\partial_t{\bf q}+\int_0^T b(t,{\bf u},{\bf u},{\bf q})\big )+
2\mu \int_0^T \int_{\Omega_{\eta}(t)}{\bf D}({\bf u}):{\bf D}({\bf q})}
\\ \\
\displaystyle{-\frac{\rho_f}{2}\int_0^T\int_0^L(\partial_t\eta)^2\psi
-\rho_sh\int_0^T\int_0^L \partial_t\eta\partial_t\psi+\int_0^T a_S(\eta,\psi)}
\\ \\
\displaystyle{=\int_0^T\langle F(t),{\bf q}\rangle_{\Gamma_{in/out}}+\rho_f\int_{\Omega_{\eta_0}}{\bf u}_0\cdot{\bf q}(0)+\rho_s h\int_0^Lv_0\psi(0).}
\end{array}
\label{VF_elastic}
\end{equation}
\end{definition}
For the problem defined on a fixed, reference domain, the definition of a weak FSI solution for the purely elastic structure case is
given by the following.
\begin{definition}\label{DefWSRefElastic}
We say that $({\bf u},\eta)\in{\cal W}^{\eta}_E(0,T)$ is a weak solution of Problem \ref{FSIref} with $D_0=D_1=D_2=0$,
if for every $({\bf q},\psi)\in C^1_c([0,T);{\cal V}^{\eta}_F\times{\cal V}_{SE})$ such that 
${\bf q}(t,z,1)=\psi(t,z){\bf e}_r$, the following equality holds:
\begin{equation}
\begin{array}{c}
\displaystyle{\rho_f\big (-\int_0^T\int_{\Omega}(R+\eta){\bf u}\cdot\partial_t{\bf q}+\int_0^T b^{\eta}({\bf u},{\bf u},{\bf q})\big )}\\ \\
+ \displaystyle{2\mu\int_0^T\int_{\Omega}(R+\eta){\bf D}^{\eta}({\bf u}):{\bf D}^{\eta}({\bf q})}
\displaystyle{-\frac{\rho_f}{2}\int_0^T\int_{\Omega}(\partial_t\eta){\bf u}\cdot{\bf q}}\\ \\
\displaystyle{-\rho_s h\int_0^T\int_0^L\partial_t\eta\partial_t\psi+\int_0^Ta_S(\eta,\psi)}
\\ \\
\displaystyle{=R\int_0^T\big (P_{in}(t)\int_0^1(q_z)_{|z=0}-P_{out}(t)\int_0^1(q_z)_{|z=L}\big )}
\\ \\
\displaystyle{+\rho_f\int_{\Omega_{\eta_0}}{\bf u}_0\cdot{\bf q}(0)+\rho_s h\int_0^Lv_0\psi(0).}
\end{array}
\label{VFRef_elastic}
\end{equation}
\end{definition}

\begin{theorem}{\rm{\bf{(Main Theorem - The Elastic Case)}}}
Let $\varrho_f$, $\varrho_s$, $\mu$, $h$, $C_i> 0$, $D_i=0$, $i=1,2,3$. Suppose that the initial data   $v_0\in L^2(0,L)$, ${\bf u}_0\in L^2(\Omega_{\eta_0})$,
and $\eta_0\in H^2_0(0,L)$ are
such that $(R+\eta_0(z))>0$, $z\in [0,L]$. Furthermore, let $P_{in}$, $P_{out}\in L^{2}_{loc}(0,\infty)$. 

Then, there exist a $T>0$ and a weak solution $(\mathbf u,\eta)$ on $(0,T)$ of
Problem \ref{FSIref} (or equivalently Problem \ref{FSI}) with $D_0=D_1=D_2=0$,  
in the sense of Definition \ref{DefWSRefElastic} (or, equivalently, Definition \ref{DefWSElastic}),
which satisfies the following energy estimate:
\begin{equation}
E(t)+\int_0^tD_F(\tau)d\tau \leq E_0+C(\|P_{in}\|_{L^{2}(0,t)}^2+\|P_{out}\|_{L^{2}(0,t)}^2),\quad t\in [0,T],
\label{EE1Elastic}
\end{equation}
where $C$ depends only on the coefficients in the problem, $E_0$ is the kinetic energy of initial data,
and $E(t)$ and $D_F(t)$ are given by
\begin{equation*}
\begin{array}{rcl}
E(t)&=&\displaystyle{\frac{\rho_f}{2}\|{\bf u}\|^2_{L^2(\Omega_{\eta}(t))}+\frac{\rho_sh}{2}\|\partial_t\eta\|^2_{L^2(0,L)}}\\
&+&\displaystyle{ \frac{1}{2}\big (C_0\|\eta\|^2_{L^2(0,L)}+C_1\|\partial_z\eta\|^2_{L^2(0,L)}+C_2\|\partial^2_z\eta\|^2_{L^2(0,L)}\big )},\\
D_F(t)&=&\displaystyle{\mu\|{\bf D}({\bf u})\|^2_{L^2(\Omega_{\eta}(t)))}.}
\end{array}
\end{equation*}
\vskip 0.1in
 Furthermore, one of the following is true: either
\begin{enumerate}
\item $T=\infty$, or
\item $\displaystyle{\lim_{t\rightarrow T}\min_{z\in [0,L]}(R+\eta(z))=0}.$
\end{enumerate}
\label{ExistenceWSE}
\end{theorem}
\proof
The main steps of the proof are analogous to the proof of Theorem \ref{ExistenceWS}. We summarize the main steps
and emphasize the main points where the proofs differ.

The problem is again {\bf split} into the structure elastodynamics problem, i.e., problem A1, and the fluid problem with the
Robin-type boundary condition, i.e., problem A2. Problem A1 is the same as in the viscoelastic case.
Problem A2 is different since only the structure inertia (the structure acceleration) 
is now coupled to the fluid motion via the boundary condition at the moving boundary, since $D_i = 0, i = 0,1,2$.

The {\bf semi-discretization} of Problem A1 is the same as before, including the analysis part in which the existence of a weak solution is proved.
The semi-discretization of Problem A2 differs by the fact that now the bi-linear functional associated with structural viscosity $a_S'={\bf 0}$.
This influences the solution spaces: a {\bf weak solution to Problem A2} now belongs to:
$$({\bf u}^{n+1},v^{n+1})\in {\cal V}_F^{\eta^n}\times L^2 (0,L),
{\rm (instead \ of\  } {\cal V}_F^{\eta^n}\times H^2_0 (0,L)).$$
The proof of the existence of a weak solution follows the same arguments based on the Lax-Milgram Lemma as before.

Furthermore, an analogous {\bf energy estimate} to \eqref{DEE} holds, with the  semi-discrete viscous dissipation $D_N^{n+1}$ replaced by
$$
\displaystyle{{(D_F)}_N^{n+1}=\Delta t\mu\int_{\Omega}(R+\eta^n)|D^{\eta^n}({\bf u}_N^{n+1})|^2.}
$$

Thus, up to this point, as before, we have a time-marching splitting scheme which defines an approximate solution to our FSI problem,
in which the structure is modeled as a linearly elastic Koiter shell. Furthermore, for each $\Delta t$, the approximate FSI solution
satisfies a discrete version of the energy estimate for the continuous problem.

The {\bf uniform energy bounds} for the approximate solutions follow
from  Lemma \ref{stabilnost} where ${(D_F)}_N^{n+1}$ replaces $D_N^{n+1}$.
The same method of proof provides the statements 1-4 of Lemma \ref{stabilnost},
with ${(D_F)}_N^{n+1}$ replaced by $D_N^{n+1}$
in the purely elastic case.

Therefore, we can proceed to get a {\bf weak convergence} result, i.e., an analogue of Lemma~\ref{weak_convergence}. 
Of course, we will not have the weak convergence 
$v_N \rightharpoonup  v\; {\rm weakly}\; {\rm in}\; L^2(0,T;H^2_0(0,L))$, however, this convergence property was only
used in passing to the limit in the bi-linear form $a_S'$, which no longer exists in the purely elastic case, and is therefore not needed in the proof.
Thus, a summary of the weak and weak* convergence results for the purely elastic case is given by the following:
\begin{equation}\label{weakconv_elastic}
\begin{array}{rcl}
\eta_N &\rightharpoonup & \eta \; {\rm weakly*}\; {\rm in}\; L^{\infty}(0,T;H^2_0(0,L)),
\\
v_N &\rightharpoonup & v\; {\rm weakly*}\; {\rm in}\; L^{\infty}(0,T;L^2(0,L)),
\\
v^*_N &\rightharpoonup & v^*\; {\rm weakly*}\; {\rm in}\; L^{\infty}(0,T;L^2(0,L)),
\\
{\bf u}_N &\rightharpoonup & {\bf u}\; {\rm weakly*}\; {\rm in}\; L^{\infty}(0,T;L^2(\Omega)),
\\
{\bf u}_N &\rightharpoonup & {\bf u}\; {\rm weakly}\; {\rm in}\; L^2(0,T;H^1(\Omega)).
\end{array}
\end{equation}

The next step is to obtain a {\bf  compactness result} which will be used to show that the limiting functions satisfy the weak formulation of
the FSI problem, i.e., we would like to obtain an analogue of Theorem \ref{u_v_convergence}. 
This is the place where the viscoelastic and purely elastic case differ the most.
Namely, in the proof of Theorem \ref{u_v_convergence} we used the fact that $(\partial_z v_N)_{n\in\N}$ is bounded in $L^2(0,T;L^2(0,L))$
(which followed from the uniform bound of the structural viscous dissipation terms).
This allowed us to have equicontinuity with respect to the spatial variable of the sequence $(v_N)_{n\in\N}$ in $L^2(0,L)$,
and use Corollary~\ref{Corollary}  to obtain the corresponding compactness result in $L^2(0,T;L^2(0,L))$
by considering only translations in time to show an integral-type ``equicontinuity'' result in time for the sequence $(v_N)_{n\in\N}$.
Since the viscoelastic structural terms are no longer present, the regularization of $v_N$ due to structural viscosity can no longer be used.

To get around this difficulty, we first recall that $v_N=({\bf u}_N)_{|\Gamma}\cdot{\bf e}_r$.
This follows from the definition of approximate functions \eqref{aproxNS} in which $v_N$ is defined
as a solution of Step A2 for a given $N\in\N$
(see the comment below \eqref{aproxNS}). 
Using the Trace Theorem we then have
$$
\|v_N\|_{L^2(0,T;H^{1/2}(0,L))}\leq C\|u_N\|_{L^2(0,T;H^1(\Omega))},\; N\in\N,
$$
where $C$ does not depend on $N$ (it depends only on $\Omega$). Therefore, we  have that 
$(v_N)_{n\in\N}$ is uniformly bounded in $L^2(0,T;H^{1/2}(0,L))$. Now, because of the
compactness of the embedding $H^{1/2}(0,L)\subset\subset L^2(0,L)$, the compactness in $L^2(0,L)$ of the sequence 
$(v_N)_{n\in\N}$ with respect to the spatial variable, follows immediately. Therefore, similarly as before, 
we only need to consider translations in time to prove integral ``equicontinuity'' in time, 
and use Corollary~\ref{Corollary} to obtain that $(v_N)_{n\in\N}$ is relatively compact in $L^2(0,T;L^2(0,L))$.
This part of the proof is identical to that of Theorem \ref{u_v_convergence}.

Showing relative compactness of $({\bf u}_N)_{n\in\N}$ is the same as in the viscoelastic case.

Therefore, at this point we have the following {\bf strong convergence} results:
\begin{equation}
\begin{array}{rcl}
\mathbf u_N &\to& u \  {\rm in}\  L^2(0,T;L^2(\Omega)),\\
v_ N &\to& v \  {\rm in}\  L^2(0,T;L^2(0,L)),\\
\tau_{\Delta t} \mathbf  u_N &\to& u \  {\rm in}\  L^2(0,T;L^2(\Omega)),\\
\tau_{\Delta t} v_ N &\to& v \  {\rm in}\  L^2(0,T;L^2(0,L)).
\end{array}
\end{equation}
To show strong convergence of $\eta_N$ we used only the information from the elastodynamics of the structure (Step A1),
and  the uniform bounds provided by Lemma~\ref{stabilnost} involving 
the kinetic energy of the fluid, the kinetic energy of the structure, and the structural elastic energy. See (\ref{derivativeeta}) and the proof below.
All this is the same as in the purely elastic case, and so the same strong convergence result holds for $\eta_N$ as before:
\begin{equation}
\eta_N \rightarrow \eta  \  {\rm in}\  L^{\infty}(0,T;H^s(0,L)), \ 0\leq s<2.
\end{equation}

The rest of the proof is the same, since the {\bf construction of the appropriate test functions} only uses convergence properties of
$\eta_N$, which are the same as in the viscoelastic case.

{\bf Passing to the limit} is performed in an analogous way as before.

%boundedness of $v_N^*$ (see (\ref{derivativeeta}). However,
%$v_N^*$ is defined from $v_N^{n+\frac 1 2}$ which are obtained in Step A1 which is the same as in viscoelastic case. 
%Therefore this part of the proof stays the same.
\qed

Therefore, slight modifications of the proof in the viscoelastic case, summarized above, lead to the proof of the existence of a weak solution
to a FSI problem between the flow of a viscous, incompressible fluid, and a linearly elastic Koiter shell model, as stated in Theorem~\ref{ExistenceWSE}.
This shows, among other things, the robustness of the method, presented in this manuscript.

\section{Appendix: The linearly viscoelastic Koiter shell model}
In this section we summarize the derivation of the linearly viscoelastic cylindrical Koiter shell model, presented in equation~\eqref{Koiter},
and list the expressions for the coefficients $C_i, D_i, i = 0,1,2,$ in terms of Youngs modulus of elasticity and Poisson ratio.
More details of the derivation can be found in \cite{MarSun,SunTam}.

Consider a clamped cylindrical shell with the reference radius of the middle surface equal to $R>0$, with the shell thickness $h>0$,
and cylinder length $L$. The basic assumptions under which the Koiter shell model holds are:
%\begin{itemize}
%\item 
(1) the shell is thin ($h/R << 1$);
%\item 
(2) the strains are small, although large deflections are admitted, and the strain energy per unit volume of the 
undeformed body is represented by a quadratic function of strain for an isotropic solid (Hooke's law);
%\item 
and (3) the state of stress is approximately plane.
%\end{itemize}

Denote by $\xib(z)=(\xi_z(z),\xi_r(z))$ 
the displacement of the middle surface at $z$, 
where $\xi_z(z)$ and $\xi_r(z)$ denote the longitudinal and radial 
component of displacement respectively.
Here, the axial symmetry of the problem has already been taken into
account assuming that the displacement in the $\theta$-direction is zero,
and that nothing in the problem depends on $\theta$. 
The change of metric and the change of
curvature tensors for a cylindrical shell are given respectively by \cite{PGC}
\begin{eqnarray*}
\gamb(\xib) = \left[\begin{array}{cc}
\xi_z' & 0\\
0 & R \xi_r
\end{array}\right],\quad
\varrhob(\xib) = \left[\begin{array}{cc}
-\xi_r'' & 0\\
0 &  \xi_r
\end{array}\right].
\end{eqnarray*}
Here $'$ denotes the derivative with respect to the longitudinal
variable $z$.
Introduce the following functional space
\begin{gather*}
V_c= H^1_0(0,L)\times H^2_0(0,L) =    
\left\{ (\xi_z,\xi_r) \in H^1(0, L)\times H^2(0, L):\right. \\
\left. \xi_z(0)=\xi_z(L)=\xi_r(0)=\xi_r(L)=0, \xi_r'(0)=\xi_r'(L) =0\right\}.
\end{gather*}
Then a weak formulation of
a linearly {\bf elastic} cylindrical Koiter shell is given by 
the following:
find $\etab=(\eta_z,\eta_r)\in V_c$ such that
\begin{equation*}\label{Koiter1}
\frac{h}{2}  \int_0^L \cA \gamb(\etab) \cdot
\gamb(\xib) Rdz + \frac{h^3}{24} \int_0^L \cA \varrhob(\etab) \cdot
\varrhob(\xib) Rdz = \int_0^L \fb \cdot\xib Rdz, \quad \xib \in V_c,
\end{equation*}
where $\cdot$ denotes the scalar product
$
A\cdot B := {\rm Tr} \left(A  B^T\right), \ A,B \in M_2(\ZR) \cong \ZR^4.
$
Here $\fb$ is the surface density of the force applied to the shell,
and ${\cal{A}}$ is the elasticity tensor given by \cite{PGC}:
\begin{eqnarray*}
\cA \Ebb &=& \frac{4\lambda \mu}{\lambda + 2 \mu} (\Abb^c \cdot
\Ebb) \Abb^c +4 \mu \Abb^c \Ebb \Abb^c,\quad  \Ebb \in
\Sym{(\ZR^2)}, \quad {\rm with}\\
\Abb_c &=&\left[ \begin{array}{cc} 1 & 0 \\ 0 &
R^2\end{array}\right], \qquad \Abb^c =\left[ \begin{array}{cc} 1 &
0 \\ 0 & \frac{1}{R^2}\end{array}\right],
\end{eqnarray*}
where $\lambda$ and $\mu$ are the Lam\'{e} constants.
Using the following relationships between the Lam\'{e} constants and the Young's modulus
of elasticity $E$ and Poisson ratio $\sigma$
$$
\frac{2\mu \lambda}{\lambda+2\mu} +2\mu =
4\mu\frac{\lambda+\mu}{\lambda+2\mu}= \frac{E}{1-\sigma^2},\qquad
\frac{2\mu \lambda}{\lambda+2\mu} =
4\mu\frac{\lambda+\mu}{\lambda+2\mu}\frac{1}{2}\frac{\lambda}{\lambda+\mu}
=  \frac{E}{1-\sigma^2} \sigma,
$$
the elasticity tensor $\cA$ reads
\begin{eqnarray*}
\cA \Ebb = \frac{2E \sigma}{1-\sigma^2} (\Abb^c \cdot \Ebb) \Abb^c
+ \frac{2E}{1+\sigma}  \Abb^c \Ebb \Abb^c,\quad  \Ebb \in
\Sym{(\ZR^2)}.
\end{eqnarray*}
Integration by parts in the weak formulation gives rise to the following static equilibrium equations
in differential form:
\begin{equation*}
\fbox{$
\begin{array}{rcl}
\displaystyle{-\frac{hE}{1-\sigma^2} \left(\eta_z'' + \sigma \frac{1}{R} \eta_r' \right)}&=& f_z, \\
\displaystyle{ \frac{hE}{R(1-\sigma^2)} 
\left(\sigma \eta_z' + \frac{\eta_r}{R}\right)  
+  \frac{h^3E}{12(1-\sigma^2)}
\left(  \eta_r'''' - 2 \sigma \frac{1}{R^2} \eta_r''  
+   \frac{1}{R^4} \eta_r \right)} &=& f_r.
\end{array}
$}
\end{equation*}
\centerline{\sc The Linearly Elastic Cylindrical Koiter Shell Model}
\vskip 0.05in
The terms multiplying $h/2$ account for the stored 
energy density due to stretching (membrane effects) and
the terms multiplying $h^3/12$ account for the stored energy density 
due to bending (flexural shell effects).

To include the {\bf viscoelastic} effects, we assume that $\etab$ is also a function of time.
Viscoelasticity will be modeled by the Kelvin-Voigt model 
in which the total stress is linearly proportional to strain and to the time-derivative of strain. 
For this purpose we define the {\bf viscosity} tensor $\mathbf{\mathcal{B}}$ by:
\begin{equation*}\label{vis}
\mathbf{\mathcal{B}} \textbf{E}=\frac{2 E_v \sigma_v}{1-\sigma_v^2}(\textbf{A}^c \cdot \textbf{E})\textbf{A}^c+ \frac{2 E_v}{1+\sigma_v} \textbf{A}^c \textbf{E} \textbf{A}^c, \quad \textbf{E} \in \textrm{Sym}(\mathbb{R}^2).
\end{equation*}
Here $E_v$ and $\sigma_v$ correspond to the viscous counterparts of the Young's modulus $E$ and Poisson's ratio $\sigma$. 
Then, for a linearly viscoelastic Koiter shell model we define the internal (stretching) force
\begin{equation*}
N := \frac{h}{2} \mathbf{\mathcal{A}} \boldsymbol\gamma(\boldsymbol\eta)+\frac{h}{2} \mathbf{\mathcal{B}} \boldsymbol\gamma(\dot{\boldsymbol\eta}),
\end{equation*}
and bending moment
\begin{equation*}
M := \frac{h^3}{24} \mathbf{\mathcal{A}} \boldsymbol\varrho(\boldsymbol\eta)+\frac{h^3}{24} \mathbf{\mathcal{B}} \boldsymbol\varrho(\dot{\boldsymbol\eta}).
\end{equation*}
The weak formulation of the {\bf dynamic equilibrium} problem for 
a linearly {\bf viscoelastic} Koiter shell is then given by the following: for each $t>0$ find $\boldsymbol\eta(\cdot,t) \in V_c$ such that $\forall \boldsymbol\xi \in V_c$
\begin{equation*}
\begin{array}{c}
\displaystyle{\frac{h}{2} \int_0^L (\mathbf{\mathcal{A}} \boldsymbol\gamma(\boldsymbol\eta) + \mathbf{\mathcal{B}} \boldsymbol\gamma(\dot{\boldsymbol\eta})) \cdot \boldsymbol\gamma(\boldsymbol\xi) R dz + \frac{h^3}{24} \int_0^L (\mathbf{\mathcal{A}} \boldsymbol\varrho(\boldsymbol\eta) + \mathbf{\mathcal{B}} \boldsymbol\varrho(\dot{\boldsymbol\eta})) \cdot \boldsymbol\varrho(\boldsymbol\xi) R dz}\\
\\
+ \displaystyle{\rho_s h \int_0^L \frac{\partial^2 \boldsymbol\eta}{\partial t^2} \cdot \boldsymbol\xi R dz= \int_0^L \mathbf{f} \cdot \boldsymbol\xi R dz,}
\end{array}
\end{equation*}
where $\rho_s$ denotes the volume shell density.

After integration by parts, the weak formulation above implies
the following dynamic equilibrium equations in differential form:
\begin{equation*}
\fbox{$
\begin{array}{r}
\displaystyle{\rho_s h \frac{\partial^2 \eta_z}{\partial t^2}-\tilde C_2 \frac{\partial \eta_r}{\partial z}-\tilde C_3 \frac{\partial^2 \eta_z}{\partial z^2}-\tilde D_2 \frac{\partial ^2 \eta_r}{\partial t \partial z} - \tilde D_3 \frac{\partial ^3 \eta_z}{\partial t \partial z^2} = f_z,}   \\
\displaystyle{\rho_s h \frac{\partial^2 \eta_r}{\partial t^2}+ \tilde C_0 \eta_r - \tilde C_1 \frac{\partial^2 \eta_r}{\partial z^2}  + \tilde C_2 \frac{\partial \eta_z}{\partial z} + \tilde C_4 \frac{\partial^4 \eta_r}{\partial z^4}}
\qquad \qquad \\
\displaystyle{+\tilde D_0 \frac{\partial \eta_r}{\partial t}- \tilde D_1 \frac{\partial^3 \eta_r}{\partial t \partial z^2} 
 +\tilde D_2 \frac{\partial ^2 \eta_z}{\partial t \partial z}+\tilde D_4 \frac{\partial ^5 \eta_r}{\partial t \partial z^4} = f_r,}
\end{array}
$}
\end{equation*}
\centerline{\sc { The Viscoelastic Cylindrical Koiter Shell  Model}}
where 
\begin{equation}
\small{
 \begin{array}{rlrlrl}
 \tilde C_0 &= \displaystyle{\frac{h E}{R^2(1-\sigma^2)}(1+\frac{h^2}{12 R^2}),}  & \tilde C_1 
 &\displaystyle{= \frac{h^3}{6} \frac{E \sigma}{R^2 (1-\sigma^2)}, } & \tilde C_2 
 &\displaystyle{=\frac{h}{R}\frac{E \sigma}{1-\sigma^2}, } \\
\tilde C_3&\displaystyle{=\frac{h E}{1-\sigma^2},} & \tilde C_4 &\displaystyle{=\frac{h^3}{12}\frac{E}{1-\sigma^2},}\\
\\
\tilde D_0 &= \displaystyle{\frac{h}{R^2} C_v(1+\frac{h^2}{12 R^2}),}  & \tilde D_1 &=\displaystyle{ \frac{h^3}{6} \frac{D_v}{R^2}, } &
 \tilde D_2 &= \displaystyle{\frac{h D_v}{R}, }\\
 \tilde D_3 &= h C_v, & \tilde D_4 &=\displaystyle{ \frac{h^3}{12}C_v,}
\label{coeff}
\end{array}}
\end{equation}
and
$$C_v:= \frac{E_v}{1-\sigma_v^2}, \quad D_v:= \frac{E_v \sigma_v}{1-\sigma_v^2}.$$

Mathematical justification of the Koiter shell model can be found in \cite{ciarlet1996asymptoticIII,Xiao,Li,Li2,Li3}.

In our problem we assume that  longitudinal displacement is negligible, i.e., $\eta_z = 0$, as is commonly used in 
modeling blood flow through human arteries. Thus, the resulting dynamics equilibrium equations in differential form are given by
\begin{equation*}
\fbox{$
\displaystyle{
\rho_s h \frac{\partial^2 \eta_r}{\partial t^2}+C_0 \eta_r -C_1 \frac{\partial^2 \eta_r}{\partial z^2}  + C_2 \frac{\partial^4 \eta_r}{\partial z^4}+D_0 \frac{\partial \eta_r}{\partial t}-D_1 \frac{\partial^3 \eta_r}{\partial t \partial z^2} 
+D_2 \frac{\partial ^5 \eta_r}{\partial t \partial z^4} = f_r,
}
$}
\end{equation*}
where 
\begin{equation*}
\small{
 \begin{array}{rlrlrl}
C_0 &= \displaystyle{\frac{h E}{R^2(1-\sigma^2)}(1+\frac{h^2}{12 R^2}),}  &  C_1 
 &\displaystyle{= \frac{h^3}{6} \frac{E \sigma}{R^2 (1-\sigma^2)}, } & C_2 &\displaystyle{=\frac{h^3}{12}\frac{E}{1-\sigma^2},}\\
\\
 D_0 &= \displaystyle{\frac{h}{R^2} C_v(1+\frac{h^2}{12 R^2}),}  & D_1 &=\displaystyle{ \frac{h^3}{6} \frac{D_v}{R^2}, } &
  D_2 &=\displaystyle{ \frac{h^3}{12}C_v,}
\end{array}}
\end{equation*}
This is exactly \eqref{Koiter}.

\vskip 0.2in
\noindent
{\bf {Acknowledgements.}}
The authors would like to thank Prof. Roland Glowinski for sharing his knowledge on numerical methods with us,
and for pointing out a couple of references. Special thanks are extended to Prof. Misha Perepelitza for his helpful discussions. 
The authors also acknowledge post-doctoral support for B. Muha provided by the Texas Higher Education Coordinating Board, 
Advanced Research Program (ARP) grant number 003652-0023-2009. \v{C}ani\'{c} also acknowledges partial 
research support  by the National Science Foundation by the grants  DMS-1109189, NSF DMS-0806941.

%\bibliography{c:/Users/Boris/Dropbox/myrefs}

\end{document}